%% file: main.tex
\documentclass[11pt]{extarticle}
\usepackage[utf8]{inputenc}

\usepackage[T1]{fontenc}
\usepackage{lmodern}
\usepackage{mathtools}
\usepackage[margin=0.5cm]{geometry}
\usepackage{lscape}
\usepackage{longtable}
\usepackage{array}
\usepackage{ragged2e}
\usepackage{booktabs}
\usepackage{xcolor}
\usepackage{pifont}
\usepackage{amssymb}
\usepackage{microtype}
\usepackage{setspace}
\newcolumntype{L}[1]{>{\RaggedRight\arraybackslash}p{#1}}
\newcommand{\NoMark}{\textcolor{red!70!black}{\ding{55}}}
\newcommand{\WarningMark}{\textcolor{orange!85!black}{\ensuremath{\blacktriangle}\kern-0.72em\raisebox{0.16ex}{\scriptsize\textcolor{white}{!}}\kern0.25em}}
\newcommand{\YesMark}{\textcolor{green!50!black}{\ensuremath{\checkmark}}}

\usepackage{amsthm}

\newtheorem{theorem}{Theorem}[section]
\newtheorem{proposition}[theorem]{Proposition}
\newtheorem{lemma}[theorem]{Lemma}
\newtheorem{corollary}[theorem]{Corollary}

\usepackage[linesnumbered,ruled,vlined]{algorithm2e}

\usepackage{graphicx}
\usepackage{multirow}
\usepackage{rotating}
\usepackage{placeins}

\usepackage[round,authoryear]{natbib}

\usepackage[nottoc]{tocbibind}
\usepackage{appendix}
\usepackage{geometry}
\usepackage{caption}
\usepackage{subcaption}
\usepackage{setspace}
\usepackage{float}
\usepackage{adjustbox}
\usepackage{longtable}
\usepackage{booktabs}
\usepackage{ltablex}
\keepXColumns

\usepackage{tabularx}
\usepackage{makecell}

\usepackage{algorithmicx}
\usepackage{algpseudocode}

\usepackage{microtype}
\usepackage{hyperref}

\usepackage{enumerate}   

\DeclareCaptionLabelFormat{subfiglabel}{Figure~\thefigure#2}
\newcommand{\best}[1]{\textbf{#1}}

\begin{document}
{\centering
{\fontsize{18}{22}\selectfont\textbf{Train Unit Scheduling with Unit Ordering under Platform-Feasible Operations}\par}

\vspace{16pt}

{\large
Yunjian Luo,\quad
Zhiyuan Lin\textsuperscript{*},\quad
Ronghui Liu\par
}

\vspace{8pt}

{\itshape Institute for Transport Studies, University of Leeds, Leeds, UK, LS2 9JT\par
}

\vspace{4pt}
{\small
tsylu@leeds.ac.uk; z.lin@leeds.ac.uk; r.liu@its.leeds.ac.uk\par
}
\vspace{6pt}
{\small * Corresponding author\par}
\vspace{12pt}
}

\textbf{Abstract}: In passenger railways where coupling and decoupling occur at platforms, a rolling-stock plan may be circulation-feasible but station-infeasible when the within-formation order of identified units causes blockage. We study the Train Unit Scheduling Problem under platform-feasible operations, where units cannot overtake or be resequenced without authorised shunting or resequencing. We formulate, to the best of our knowledge, the first single-stage unit-level integer linear programming model for this setting. It tracks identified units on a unit-indexed connection network and jointly determines the movements, coupling and decoupling decisions, and within-formation positions of train units, so every feasible integer solution provides a blockage-free schedule under the modelled restrictions. We further derive an exact fixed-assignment characterisation of orderability. Active coupling and decoupling requirements induce trip-wise precedence digraphs, while continuation arcs impose pairwise carry-over consistency. An assignment is orderable if and only if these digraphs admit a continuation-consistent family of topological orders. This yields a Train Unit Scheduler with Ordered Units (TUSOU), an exact branch-and-bound-and-cut train unit scheduler developed by us using ordering certification, lazy recovery of ordering constraints and activated cycle inequalities. Experiments on five real-world-derived TransPennine Express instances show that TUSOU produces certified blockage-free schedules, solves all instances to zero reported gap under solver tolerances, and outperforms direct full-model Gurobi baselines. Certification rejects 39 of 59 integer assignment-candidate encounters, showing that orderability should be embedded in optimisation rather than treated as post-processing.



\section{Introduction}\label{sec:introduction}
\input{paper/1_introduction}

\section{Problem description}\label{sec:Problem description}
\input{paper/2_Problem_description}

\section{Literature Review}\label{sec:Literature Review}
\input{paper/3_related_work}

\section{Integer linear programming model with unit ordering and blockage elimination}\label{sec:Enhanced formulation}
\input{paper/4_model}

\section{TUSOU: an exact solution framework}\label{sec:solver_framework}
\input{paper/5_alg}

\section{Implementation and computational study}\label{sec:Implementations}
\input{paper/6_Experiment}

\section{Conclusions}\label{sec:Conclusions}
\input{paper/7_conclusion}

\medskip
\newpage
\appendix
\appendixpage
\section{Main notation}
\label{app:notation}
\input{paper/12_notation}

\section{Omitted proofs}
\label{app:proof}
\input{paper/13_proof}

\section{Representative local ordering certificates}
\label{app:certificates}
\input{paper/16_certificates}
\newpage
\bibliography{main}
\end{document}

%% file: paper/1_introduction.tex
The train unit scheduling problem (TUSP) in passenger railways determines how train units are assigned and connected to operate a fixed timetable over a planning horizon, typically a single operating day \citep{lin2014two,lin2016branch}. Within the railway planning hierarchy, it lies between timetabling and crew scheduling and directly affects fleet utilisation, empty-running movements, operating cost, and operational robustness. In many UK and European passenger railway systems, capacity can be strengthened or weakened by coupling and decoupling train units between consecutive services. A daily rolling-stock plan must therefore decide not only which units cover which trips, but also how units are connected, where empty-running movements are used, and how coupled formations are reconfigured across the day.

A central operational difficulty is that platform coupling and decoupling are not determined solely by network connectivity. When unit coupling/decoupling is  handled directly at a platform, units  usually cannot overtake one another or exchange positions unless an explicit shunting or resequencing operation is possible. Hence the front–rear order of identified units within a formation becomes part of station feasibility. A rolling-stock plan may satisfy trip coverage, fleet capacity, timing, connection, and type-compatibility constraints, but may still be infeasible at the station level because the unit that must depart, split, or join first is trapped behind another unit. We refer to such a conflict as a platform blockage. A blockage can sometimes be repaired by shunting, but such recovery requires suitable infrastructure, operating resources, and sufficient turn-round time. When these are unavailable, the trapped unit cannot reach its planned departure position in time and the planned service may be delayed. In this paper, shunting and platform-side resequencing are therefore treated as explicit operational permissions rather than implicit repair options; when such permissions are absent, the schedule must be feasible under direct platform operations. We call this setting platform-feasible operations. It is typical in the UK railway system, where in many cases, coupling/decoupling is conducted directly at a platform without further shunting, especially during short turn-rounds.

The modelling implication is that platform blockage arises from the precise order of individual units within a train formation. A scheduling result that merely specifies which units operate which trips—without also determining each unit's within-formation position and the evolution of that order through continuations, couplings, and decouplings—cannot guarantee platform-feasible execution. The ordering information must be an explicit part of the decision space, not an afterthought. For example, if two units continue together, their relative order is preserved or reversed depending on the direction of travel. If units couple or decouple, the resulting order must be compatible with approach sides, departure sides, arrival sequence, or departure sequence. For multi-unit formations, these local ordering requirements may interact across several stations, and a feasible schedule exists only if all such requirements admit a globally consistent set of trip-wise formation orders.

Most existing approaches for the TUSP do not model this unit ordering information explicitly, or only deal with it in a sequential manner. Type-aggregated multicommodity-flow formulations can determine how many units of each type operate a trip and how type-level flows connect consecutive trips, but they do not decide which identified train unit occupies which position within a coupled formation. This creates a representation gap: the model can produce a network-feasible allocation, but the allocation remains ambiguous or infeasible as it does not contain information about unit  within-formation orders. Several approaches address this issue by dividing the original problem into two stages. \cite{lin2014two,lin2016branch} use a two-phase framework in which unit allocation is determined first and station feasibility is resolved afterwards, while \cite{lei2022resolution} strengthen this logic through iterative refinement between assignment and ordering resolution. These methods recognise that network-feasible rolling-stock assignments may fail at the station level, but the ordering decision is still not optimised simultaneously with identified unit assignment. If a selected assignment admits no consistent within-formation ordering, the method must repair the solution, impose an additional exclusion, or search for another assignment. One significant disadvantage is that global optimality often cannot be guaranteed by such sequential approaches.

To address these gaps, we formulate a single-stage integer linear programming model for the Train Unit Scheduling Problem under platform-feasible operations. The model lifts the standard connection-network representation from type-level flows to identified train-unit solutions  and jointly determines individual unit movements, trip coverage, feasible turn-round and empty-running connections, coupling and decoupling decisions, and within-formation order. Each feasible integer solution of the full model therefore specifies not only a circulation-feasible rolling-stock plan, but also an explicit unit-level formation order on every operated trip. Under the modelled platform restrictions, this solution is directly interpretable as an executable blockage-free schedule, rather than as an assignment-feasible plan that still requires post hoc ordering repair.

We further show that the ordering layer admits an exact fixed-assignment characterisation. Given an integral assignment \(\bar{x}\), the active coupling and decoupling instances induce trip-wise precedence digraphs on the units assigned to each trip, while active continuation arcs impose pairwise carry-over equalities across trips. Local acyclicity of the trip-wise precedence digraphs is necessary for orderability, but it is not sufficient by itself because continuation consistency must also hold. A fixed assignment is orderable if and only if the induced precedence digraphs admit a continuation-consistent family of topological orders. Equivalently, the projected orderable assignment set is \(X_{\mathrm{ord}}=\{\bar{x}\in X_{\mathrm{net}}:\Omega(\bar{x})\neq\emptyset\}\), where \(\Omega(\bar{x})\) is the set of blockage-free ordering witnesses for the selected unit movements. Thus ordering certification tests the existence of at least one admissible ordering witness, rather than the success of a particular constructive ordering rule.

Building on this characterisation, we develop the Train Unit Scheduler with Ordered Units (TUSOU), an exact certification-equipped branch-and-bound-and-cut framework. TUSOU searches over network-feasible unit-level assignments, but accepts an integer assignment only when ordering certification returns a blockage-free witness. If certification proves \(\Omega(\bar{x})=\emptyset\), the candidate assignment is rejected and missing original ordering instances from an infeasible active subsystem are recovered lazily. This separates exactness from uniform enforcement of all local-ordering constraints during LP bounding. We also introduce activated cycle inequalities which  provide auxiliary strengthening for the coupling/decoupling precedence layer. Exactness, however, is enforced by ordering certification and lazy recovery of original ordering instances, rather than by the cuts themselves.

The main contributions of this paper are as follows.

\begin{enumerate}[(i)]
\item We formulate, to the best of our knowledge, the first single-stage unit-level integer linear programming model for the Train Unit Scheduling Problem under platform-feasible operations. Unlike type-aggregated or two-phase formulations, the model jointly determines individual unit connectivity sequences, coupling and decoupling decisions, and unit within-formation positions. Hence each feasible integer solution directly yields an executable unit-level schedule under the modelled platform restrictions, rather than a network-feasible assignment that still requires post hoc ordering repair.

\item We provide an exact fixed-assignment characterisation of blockage-free orderability. For a fixed unit assignment, active coupling and decoupling decisions induce trip-wise precedence digraphs, while active continuation arcs impose pairwise carry-over consistency across trips. An assignment is orderable if and only if these local precedence digraphs admit a continuation-consistent family of topological orders. This result yields an exact ordering certification test for whether \(\Omega(\bar{x})\neq\emptyset\), i.e., whether the selected unit movements admit at least one blockage-free ordering witness.

\item We embed ordering certification in TUSOU, an exact certification-equipped branch-and-bound-and-cut framework. TUSOU accepts an integer assignment only if certification returns an ordering witness; otherwise, the assignment is rejected and missing original ordering instances from an infeasible active subsystem are recovered lazily. This allows exactness to be maintained without imposing the full local-ordering layer uniformly at every LP node and leads to the Shallow/Deep template design used in the computational study. Here, the Shallow template omits the unrecovered coupling/decoupling and continuation constraints during LP bounding, whereas the Deep template includes the full local-ordering layer to obtain ordering-enriched bounds. Under either template, every integer assignment candidate must pass ordering certification before it can be accepted.
\end{enumerate}

Computational experiments on five real-world-derived TransPennine Express corridor instances demonstrate the practical implications of these contributions. The Shallow TUSOU template solves all five instances to zero reported gap and is faster than the direct full-model Gurobi baselines on all five instances. Within TUSOU, the Shallow template is the most reliable configuration among the tested variants, and ordering certification provides explicit diagnostics on station infeasibility in the relaxed assignment space. Across the nine terminated TUSOU configurations with candidate diagnostics, ordering certification rejects 39 of 59 integer assignment-candidate encounters. In the Shallow configurations alone, it rejects 28 of 38 encounters. These diagnostics indicate that many network-feasible candidates admit no blockage-free ordering witness and that orderability must be embedded in optimisation rather than treated as post-processing.

The remainder of this paper is organised as follows. \S~\ref{sec:Problem description} describes the operational setting and formalises the problem. \S~\ref{sec:Literature Review} reviews related work. \S~\ref{sec:Enhanced formulation} presents the ordered-unit formulation and the fixed-assignment orderability results. \S~\ref{sec:solver_framework} describes TUSOU, including ordering certification, lazy recovery of original ordering instances, activated cycle inequalities, and the Shallow/Deep LP-template design. \S~\ref{sec:Implementations} reports the computational results. \S~\ref{sec:Conclusions} concludes and outlines directions for future research.

%% file: paper/2_Problem_description.tex
This section describes the operational setting and establishes why the Train Unit Scheduling Problem with platform-feasible operations presents a fundamental challenge that cannot be adequately addressed by existing two-phase approaches. We define the timetable network, the ordering conventions, the operational permissions that determine station-ordering restrictions, and the decision scope of the model.

\subsection{Timetable, units, and connection network}
Following \cite{cacchiani2010solving} and \cite{lin2014two,lin2016branch}, the timetable is represented by a directed acyclic graph. Each scheduled passenger trip is a node. A feasible transition of a train unit from trip $i$ to trip $j$ is represented by a directed turn-round arc $a = (i,j)$. Such an arc may represent a direct turn-round at the same station or an empty-running movement if the arrival station of $i$ differs from the departure station of $j$. Source and sink nodes represent sign-on and sign-off of train units. A source-to-sink path therefore corresponds to the daily sequence of trips and empty-running movements operated by one identified train unit.

\label{sec:formation_operations}
\begin{table}[!ht]
\centering
 \resizebox{0.7\linewidth}{!}{
 \begin{tabular}{|c|c|c|c|c|c|} 
 \hline
 Trip & Departure station & Arrival station & Departure time & Arrival time & Passenger demand \\ [0.5ex] 
 \hline
 1E06 & London Bridge & Uckfield & 8:30 & 10:40 & 250\\ 
 \hline
 2E32 & London Bridge & Uckfield & 9:20 & 11:20 & 300\\ 
 \hline
 2E11 & Uckfield & Oxted & 11:30 & 13:40 & 700\\ 
 \hline
 1E09 & East Croydon & London Bridge & 15:20 & 15:50 & 190\\ 
 \hline
 2E03 & Oxted & Uckfield & 14:00 & 15:40 & 310\\ 
 \hline
 \end{tabular}}
 \caption{A given timetable}
 \label{table:A given timetable}
\end{table}

\begin{figure}[ht]
    \centering
    \includegraphics[width=0.6\linewidth]{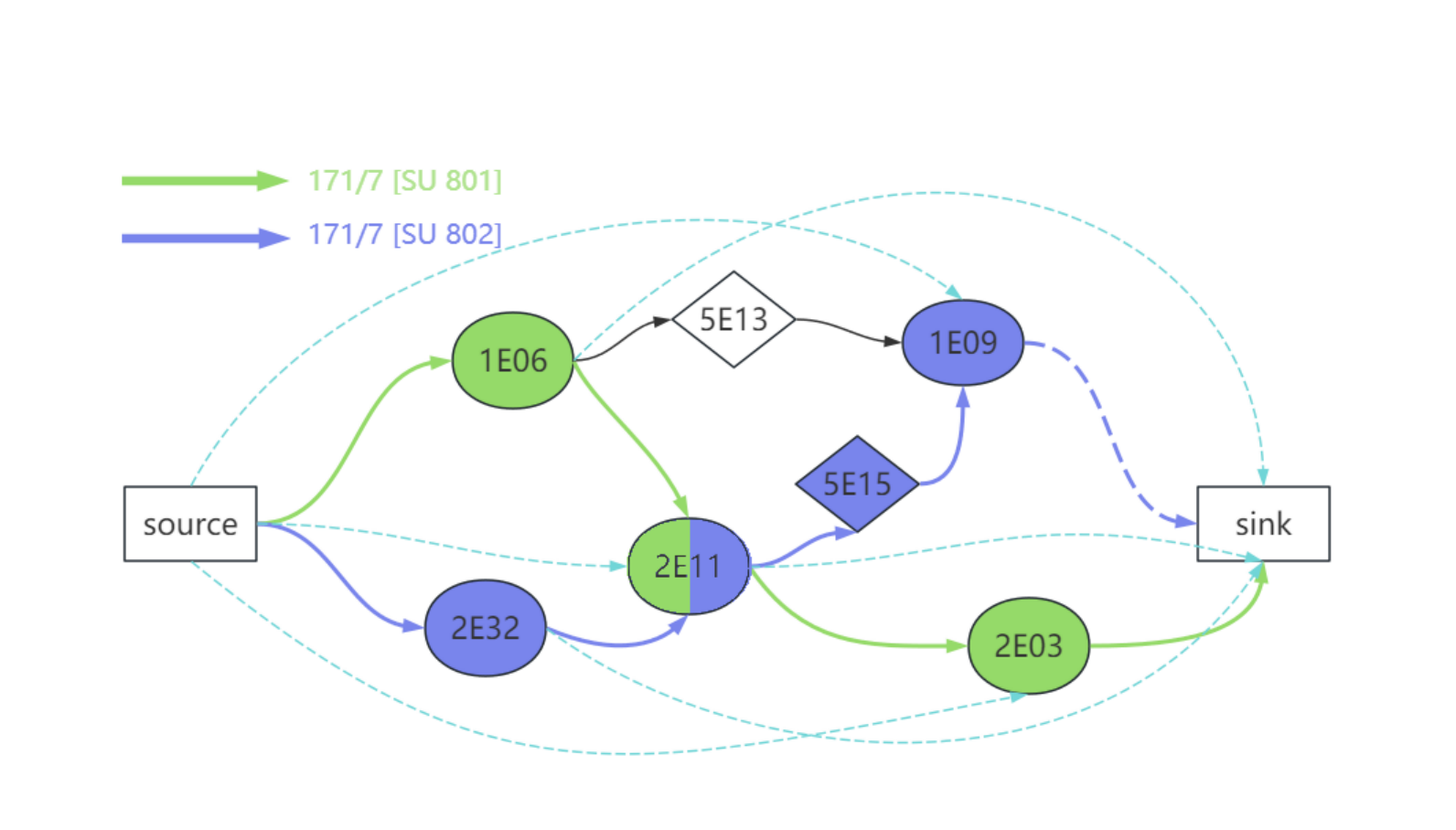}
    \caption{The corresponding directed acyclic graph based on the timetable in Table~\ref{table:A given timetable}}
    \label{fig:a solution using two units}
\end{figure}

Table~\ref{table:A given timetable} gives a small illustrative timetable, and Figure~\ref{fig:a solution using two units} shows the corresponding connection network. The same network structure underlies the optimisation model in \S~\ref{sec:Enhanced formulation}, but the model used in this paper is defined at the level of identified train units rather than type-aggregated flows as in prior work.

Let $T$ denote the set of train unit types and let $H_t$ be the set of available individual units of type $t\in T$. We write $H = \bigcup_{t \in T} H_t$ for the set of all identified units. Each unit $h \in H$ has a feasible arc set determined by type compatibility, timing, infrastructure restrictions, and operational permissions. The assignment of unit $h$ to a sequence of arcs defines its unit-level daily schedule.

\subsection{The platform blockage problem}
In passenger railway systems where coupling and decoupling operations are performed at station platforms, a rolling-stock plan that is circulation-feasible at the network level may nevertheless be station-infeasible. This occurs because network-level models abstract stations as single points and leave the within-formation order of identified train units undetermined \citep{lei2022resolution}. The resulting plan may appear valid at the level of unit flows, yet fail operationally when the physical order of units on the platform is considered.

\begin{table}[!ht]
\centering
 \resizebox{0.7\linewidth}{!}{
 \begin{tabular}{|c|c|c|c|c|} 
 \hline
 Trip & Departure station & Arrival station & Departure time & Arrival time  \\ [0.5ex] 
 \hline
 $i$ & B & A & 9:40 & 10:00 \\ 
 \hline
 $j_1$ & A & C & 10:20 & 10:50 \\ 
 \hline
 $j_2$ & A & B & 10:30 & 10:50 \\ 
 \hline
 \end{tabular}}
 \caption{A simple timetable}
 \label{table:A simple timetable}
\end{table}

\begin{figure}[htbp]
    \centering
    \begin{subfigure}[t]{0.4\textwidth}
        \centering
        \includegraphics[width=\linewidth]{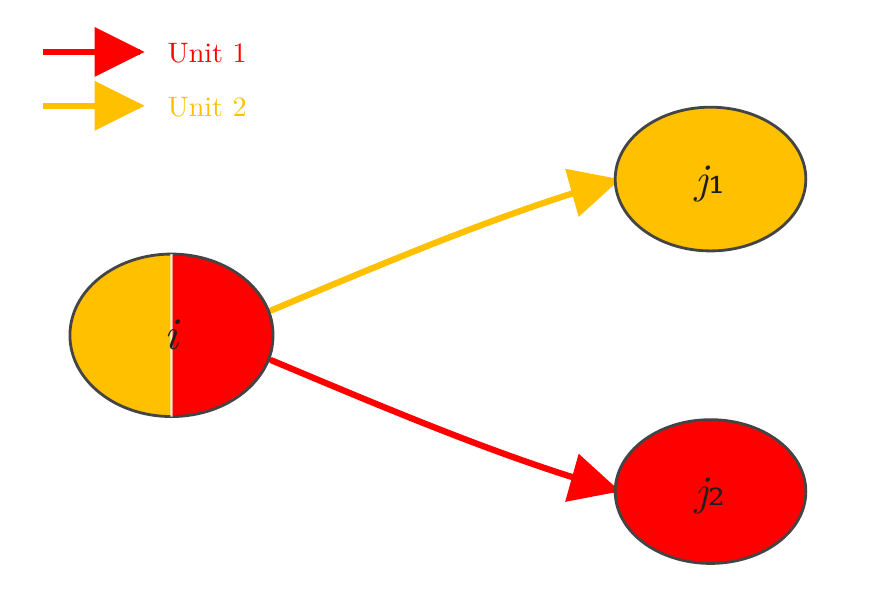}
        \caption{The corresponding directed acyclic graph based on the timetable in Table~\ref{table:A simple timetable}}
        \label{fig:figure1b blockage}
    \end{subfigure}\hfill
    \begin{subfigure}[t]{0.6\textwidth}
        \centering
        \includegraphics[width=\linewidth]{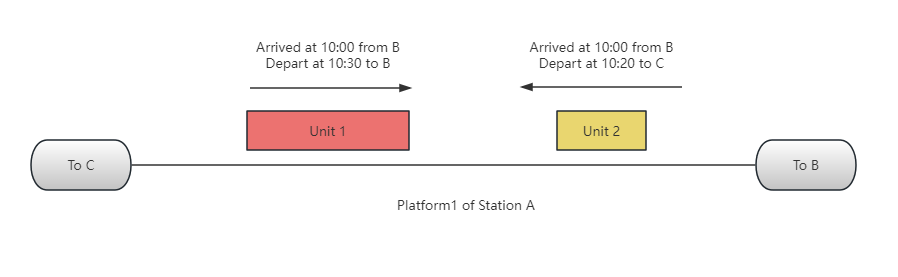}
        \caption{An example of station blockage}
        \label{fig:An example of the blockage on station level}
    \end{subfigure}
\end{figure}

Figure~\ref{fig:An example of the blockage on station level} illustrates the canonical blockage scenario as Table~\ref{table:A simple timetable}. Figure~\ref{fig:figure1b blockage} shows the corresponding connection network. The two train units initially operate together as a coupled formation and arrive at Station A from side B at 10:00. They are then scheduled to decouple and operate different subsequent trips: Unit 2 is due to depart towards side C at 10:20, whereas Unit 1 is due to depart towards side B at 10:30. However, Unit 1 is positioned between Unit 2 and the C-side exit. Consequently, Unit 2 cannot reach its required departure direction without passing Unit 1 on the platform. The connection plan is feasible at the timetable-network level, yet the identified-unit plan is station-infeasible because the unit that must leave first is trapped by the within-formation order. It is only one illustrative example of platform blockage. The ordering restrictions considered in this paper are intended to preclude all blockage scenarios that may arise from incompatible train-unit orderings at stations. This is not an edge case—it is a structural consequence of treating stations as points rather than as physical spaces with ordering constraints. 

Prior approaches have attempted to resolve this difficulty through two-phase methods: first solving a network-level integer multi-commodity flow problem that treats stations as points, then attempting to resolve station-level conflicts in a second phase \citep{lin2014two,lei2022resolution}. However, this separation has significant limitations. The station-level ordering constraints depend on the coupling and decoupling structure of the assignment, which in turn depends on which units are assigned to which trips and how they are connected. The two problems are coupled: one cannot determine whether an assignment is orderable without knowing the ordering, and one cannot choose an ordering without knowing the assignment. Post-processing an assignment to find a feasible ordering is a non-trivial feasibility problem in its own right. Indeed, as we shall demonstrate in \S~\ref{sec:solver_framework}, ordering certification rejects 39 of 59 integer assignment-candidate encounters, indicating that orderability must be embedded in optimisation rather than treated as an afterthought if global optimality is to be ensured.

\subsection{The structural sources of ordering difficulty}
The platform blockage problem arises from the interaction of three classes of ordering restrictions. Understanding their interaction is essential to appreciating why the problem is structurally complex.

\paragraph{Continuation carry-over} When two or more units continue together through the same feasible arc $a = (i,j)$, their relative physical order is carried over from the arrival formation to the departure formation, unless the arc is marked as allowing off-platform resequencing. For each trip $i$, let $d_i \in \{-1,1\}$ denote its travel direction relative to a fixed reference direction. Within a formation operating trip $i$ served by multiple units, unit positions are counted from the leading side in the direction of travel: position 1 denotes the front unit, and larger indices denote progressively rearward units. A unit $h_1$ is said to be \textit{ahead} of unit $h_2$ on trip $i$ if its position index is smaller. Under continuation, the pairwise position difference is preserved if $d_i = d_j$ and reversed if $d_i = -d_j$.

\paragraph{Coupling-induced precedence} Coupling events impose precedence requirements on the departing train formation. If incoming units or sub-formations approach from opposite directions, their relative order is determined by the accessible platform sides. If they approach from the same direction, their feasible order is determined by the arrival sequence, after applying any exogenous priority rule for tied arrival times.

\paragraph{Decoupling-induced precedence} Decoupling events impose precedence requirements on the arriving train formation. If the resulting units depart in opposite directions, each unit must already be positioned on the appropriate side. If they depart in the same direction, the unit that must depart first must be accessible from the relevant departure side, again after applying any exogenous priority rule for tied departure times.

These restrictions are enforced only where direct platform operations apply. If a turn-round, station event, or connection is marked as permitting shunting, off-platform resequencing, or another operation that can change the unit order, then the corresponding local ordering restriction is deactivated. Thus the model can represent mixed operating environments.

The difficulty lies in the \textit{global} consistency of these local restrictions. Coupling events create trip-wise precedence requirements: for a trip served by multiple units, the within-formation order must satisfy the precedence constraints induced by how those units arrived and how they will depart. Continuation arcs impose pairwise carry-over consistency across trips: the order of two units on a departure trip must match (or be the reverse of) their order on the arrival trip, depending on direction. The feasibility of an assignment therefore reduces to whether these trip-wise precedence digraphs admit a family of topological orders that are consistent across continuation arcs. This is not a trivial property—it is a combinatorial consistency problem that can fail even when each individual trip's precedence constraints are locally satisfiable. The characterisation of this property, which we term \textit{blockage-free orderability}, is given in \S~\ref{sec:Enhanced formulation} through the fixed-assignment certification test.

\subsection{Decision scope and input settings}
The Train Unit Scheduling Problem studied in this paper decides three interrelated elements:

\begin{enumerate}[(i)]
\item \textbf{Unit assignment.} Individual train units are assigned to passenger trips so that demand, type compatibility, and admissible formation requirements are satisfied.

\item \textbf{Connection selection.} Feasible turn-round and empty-running connections are selected, thereby constructing one source-to-sink path for each deployed unit.


\item \textbf{Within-formation ordering.} For every trip operated by a coupled formation under platform-feasible operations, the within-formation order of the identified units is determined explicitly.
\end{enumerate}

These three decisions are coupled. The choice of unit assignment determines which coupling and decoupling events occur, which in turn determines the precedence constraints that any ordering must satisfy. The ordering, conversely, can render an otherwise feasible assignment inexecutable. Separating these decisions —as two-phase approaches do— risks producing network-level plans that cannot be made station-feasible, requiring costly iteration or, worse, producing plans that appear optimal but are not operationally executable.

By contrast, the single-stage model presented in \S~\ref{sec:Enhanced formulation} jointly determines all three elements. It tracks individual units on a unit-indexed connection network and simultaneously determines unit movements, coupling and decoupling decisions, and within-formation positions. Each feasible integer solution yields an executable unit-level schedule under the modelled platform restrictions. The model is augmented with an exact ordering certification test that verifies whether a selected assignment admits at least one blockage-free ordering witness, and a branch-and-bound-and-cut framework that embeds this certification to recover missing ordering instances and strengthen the formulation with activated cycle inequalities.

The input data consist of both network-level and station-level information. Network-level inputs include the timetable, unit types, available identified units, feasible trip-unit compatibility, unit capacities, passenger demand, feasible families, connection arcs, travel times, turn-round times, and empty-running possibilities. Station-level inputs include trip directions, coupling and decoupling permissions, and the shunting or resequencing permissions that determine which ordering restrictions are active.

The objective reflects the main operating costs: deploying train units, operating passenger-service unit mileage, and using empty-running movements. The full integer linear programming formulation is given in \S~\ref{sec:Enhanced formulation}.

%% file: paper/3_related_work.tex
Train unit scheduling models differ substantially in their representation of rolling stock, the degree of integration between assignment and station-level feasibility, and whether within-formation ordering is handled explicitly. For the purposes of this paper, the key distinction is whether a model tracks identified individual units and can determine blockage-free platform operations within the optimisation model itself. From this perspective, the literature can be grouped into three main streams: (i) composition-based models, (ii) type-aggregated multicommodity-flow models (including two-phase frameworks that assign units first and resolve station-level ordering afterwards), and (iii) hypergraph-based unit-connectivity models.

\subsection{Composition-Based Models}
This stream, developed mainly for NS Reizigers, models rolling-stock scheduling at the level of train compositions rather than identified individual units. In this setting, the successor trip for each train is predetermined, which reflects the dense, high-frequency timetable structure of the Dutch railway network.

\cite{schrijver1993minimum} proposed a multicommodity flow model minimising the number of train units required on a single line, with at most two compatible unit types. \cite{alfieri2006efficient} extended this by introducing a transition graph that explicitly captures feasible composition changes and unit permutations, although the objective remained focused on fleet size and carriage-kilometres. \cite{fioole2006rolling} further enriched the model by allowing underway combining and splitting of trains, incorporating three evaluation criteria: seat-shortage kilometres, shunting operations, and carriage-kilometres. \cite{peeters2008circulation} developed a branch-and-price algorithm for the same problem class, demonstrating that the LP relaxation provides tight bounds and that optimal solutions can be obtained efficiently for real-world NS instances.

These composition-based models are highly relevant for rolling-stock circulation with predetermined successor relations. However, their decision objects are compositions and composition transitions rather than identified units operating under platform-order restrictions. Moreover, because each trip's successor is fixed in advance, these models do not need to determine how trips are linked. This reflects a different operational setting—timetables in the Netherlands are sufficiently regular and dense that pre-sequencing is practical—rather than a weakness of the models. Crucially, however, these models are not designed to determine whether a specific unit-level coupling or decoupling operation can be executed blockage-free under the platform-feasible operations studied here. This difference arises because coupling and decoupling operations are performed at depots in the Netherlands, rather than directly at platforms. This exemplifies how railway operational requirements vary across countries and influence the choice of modelling framework.

\subsection{Type-Aggregated Multicommodity-Flow Models}
A second major stream, developed mainly in Italy and the UK, formulates train unit assignment and scheduling as type-aggregated integer multicommodity-flow models. In this framework, timetabled trips are represented as nodes in a directed acyclic graph, feasible turn-round or empty-running connections as arcs, and train-unit types as integer commodities. Compared with composition-based models that rely on predefined successor relations, this modelling paradigm offers substantially greater flexibility in selecting turn-round connections and empty-running moves.

\cite{cacchiani2010solving} established this line with a path-based ILP formulation for the Train Unit Assignment Problem (TUAP), solved using an LP-based diving heuristic. Their computational results on real-world Italian regional instances demonstrated significant improvements over manual schedules (10–20\% reduction in fleet size).\cite{cacchiani2019effective} later proposed a peak-period heuristic for more efficient rolling stock planning. Subsequent work by \cite{cacchiani2013lagrangian} developed a Lagrangian heuristic that could find feasible solutions in much shorter time, making the approach suitable for real-time applications. 

Closely related UK studies by \cite{lin2014two,lin2016branch} incorporate practical operational constraints such as coupling/decoupling permissions, time allowances, type compatibility, and locations banned for coupling/decoupling. These move closer to the real train unit scheduling problem faced by UK operators. \cite{lin2016branch} developed a branch-and-price solver capable of handling instances of up to about 500 train trips, with customised branching strategies for train-family and banned-location constraints. \cite{lin2017train} further extended the model to handle bi-level capacity requirements, distinguishing between target (hard) and desirable (soft) capacity provisions, guided by historic capacity provisions and passenger count surveys.

However, the limitation of this stream is representational rather than algorithmic. The unit flow variables are type-aggregated or allocation-focused: they determine how many units of each type operate a trip or traverse a connection, but not which identified unit occupies which position in a coupled formation. This distinction matters because platform blockage is triggered by the ordering of specific units, not only by the number or type of units assigned to a trip. A type-level or allocation-level solution can therefore be feasible with respect to fleet circulation, timing, demand, and coupling compatibility while still failing to admit any within-formation order that executes all platform coupling/decoupling operations and continuations without overtaking or resequencing.

Within this second stream, the closest contributions to the present setting are those that explicitly address station-level ordering, albeit outside the main assignment model. \cite{lin2014two} propose a two-phase framework in which unit assignment is decided first in a network-level phase, and station-level feasibility—including unit blockage resolution and coupling-order assignment—is handled afterwards in a second phase. The second phase is modelled as a multidimensional matching problem that determines feasible shunting plans, parking berths, and parking methods at each station, resolving conflicts such as crossing linkages and platform overcapacity.

\cite{lei2022resolution} later introduce an iterative refinement between assignment and ordering resolution. Their approach detects unresolvable station-level conflicts in a Phase I solution, formulates them as valid cuts, and feeds them back to the core solver to eliminate those inoperable conflicts. Feasible coupling orders are assigned simultaneously with conflict detection. This adaptive approach allows incremental insertion of station-level constraints.

These methods recognise the operational importance of station feasibility, but they remain allocation-first: within-formation ordering is not part of the first-stage feasible region, and ordering incompatibilities may be detected only after an assignment has already been selected. If the selected assignment is non-orderable, the solution must be repaired, excluded, or followed by another assignment search. This fundamental limitation motivates the present paper, which embeds ordering certification inside the exact optimisation process. An integer assignment candidate is accepted only if ordering certification proves that the selected unit movements admit at least one blockage-free ordering witness.

\subsection{Hypergraph-Based Models}
A third stream, developed mainly for German passenger railways (Deutsche Bahn), uses hypergraph-based models. \cite{borndorfer2016integrated} introduced a hypergraph formulation that explicitly represents train-unit connectivity. The key innovation is that hyperarcs model the simultaneous movement of multiple train units, thereby representing train compositions and composition changes in a unified framework. This is not merely an alternative formulation: the hypergraph framework captures coupling and decoupling operations as hyperflows, where a hyperarc connects multiple arrival nodes to multiple departure nodes, representing the simultaneous movement of a coupled composition. A coarse-to-fine column-generation approach was developed to solve the resulting large-scale integer programs.

Subsequent works by \cite{grimm2017propagation,grimm2019cut,grimm2023assignment} developed propagation, cutting-plane, and column-generation variants of this framework, including handling of maintenance constraints through resource flows and cutting planes. The Deutsche Bahn implementation, reported by \cite{borndorfer2021deutsche}, demonstrated annual savings of €74 million and a reduction of 34,000 tonnes of CO$_2$ emissions, highlighting the practical impact of this modelling approach.

These models are closer to the present paper in modelling granularity than aggregate flow formulations, because they can represent ordered compositions and, in the DB setting, even the orientation of individual units. In the hypergraph model, trip nodes represent individual train units in specific positions within a composition, and hyperarcs represent the movement of coupled units. 

\cite{grimm2025comparison} provide a direct analytical and numerical comparison between the composition model \citep{fioole2006rolling} and the hypergraph model \citep{borndorfer2016integrated}. They prove that, in both the NS setting and a simplified DB setting, the LP relaxations of the two models provide identical bounds as long as certain assumptions are met. Their numerical experiments show that the composition model is more compact and faster for NS instances with predefined connections, while the hypergraph model outperforms the composition model for DB-light instances with flexible turnings and deadheading options. 

Nevertheless, the operational logic studied here remains different. In the German hypergraph literature, unit orientation is introduced primarily because infrastructure restrictions and coach-facing requirements make the facing of a unit operationally relevant (e.g., first-class coaches must be positioned correctly at platforms). By contrast, the central issue in the present paper is direction-aware platform blockage during coupling and decoupling when no resequencing is allowed on the platform. The hypergraph framework is structurally rich, but the published formulations do not explicitly implement the direction-aware blockage-handling ordering constraints that are central here.

\subsection{Other Related Contributions}
Several adjacent studies are relevant but address different planning settings or problem variants.

In the Chinese high-speed railway context, \cite{gao2020branch} formulated path-based and arc-segment-based ILP models for trip sequence planning, and \cite{gao2022weekly} later developed a weekly rolling-stock scheduling model. Owing to differences in operational procedures, however, these studies do not focus on within-formation platform-order feasibility.

In the context of railway shunting, \cite{freling2005shunting} proposed a two-stage approach that accounts for directional accessibility of shunting tracks and prevents overcapacity and crossing at sidings. \cite{kroon2008shunting} introduced an integrated formulation using virtual lanes to jointly address matching and assignment decisions, with computational experiments on NS stations. These papers are operationally related to our work, but they address depot or yard execution rather than embedding platform blockage avoidance directly into the day-ahead train unit scheduling model.

Studies on robustness have also been undertaken. \cite{cacchiani2012railway} proposed a robust two-stage optimisation model for handling large disruptions. \cite{cadarso2014improving} focused on improving robustness by penalising difficult shunting operations and expected delays. These contributions are complementary to the present work but do not address the specific problem of within-formation ordering under platform-feasible operations.

\cite{nielsen2012rolling} and \cite{wagenaar2017maintenance} extended the composition model to handle disruption management and maintenance appointments, respectively. \cite{wagenaar2017rolling} incorporated deadheading trips and adjusted passenger demand in rolling stock rescheduling, bringing the composition model closer to the flexibility of other modelling paradigms.

In the context of rolling stock rescheduling, \cite{lusby2017branch} proposed a path-based branch-and-price approach in which coupling and decoupling are represented as movements of train units between train services and station depots. Units of the same type remain interchangeable when no unit-specific restrictions apply, whereas maintenance-critical units, such as those subject to individual distance-to-maintenance limits, are modelled separately in the pricing problem. But it does not embed platform blockage avoidance or formation-order into the scheduling model, mainly because coupling/decoupling is conducted at depots in their problem settings.

\subsection{Solution Techniques Across Modelling Paradigms}
The choice of solution technique in rolling stock scheduling is closely linked to the modelling framework being used. Since the formulations developed in the literature differ substantially in structure and size, the corresponding computational approaches also vary. Certain formulations remain sufficiently compact to be handled directly by off-the-shelf commercial MILP solvers. This is the case, for instance, for the composition model of \cite{fioole2006rolling}, a number of its later developments \citep{nielsen2012rolling,wagenaar2017maintenance}, and the formulation introduced by Cadarso and Marín (2011).

By contrast, models with very large variable sets, often growing exponentially with the instance size, typically require decomposition-based methods. Column generation has therefore been adopted for formulations such as the hypergraph model of \cite{borndorfer2016integrated} and the path-based approaches of \cite{lusby2017branch}, \cite{cacchiani2010solving}, and \cite{lin2014two}. For some of these path-based formulations, additional work has focused on tightening the mathematical description. In particular, \cite{cacchiani2013lagrangian} and \cite{lin2016local} studied the use of local convex hulls to strengthen the underlying formulations. \cite{lin2020avoiding} further developed a branching strategy aimed at eliminating avoidable coupling and uncoupling decisions in multicommodity-flow rolling stock scheduling models.

Lagrangian relaxation has also been applied in this context. For example, \cite{cacchiani2013lagrangian} proposed a Lagrangian heuristic for the model of \cite{cacchiani2010solving}, in which passenger-capacity requirements and fleet-availability constraints are relaxed.

\subsection{Model Comparison}
To provide a structured overview of the modelling landscape, Table~\ref{tab:literature-review} compares the key papers across multiple dimensions: unit representation (type-aggregated vs. identified individual units), ordering handling (whether within-formation ordering is explicitly determined, ignored, or handled in post-processing), platform-blockage handling, passenger demand treatment, deadheading flexibility, solution method, and operational setting.

\begin{sidewaystable}[p]
\centering
\begingroup
\scriptsize
\setlength{\tabcolsep}{1.7pt}
\renewcommand{\arraystretch}{0.95}

\caption{Literature review comparison table}
\label{tab:literature-review}

\begin{tabular}{@{}
L{2.10cm}
L{2.60cm}
L{2.20cm}
L{2.80cm}
L{4.70cm}
L{2.00cm}
L{2.10cm}
L{3.00cm}
L{2.70cm}
@{}}

\toprule
\textbf{Reference}
& \textbf{Framework}
& \textbf{Unit Identification}
& \textbf{Within-Formation Ordering}
& \textbf{Platform Blockage / Ordering Feasibility}
& \textbf{Successor Relations}
& \textbf{Passenger Demand}
& \textbf{Solution Method}
& \textbf{Context}\\
\midrule

\cite{schrijver1993minimum}
& Multi-commodity flow
& Type-aggregated
& Ignored
& \NoMark{} Not handled
& Predefined
& Hard constraint
& Polynomial / Flow
& NS (Netherlands)\\

\addlinespace[1.5pt]

\cite{schrijver1993minimum} & Multi-commodity flow & Type-aggregated & Ignored & \NoMark{} Not handled & Predefined & Hard constraint & Polynomial / Flow & NS (Netherlands) \\
\addlinespace[1.5pt]
\cite{alfieri2006efficient} & Transition-graph flow & Type-aggregated & Explicit in compositions & \NoMark{} Not handled & Predefined & Hard constraint & MILP / Preprocessing & NS (Netherlands) \\
\addlinespace[1.5pt]
\cite{fioole2006rolling} & Composition model & Type-aggregated & Explicit in compositions & \NoMark{} Not handled & Predefined & Soft penalties (seat-shortage km) & MILP (CPLEX) & NS (Netherlands) \\
\addlinespace[1.5pt]
\cite{peeters2008circulation} & Composition model (branch-and-price) & Type-aggregated & Explicit in compositions & \NoMark{} Not handled & Predefined & Hard / Soft (weighted) & Branch-and-Price & NS (Netherlands) \\
\addlinespace[1.5pt]
\cite{cacchiani2010solving} & Path-based multi-commodity flow & Type-aggregated & Ignored & \NoMark{} Not handled & Flexible & Hard constraint & LP-based diving heuristic & Italy (regional) \\
\addlinespace[1.5pt]
\cite{cacchiani2013integer} & Lagrangian relaxation & Type-aggregated & Ignored & \NoMark{} Not handled & Flexible & Hard constraint & Lagrangian heuristic & Italy (regional) \\
\addlinespace[1.5pt]
\cite{cacchiani2012nominal} & Robust optimisation & Type-aggregated & Ignored & \NoMark{} Not handled & Flexible & Hard constraint & Two-stage stochastic / MILP & NS / Italy \\
\addlinespace[1.5pt]
\cite{lin2014two,lin2016branch} & Two-phase: Phase I = multi-commodity flow; Phase II = matching & Type-aggregated (Phase I); Individual (Phase II matching) & Ignored in Phase I; handled in Phase II & \WarningMark{} Post-processing (station-level conflicts resolved after assignment via matching model) & Flexible & Hard constraint & Branch-and-Price (Phase I) + MILP (Phase II) & UK (First ScotRail / Southern) \\
\addlinespace[1.5pt]
\cite{lin2017train} & Multi-commodity flow (bi-level capacity) & Type-aggregated & Ignored & \NoMark{} Not handled & Flexible & Bi-level (target + desirable) & Branch-and-Price & UK (ScotRail) \\
\addlinespace[1.5pt]
\cite{lei2022resolution} & Adaptive two-phase with iterative cuts & Type-aggregated (Phase I); Individual (cut generation) & Ignored in Phase I; resolved via iterative cuts & \WarningMark{} Iterative repair (conflicts detected after assignment; invalid solutions cut off and reassignment sought) & Flexible & Hard constraint & Adaptive branch-and-price + cut generation & UK (TransPennine Express) \\
\addlinespace[1.5pt]
\cite{borndorfer2016integrated} & Hypergraph (hyperflow) & Identified individual units & Explicit (positions and orientation) & \NoMark{} Not handled & Flexible & Hard constraint & Column generation (Coarse-to-Fine) & Germany (DB) \\
\addlinespace[1.5pt]
\cite{borndorfer2021deutsche} & Hypergraph (FEO implementation) & Identified individual units & Explicit (positions and orientation) & \NoMark{} Not handled & Flexible & Hard constraint & Coarse-to-Fine (C2F) + MILP & Germany (DB) \\
\bottomrule
\end{tabular}

\endgroup
\end{sidewaystable}

\begin{sidewaystable}[p]
\ContinuedFloat
\centering
\begingroup
\scriptsize
\setlength{\tabcolsep}{1.7pt}
\renewcommand{\arraystretch}{0.95}

\caption[]{Literature review comparison table (continued)}

\begin{tabular}{@{}
L{2.10cm}
L{2.60cm}
L{2.20cm}
L{2.80cm}
L{4.70cm}
L{2.00cm}
L{2.10cm}
L{3.00cm}
L{2.70cm}
@{}}

\toprule
\textbf{Reference}
& \textbf{Framework}
& \textbf{Unit Identification}
& \textbf{Within-Formation Ordering}
& \textbf{Platform Blockage / Ordering Feasibility}
& \textbf{Successor Relations}
& \textbf{Passenger Demand}
& \textbf{Solution Method}
& \textbf{Context}\\
\midrule

\cite{grimm2017propagation,grimm2019cut,grimm2023assignment}
& Hypergraph variants (propagation, cuts)
& Identified individual units
& Explicit (positions and orientation)
& \NoMark{} Not handled
& Flexible
& Hard constraint
& Column generation / Cutting planes
& Germany (DB)\\

\addlinespace[1.5pt]
\cite{grimm2017propagation,grimm2019cut,grimm2023assignment} & Hypergraph variants (propagation, cuts) & Identified individual units & Explicit (positions and orientation) & \NoMark{} Not handled & Flexible & Hard constraint & Column generation / Cutting planes & Germany (DB) \\
\addlinespace[1.5pt]
\cite{grimm2025comparison} & Composition vs Hypergraph comparison & Type-aggregated (CM) / Identified (HM) & Explicit in both models (CM: compositions; HM: positions) & \NoMark{} Not handled & Predefined (NS) / Flexible (DB-light) & Hard constraint & MILP (CPLEX) & NS / DB-light \\
\addlinespace[1.5pt]
\cite{gao2020branch,gao2022weekly} & Path-based / Arc-segment ILP & Type-aggregated & Ignored & \NoMark{} Not handled & Flexible (trip sequencing) & Hard constraint & Branch-and-Price / Heuristics & China (High-speed rail) \\
\addlinespace[1.5pt]
\cite{cadarso2011robust,cadarso2014improving} & MILP / Robust optimisation & Type-aggregated & Ignored & \NoMark{} Not handled & Flexible & Hard / Soft & MILP (CPLEX) & Spain (Rapid transit) \\
\addlinespace[1.5pt]
\cite{nielsen2012rolling}; Wagenaar et al. (2017) & Composition model (extended) & Type-aggregated & Explicit in compositions & \NoMark{} Not handled & Predefined & Soft penalties & MILP (CPLEX) / Rolling horizon & NS (Disruption / Maintenance) \\
\addlinespace[1.5pt]
\cite{lusby2017branch} & Path-based (branch-and-price) & Type-aggregated & Ignored & \NoMark{} Not handled & Flexible & Soft penalties & Branch-and-Price & Denmark (DSB S-tog) \\
\addlinespace[1.5pt]
Current Manuscript & Identified-unit ILP with ordering certification & Identified individual units & Explicit within-formation ordering with certification & \YesMark{} Handled inside optimisation with embedded certification (blockage-free orderability verified via trip-wise precedence digraphs and continuation consistency; branch-and-cut with lazy recovery) & Flexible & Hard constraint & Branch-and-bound-and-cut with ordering certification & UK (TransPennine Express) \\
\bottomrule
\end{tabular}

\endgroup
\end{sidewaystable}

\subsection{Summary and Research Gap}
The literature reveals a clear modelling gap. Composition-based models capture feasible train compositions and composition changes but rely on fixed successor relations and do not track identified units for platform-order feasibility. Type-aggregated network-flow models provide flexible timetable connectivity and have been extended with many real-world constraints, but they generally omit explicit within-formation ordering or defer it to a later stage. Two-phase and iterative approaches \citep{lin2014two,lei2022resolution} recognise the operational importance of station-level ordering but treat ordering as post-processing or repair, with the consequence that non-orderable assignments may be selected and must be rejected or repaired. Hypergraph-based unit-level models are structurally richer and can represent ordered compositions and unit orientation, but the published formulations do not explicitly implement the direction-aware blockage-handling ordering constraints required under platform-feasible operations.

What remains open is a train unit scheduling model that directly determines identified unit movements, turn-round reconfiguration, and executable within-formation order in a single optimisation framework under platform-feasible operations. The present paper addresses this gap by formulating a single-stage unit-level integer linear programming model that jointly determines unit assignments, coupling and decoupling decisions, and within-formation positions. The model is augmented with an exact ordering certification test and a branch-and-cut framework that embeds certification inside the optimisation process, ensuring that each feasible integer solution yields an executable unit-level schedule under the modelled platform restrictions. The present paper also provides the first exact fixed-assignment characterisation of blockage-free orderability, showing that an assignment is orderable if and only if trip-wise precedence digraphs admit a continuation-consistent family of topological orders.

%% file: paper/4_model.tex
This section develops the single-stage ordered-unit formulation and the key theoretical results that underpin the TUSOU framework. We first present the ILP model, which jointly determines identified unit movements, turn-round reconfiguration, and within-formation positions under platform-feasible operations. We then show that the ordering layer admits an exact structural characterisation: coupling/decoupling feasibility induces trip-wise precedence digraphs whose acyclicity is necessary and sufficient for the existence of a feasible local order, whereas formation continuation is governed separately by pairwise carry-over consistency. After briefly discussing the formulation choice for the ordering layer, we derive activated cycle inequalities as valid strengthening consequences of the acyclicity structure. These inequalities strengthen the LP relaxation. We further describe the exact algorithmic logic of TUSOU  in \S~\ref{sec:solver_framework}.

\subsection{Model description}\label{sec:Model description}
Let $\widetilde{N}$ be the set of trip nodes of the network DAG. We denote the artificial source and sink nodes by $0$ and $\infty$, respectively, and define the complete node set as $N:=\widetilde{N}\cup\{0,\infty\}$. A feasible connection possibility of a unit from trip $i$ to trip $j$ is represented by a directed arc $a=(i,j)\in\widetilde{A}$, called a turn-round arc. Arcs leaving the source are sign-on arcs, denoted by $A^{\mathrm{on}}:=\{(0,j):j\in\widetilde{N}\}$, and arcs entering the sink are sign-off arcs, denoted by $A^{\mathrm{off}}:=\{(j,\infty):j\in\widetilde{N}\}$. Hence, $A:=\widetilde{A}\cup A^{\mathrm{on}}\cup A^{\mathrm{off}}$. Let $T$ denote the set of train unit types. For each unit type $t\in T$, let $\tilde N^t\subseteq\tilde N$ denote the set of trips that can be operated by type $t$, and let $\tilde A^t\subseteq\tilde A$ denote the set of feasible turn-round arcs for type $t$. Let $E^t\subseteq\tilde A^t$ denote the subset of type-$t$ arcs that involve empty-running. For each individual unit $h\in H$, let $A^h\subseteq A$ denote its feasible arc set. As illustrated earlier in Figure~\ref{fig:a solution using two units}, a source-to-sink path corresponds to the daily trip sequence operated by one train unit. For turn-round arcs that involve off-platform shunting or depot returning, the arc itself remains available in the connection network, but the corresponding platform-ordering instances can be deactivated because overtaking or resequencing is permitted off-platform.

Coupling-compatible types are grouped into families, and $F$ denotes the set of all families. For each trip $j\in\tilde N$, let $F_j\subseteq F$ denote the set of families feasible for $j$, i.e., families whose unit-type combinations are compatible with the requirements of trip $j$. We assume throughout that the feasible families $F_j$ are mutually exclusive for each trip $j$, as in the underlying family-selection formulation. For each type $t\in T$, let $H_t=\{h_t(1),h_t(2),\ldots,h_t(b^t)\}$ denote the set of available individual units of type $t$, where $b^t$ is the available fleet size of type $t$. Let $H:=\bigcup_{t\in T}H_t$ denote the set of all individual train units. When no type distinction is required, we index units by $h\in H$; when the type matters, we write $h\in H_t, t\in T$.

We introduce the direction parameter $d_i\in\{-1,1\}$ for each trip $i\in\tilde N$. It represents the travel direction of trip $i$: $d_i=1$ corresponds to the chosen reference direction (e.g., left-to-right), and $d_i=-1$ to the opposite direction. When describing a station event, the ``arrival'' and ``departure'' directions refer to the directions of the corresponding incoming and outgoing trips at that station. We also use the departure time $\tau_i^{\text{dep}}$ and arrival time $\tau_i^{\text{arr}}$ of trip $i$ as input parameters.

We define the arc-assignment variable $x_a^h\in\{0,1\}$ for all $a\in A^h$ and $h\in H$, where $x_a^h=1$ if unit $h$ traverses arc $a$. The family-selection variable $y_j^f\in\{0,1\}$, for all $j\in\tilde N$ and $f\in F_j$, indicates whether trip $j$ is operated by family $f$.

For each node $i\in N$ and unit $h\in H$, let $\delta_+^h(i):=\{a\in A^h: a=(i,\cdot)\}$ and $\delta_-^h(i):=\{a\in A^h: a=(\cdot,i)\}$ denote the sets of feasible arcs of unit $h$ departing from and arriving at node $i$, respectively.

The objective minimizes a weighted sum of three terms:
\begin{equation}
\begin{aligned}
\min \quad
& W_{1} \sum_{t \in T} \sum_{h \in H_t} \sum_{a \in \delta_{+}^h(0)} x_a^h
+ W_{2} \sum_{t \in T} \sum_{j \in \tilde N^t} \sum_{h \in H_t} \sum_{a \in \delta_{+}^h(j)} \mu_{j}^{t} x_a^h \\
& + W_{3} \sum_{t \in T} \sum_{a \in E^{t}} \sum_{h \in H_t} x_a^h
\end{aligned}
\label{objective function}
\end{equation}
The first term aims to minimize the total number of units used. The second term aims to minimize the total passenger train unit-mileage, where $\mu_{j}^{t}$ represents the mileage of one unit of type $t$ operating train $j \in N$. The third term aims to minimize the total number of empty-running trains.

Table~\ref{table:2} in Appendix~\ref{app:notation} summarises the notation. In the main text, we define only the symbols needed for the formulation as they are introduced.

The formulation ensures the platform-feasible operating requirements defined in \S~\ref{sec:Problem description}. Ordering restrictions are instantiated only for station events and arcs for which no shunting or resequencing permission is available. The constraints of the ILP formulation are given below.

\begin{equation}
\begin{aligned}
& \sum_{a\in\delta_+^h(0)} x_a^h \le 1, \quad \forall h \in H
\end{aligned}
\label{constraint:unit used}
\end{equation}

\begin{equation}
\begin{aligned}
& \sum_{a \in \delta_{+}^h(j)} x_a^h = \sum_{a \in \delta_{-}^h(j)} x_a^h, \quad \forall h \in H_t,\ \forall t \in T,\ \forall j \in \tilde{N}
\end{aligned}
\label{constraint:flow balance}
\end{equation}

\begin{equation}
\begin{aligned}
& \sum_{h \in H_t} \sum_{a \in \delta_{+}^h(0)} x_a^h \leq b^t, \quad \forall t \in T
\end{aligned}
\label{constraint:unit upper bound}
\end{equation}

\begin{equation}
\begin{aligned}
& \sum_{t \in T} \sum_{h \in H_t} \sum_{a \in \delta_{+}^h(j)} \kappa^t x_a^h \geq r_j, \quad \forall j \in \tilde{N}
\end{aligned}
\label{constraint:passenger demand}
\end{equation}

\begin{equation}
\begin{aligned}
& \sum_{t \in f} \sum_{h \in H_t} \sum_{a \in \delta_{+}^h(j)} x_a^h \leq v^f y_j^f, \quad \forall j \in \tilde{N},\ \forall f \in F_j
\end{aligned}
\label{constraint:y 1}
\end{equation}

\begin{equation}
\begin{aligned}
& \sum_{f \in F_j} y_j^f = 1, \quad \forall j \in \tilde{N}
\end{aligned}
\label{constraint:y 3}
\end{equation}

\begin{equation}
\begin{aligned}
& d_j(d_{i_1}-d_{i_2})\bigl(\theta_j^{h_1}-\theta_j^{h_2}\bigr) \ge -M_{j,i_1,i_2}^{\mathrm{oppC}}\bigl(2-x_{a_1}^{h_1}-x_{a_2}^{h_2}\bigr),\\
& \forall j\in\tilde N,\ \forall h_1,h_2\in H,\ h_1\neq h_2,\ \forall a_1=(i_1,j)\in\delta_-^{h_1}(j)\cap \tilde A,\ \forall a_2=(i_2,j)\in\delta_-^{h_2}(j)\cap \tilde A
\end{aligned}
\label{constraint:couple opposite}
\end{equation}

\begin{equation}
\begin{aligned}
& d_j(d_{i_1}+d_{i_2})(\tau_{i_1}^{arr}-\tau_{i_2}^{arr})\bigl(\theta_j^{h_1}-\theta_j^{h_2}\bigr) \ge -M_{j,i_1,i_2}^{\mathrm{sameC}}\bigl(2-x_{a_1}^{h_1}-x_{a_2}^{h_2}\bigr),\\
& \forall j\in\tilde N,\ \forall h_1,h_2\in H,\ h_1\neq h_2,\ \forall a_1=(i_1,j)\in\delta_-^{h_1}(j)\cap \tilde A,\ \forall a_2=(i_2,j)\in\delta_-^{h_2}(j)\cap \tilde A
\end{aligned}
\label{constraint:couple same}
\end{equation}

\begin{equation}
\begin{aligned}
& d_i(d_{j_1}-d_{j_2})\bigl(\theta_i^{h_1}-\theta_i^{h_2}\bigr) \le M_{i,j_1,j_2}^{\mathrm{oppD}}\bigl(2-x_{a_1}^{h_1}-x_{a_2}^{h_2}\bigr),\\
& \forall i\in\tilde N,\ \forall h_1,h_2\in H,\ h_1\neq h_2,\ \forall a_1=(i,j_1)\in\delta_+^{h_1}(i)\cap \tilde A,\ \forall a_2=(i,j_2)\in\delta_+^{h_2}(i)\cap \tilde A
\end{aligned}
\label{constraint:decouple opposite}
\end{equation}

\begin{equation}
\begin{aligned}
& d_i(d_{j_1}+d_{j_2})(\tau_{j_1}^{dep}-\tau_{j_2}^{dep})\bigl(\theta_i^{h_1}-\theta_i^{h_2}\bigr) \ge -M_{i,j_1,j_2}^{\mathrm{sameD}}\bigl(2-x_{a_1}^{h_1}-x_{a_2}^{h_2}\bigr),\\
& \forall i\in\tilde N,\ \forall h_1,h_2\in H,\ h_1\neq h_2,\ \forall a_1=(i,j_1)\in\delta_+^{h_1}(i)\cap \tilde A,\ \forall a_2=(i,j_2)\in\delta_+^{h_2}(i)\cap \tilde A
\end{aligned}
\label{constraint:decouple same}
\end{equation}

\begin{equation}
\begin{aligned}
& (\theta_j^{h_1}-\theta_j^{h_2})-d_i d_j(\theta_i^{h_1}-\theta_i^{h_2})
\ge -M_{ij}^{\mathrm{cont}}\bigl(2-x_a^{h_1}-x_a^{h_2}\bigr),\\
& \forall h_1,h_2\in H,\ h_1\neq h_2,\ \forall a=(i,j)\in\tilde A\cap A^{h_1}\cap A^{h_2}
\end{aligned}
\label{constraint:continuation lower}
\end{equation}

\begin{equation}
\begin{aligned}
& (\theta_j^{h_1}-\theta_j^{h_2})-d_i d_j(\theta_i^{h_1}-\theta_i^{h_2})
\le M_{ij}^{\mathrm{cont}}\bigl(2-x_a^{h_1}-x_a^{h_2}\bigr),\\
& \forall h_1,h_2\in H,\ h_1\neq h_2,\ \forall a=(i,j)\in\tilde A\cap A^{h_1}\cap A^{h_2}
\end{aligned}
\label{constraint:continuation upper}
\end{equation}

\begin{equation}
\begin{aligned}
& s_i^h \le \theta_i^h \le q_i,\\
& \forall i\in\tilde N,\ \forall t\in T,\ \forall h\in H_t
\end{aligned}
\label{constraint:order2}
\end{equation}

\begin{equation}
\begin{aligned}
& \theta_i^h \le M_i s_i^h,\\
& \forall i\in\tilde N,\ \forall t\in T,\ \forall h\in H_t
\end{aligned}
\label{constraint:order2b}
\end{equation}

\begin{equation}
\begin{aligned}
& 1 - M_i(1-u_i^{h_1h_2}) - M_i\bigl(2-s_i^{h_1}-s_i^{h_2}\bigr) \le \theta_i^{h_1}-\theta_i^{h_2} \\
& \hspace{2.2cm} \le M_i u_i^{h_1h_2}-1 + M_i\bigl(2-s_i^{h_1}-s_i^{h_2}\bigr),\\
& \forall i\in\tilde N,\ \forall t_1,t_2\in T,\ \forall h_1\in H_{t_1},\ \forall h_2\in H_{t_2},\ h_1\neq h_2
\end{aligned}
\label{constraint:order3}
\end{equation}

\begin{equation}
\begin{aligned}
& \theta_i^h \in \Theta_i^h,\quad \forall i\in\tilde N,\ \forall t\in T,\ \forall h\in H_t
\end{aligned}
\label{constraint:order4}
\end{equation}

\begin{equation}
\begin{aligned}
& x_a^h \in \{0,1\}, \quad \forall a \in A^h,\ \forall h \in H_t,\ \forall t \in T
\end{aligned}
\label{x}
\end{equation}

\begin{equation}
\begin{aligned}
& y_j^f \in \{0,1\}, \quad \forall j \in \tilde{N},\ \forall f \in F_j
\end{aligned}
\label{y}
\end{equation}

\begin{equation}
\begin{aligned}
& u_i^{h_1h_2} \in \{0,1\},\\
& \forall i\in\tilde N,\ \forall t_1,t_2\in T,\ \forall h_1\in H_{t_1},\ \forall h_2\in H_{t_2},\ h_1\neq h_2
\end{aligned}
\label{constraint:order1}
\end{equation}

Constraints~(\ref{constraint:flow balance})--(\ref{constraint:y 3}) capture the standard flow, fleet, demand, family-selection used in train unit scheduling models. Constraints~(\ref{constraint:flow balance}) enforce the flow balance of each unit over the network, ensuring continuity of the unit schedule between consecutive trips. Constraints~(\ref{constraint:unit upper bound}) ensure that the number of units deployed for each type remains within the corresponding upper bound. Constraints~(\ref{constraint:passenger demand}) ensure that the capacity demand of each passenger train is met. Constraints~(\ref{constraint:y 1}) define the family--train indicator variables and restrict feasible coupling combinations according to the total number of units. Constraints~(\ref{constraint:y 3}) enforce that exactly one family is assigned to each train.

The new constraints are Constraint~(\ref{constraint:unit used}), the coupling/decoupling ordering Constraints~(\ref{constraint:couple opposite})--(\ref{constraint:decouple same}), the continuation consistency Constraints~(\ref{constraint:continuation lower})--(\ref{constraint:continuation upper}), and the order-linking Constraints~(\ref{constraint:order2})--(\ref{constraint:order4}). Constraint~(\ref{constraint:order1}) introduces the binary pairwise comparison indicators $u_i^{h_1h_2}$. 

Constraint~(\ref{constraint:unit used}) is the source-departure restriction: each individual train unit may leave the source node at most once. Together with the flow-balance Constraints~(\ref{constraint:flow balance}) on trip nodes, the binary arc variables, and the acyclic network structure, this implies that each unit can traverse at most one source-to-sink path and therefore can cover each trip at most once.

Two groups of ordering constraints are required. Constraints~(\ref{constraint:couple opposite})--(\ref{constraint:decouple same}) apply to coupling and decoupling station events and are instantiated on candidate pairs of unit-specific incoming or outgoing arcs. Constraints~(\ref{constraint:continuation lower})--(\ref{constraint:continuation upper}) apply to pairwise carry-over consistency on feasible arcs and are instantiated on candidate pairs of units that traverse the same arc. In both families, an instance is active only when the relevant unit-specific arc selections are made.

In Constraints~(\ref{constraint:couple opposite})--(\ref{constraint:couple same}), $i_1$ and $i_2$ denote the predecessor trips associated with the selected incoming arcs of the two units into trip $j$ (i.e., for a selected arc $a_r=(i_r,j)$ we use its tail trip $i_r$). Analogously, in Constraints~(\ref{constraint:decouple opposite})--(\ref{constraint:decouple same}), $j_1$ and $j_2$ denote the successor trips associated with the selected outgoing arcs after decoupling from trip $i$. In Constraints~(\ref{constraint:continuation lower})--(\ref{constraint:continuation upper}), both units traverse the same arc $a=(i,j)\in\tilde A$.

For notational convenience, we define the assignment indicator $s_i^h:=\sum_{a\in\delta_+^h(i)}x_a^h\in\{0,1\}, i\in\tilde N,$ which equals $1$ if unit $h$ operates trip $i$ and $0$ otherwise.

To support pairwise order comparisons, we introduce binary variables $u_i^{h_1h_2}\in\{0,1\}$ for each trip $i\in\tilde N$ and distinct units $h_1,h_2\in H$. When both units are active on trip $i$, $u_i^{h_1h_2}=1$ indicates that unit $h_1$ is behind unit $h_2$, i.e., $\theta_i^{h_1}>\theta_i^{h_2}$; the order-linking constraints handle the activation logic and the inactive cases.

In addition, for each trip $i\in\tilde N$ and unit $h\in H$, we introduce an integer order variable $\theta_i^h$ to record the position of unit $h$ within the formation operating trip $i$. Let $q_i:=\sum_{t\in T}\sum_{h\in H_t}s_i^h$ denote the number of units assigned to trip $i$, i.e., the formation size on trip $i$. We also define $\bar q_i:=\max_{f\in F_i}v^f, \Theta_i^h:=\{0,1,\ldots,\bar q_i\}$. We set $\theta_i^h=0$ if unit $h$ is not assigned to trip $i$. Otherwise, $\theta_i^h\in\{1,\ldots,q_i\}\subseteq \Theta_i^h$ records its position counted from the leading end in the travel direction of trip $i$: $\theta_i^h=1$ denotes the front unit, $\theta_i^h=2$ the second unit, ..., and $\theta_i^h=q_i$ the rear unit. The linking constraints enforce this interpretation. We set the trip-specific Big-$M$ parameter $M_i:=\bar q_i$ for the order-linking constraints.

Constraints~(\ref{constraint:order2}) and~(\ref{constraint:order2b}) link the order variables to the trip-assignment decisions so that $\theta_i^h=0$ for inactive units and $\theta_i^h\in\{1,\dots,q_i\}$ for active units. Constraint~(\ref{constraint:order3}) ensures that any two active units on the same trip receive distinct positive orders, while inactive pairs are relaxed by the Big-$M$ terms. Constraint~(\ref{constraint:order4}) restricts $\theta_i^h$ to the explicit domain $\Theta_i^h=\{0,1,\ldots,\bar q_i\}$. Setting $M_i:=\bar q_i$ in Constraints~(\ref{constraint:order2b})--(\ref{constraint:order3}) gives a tight trip-dependent bound for the global model. Together, these constraints guarantee that $\theta_i^h=0$ if and only if unit $h$ is not assigned to trip $i$, and otherwise $\theta_i^h$ is a distinct positive position in the coupled formation.

As for the big-M  for activated ordering instances, one valid choice is $M_{j,i_1,i_2}^{\mathrm{oppC}}:=2\bar q_j$ and $M_{i,j_1,j_2}^{\mathrm{oppD}}:=2\bar q_i$, while for the time-weighted cases one may take $M_{j,i_1,i_2}^{\mathrm{sameC}}:=|d_j(d_{i_1}+d_{i_2})(\tau_{i_1}^{arr}-\tau_{i_2}^{arr})|\bar q_j$
and $M_{i,j_1,j_2}^{\mathrm{sameD}}:=|d_i(d_{j_1}+d_{j_2})(\tau_{j_1}^{dep}-\tau_{j_2}^{dep})|\bar q_i$. For the continuation constraints, one valid choice is $M_{ij}^{\mathrm{cont}}:=2 \cdot \max\{\bar q_i,\bar q_j,1\}$.

\begin{figure}[htbp]
    \centering
    \begin{subfigure}[t]{0.32\textwidth}
        \centering
        \includegraphics[width=\linewidth]{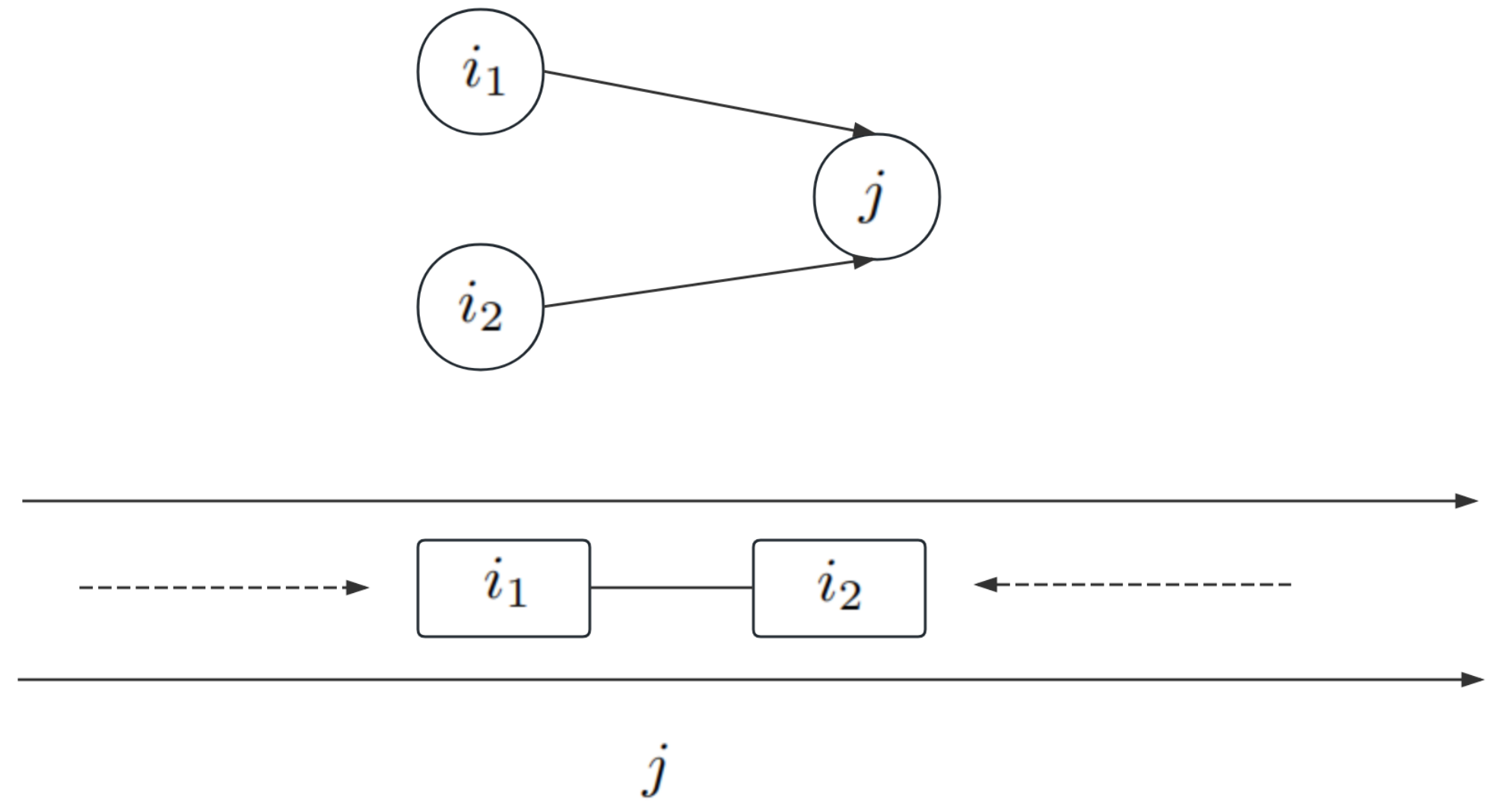}
        \caption{Opposite-direction coupling}
        \label{fig:Constraint (36)}
    \end{subfigure}\hfill
    \begin{subfigure}[t]{0.32\textwidth}
        \centering
        \includegraphics[width=\linewidth]{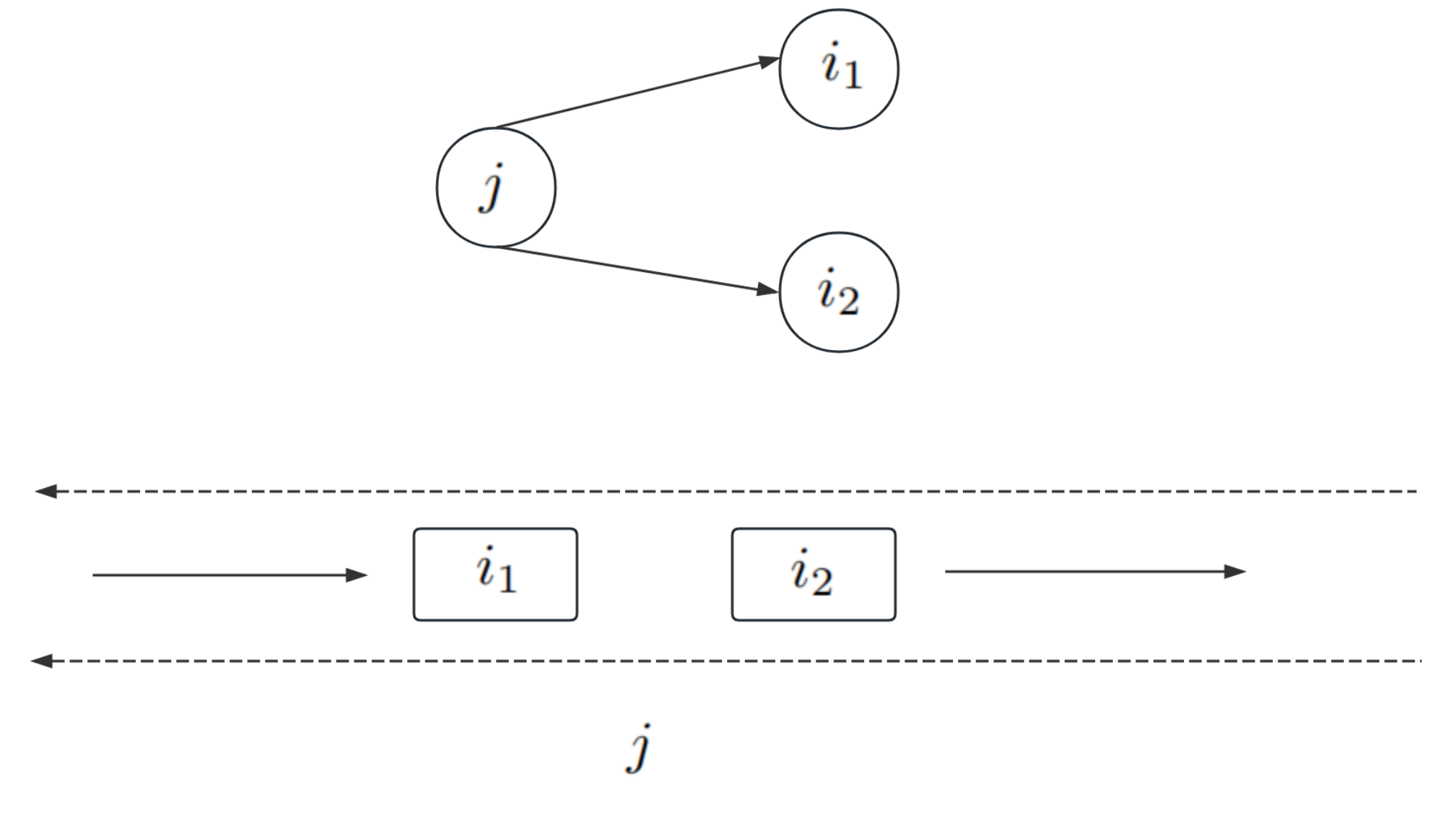}
        \caption{Same-direction decoupling}
        \label{fig:Constraint (39)}
    \end{subfigure}\hfill
    \begin{subfigure}[t]{0.32\textwidth}
        \centering
        \includegraphics[width=\linewidth]{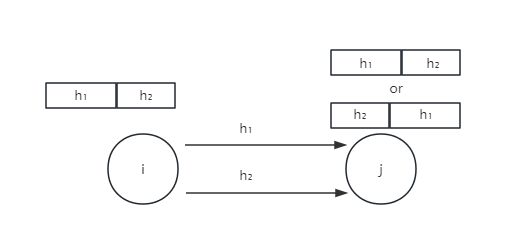}
        \caption{Formation  continuation}
        \label{fig:An example of no coupling/decoupling}
    \end{subfigure}
\end{figure}

We use an illustrative example to show how correct unit orders can be determined and blockages avoided. In Figure~\ref{fig:Constraint (36)} and ~\ref{fig:Constraint (39)}, dashed lines indicate trips already completed, whereas solid lines indicate trips to be operated next. Travelling from left to right is represented by $d_i=1$, while travelling from right to left is represented by $d_i=-1$. For any trip $i$, unit $h_1$ is ahead of unit $h_2$ if and only if $\theta_i^{h_1}<\theta_i^{h_2}$. The order variable $\theta_i^h$ is counted from the leading side in direction $d_i$, so smaller values correspond to positions closer to the leading side.

The ordering layer covers six local operating cases grouped into two families. The first family consists of four coupling/decoupling cases: opposite-direction coupling, opposite-direction decoupling, same-direction coupling, and same-direction decoupling. The second family consists of formation continuation on a feasible turn-round arc, with two directional realisations: same-direction continuation and opposite-direction continuation. In the active case, opposite-direction coupling and decoupling reduce to pure pairwise precedence requirements determined by platform accessibility from the relevant approach or departure sides; same-direction coupling and decoupling reduce to time-weighted precedence conditions determined respectively by arrival-time and departure-time order; and continuation reduces to preservation or sign reversal of the pairwise position difference according to whether the two trips have the same or opposite directions. In noncorresponding directional or timing configurations, the relevant coefficient becomes zero and the associated Big-$M$ term relaxes the constraint to a tautological form, so that the constraint is inactive unless the arc selections and operating case coincide. Detailed coefficient substitutions for all active and inactive cases are reported in Table~\ref{table:Situations in which constraints apply}.

Figures~\ref{fig:Constraint (36)} and~\ref{fig:Constraint (39)} illustrate the two representative coupling/decoupling mechanisms. In opposite-direction interactions, the active inequality fixes which unit must occupy the platform-accessible side of the formation so that a decoupling or coupling can be executed without crossing on the platform. In same-direction interactions, the active inequality links formation order to arrival-time or departure-time precedence, so that the unit that must move first is placed in the accessible position. The remaining active cases follow by sign reversal or by exchanging the roles of incoming and outgoing trips.

Formation continuation is handled separately by Constraints~(\ref{constraint:continuation lower})--(\ref{constraint:continuation upper}). Whenever two units $h_1$ and $h_2$ traverse the same feasible arc $a=(i,j)\in \tilde A$, the active case imposes $\theta_j^{h_1}-\theta_j^{h_2}=d_i d_j(\theta_i^{h_1}-\theta_i^{h_2})$.
Thus the pairwise position difference is preserved when the two trips have the same direction and reversed when they have opposite directions. Under the standing continuation assumption of this paper, this exactly captures feasible carry-over without overtaking, resequencing, or insertion/removal between the two units. Figure~\ref{fig:An example of no coupling/decoupling} illustrates this mechanism.

Overall, Constraints~(\ref{constraint:couple opposite})--(\ref{constraint:decouple same}) encode pairwise local feasibility for coupling and decoupling events, whereas Constraints~(\ref{constraint:continuation lower})--(\ref{constraint:continuation upper}) encode pairwise consistency on feasible arcs. Under the platform feasible assumptions, violations of either family correspond to blockage-causing local conflicts.

\begin{sidewaystable}[p]
\centering
\small
\renewcommand{\arraystretch}{1.2}
\setlength{\tabcolsep}{4pt}

\begin{adjustbox}{max totalsize={0.98\textheight}{0.98\textwidth},center}
\begin{tabular}{
@{}
>{\raggedright\arraybackslash}p{3.2cm}
>{\centering\arraybackslash}p{1.8cm}
>{\centering\arraybackslash}p{2.5cm}
>{\centering\arraybackslash}p{1.4cm}
>{\raggedright\arraybackslash}p{5.8cm}
>{\raggedright\arraybackslash}p{12.8cm}
@{}
}
\toprule
Constraint & Operation & \makecell{Relative direction} & \makecell{Coeff.\\product} & \makecell{Reduced\\inequality} & Scenario description \\
\midrule

Constraint~(\ref{constraint:couple opposite}) 
& coupling 
& same direction 
& 0 
& $0 \geq 0$ 
& Inactive in this scenario. \\

Constraint~(\ref{constraint:couple opposite}) 
& coupling 
& opposite direction 
& 2 
& $2(\theta_j^{h_1}-\theta_j^{h_2}) \geq 0$ 
& Enforces $\theta_j^{h_2}<\theta_j^{h_1}$, i.e., $h_2$ must be ahead of $h_1$ on trip $j$; otherwise the units would need to overtake on the platform. \\

Constraint~(\ref{constraint:couple opposite}) 
& coupling 
& opposite direction 
& $-2$ 
& $-2(\theta_j^{h_1}-\theta_j^{h_2}) \geq 0$ 
& Enforces $\theta_j^{h_1}<\theta_j^{h_2}$, i.e., $h_1$ must be ahead of $h_2$ on trip $j$; otherwise the units would need to overtake on the platform. \\

\addlinespace[2pt]

Constraint~(\ref{constraint:couple same}) 
& coupling 
& same direction 
& 2 
& $2(\theta_j^{h_1}-\theta_j^{h_2})(\tau^{arr}_{i_1}-\tau^{arr}_{i_2}) \geq 0$ 
& Enforces $\operatorname{sign}(\theta_j^{h_1}-\theta_j^{h_2})=\operatorname{sign}(\tau^{arr}_{i_1}-\tau^{arr}_{i_2})$, i.e., the earlier-arriving unit must have the smaller $\theta$ on trip $j$. \\

Constraint~(\ref{constraint:couple same}) 
& coupling 
& same direction 
& $-2$ 
& $-2(\theta_j^{h_1}-\theta_j^{h_2})(\tau^{arr}_{i_1}-\tau^{arr}_{i_2}) \geq 0$ 
& Enforces $\operatorname{sign}(\theta_j^{h_1}-\theta_j^{h_2})=-\operatorname{sign}(\tau^{arr}_{i_1}-\tau^{arr}_{i_2})$, i.e., the earlier-arriving unit must have the larger $\theta$ on trip $j$. \\

Constraint~(\ref{constraint:couple same}) 
& coupling 
& opposite direction 
& 0 
& $0 \geq 0$ 
& Inactive in this scenario. \\

\addlinespace[2pt]

Constraint~(\ref{constraint:decouple opposite}) 
& decoupling 
& same direction 
& 0 
& $0 \leq 0$ 
& Inactive in this scenario. \\

Constraint~(\ref{constraint:decouple opposite}) 
& decoupling 
& opposite direction 
& 2 
& $2(\theta_i^{h_1}-\theta_i^{h_2}) \leq 0$ 
& Enforces $\theta_i^{h_1}<\theta_i^{h_2}$, i.e., $h_1$ must be ahead of $h_2$ in the arriving formation on trip $i$ to allow the unit assigned to $j_1$ to depart without overtaking. \\

Constraint~(\ref{constraint:decouple opposite}) 
& decoupling 
& opposite direction 
& $-2$ 
& $-2(\theta_i^{h_1}-\theta_i^{h_2}) \leq 0$ 
& Enforces $\theta_i^{h_2}<\theta_i^{h_1}$, i.e., $h_2$ must be ahead of $h_1$ in the arriving formation on trip $i$ to allow the unit assigned to $j_2$ to depart without overtaking. \\

\addlinespace[2pt]

Constraint~(\ref{constraint:decouple same}) 
& decoupling 
& same direction 
& 2 
& $2(\theta_i^{h_1}-\theta_i^{h_2})(\tau^{dep}_{j_1}-\tau^{dep}_{j_2}) \geq 0$ 
& Enforces $\operatorname{sign}(\theta_i^{h_1}-\theta_i^{h_2})=\operatorname{sign}(\tau^{dep}_{j_1}-\tau^{dep}_{j_2})$, i.e., the earlier-departing unit must have the smaller $\theta$ on trip $i$. \\

Constraint~(\ref{constraint:decouple same}) 
& decoupling 
& same direction 
& $-2$ 
& $-2(\theta_i^{h_1}-\theta_i^{h_2})(\tau^{dep}_{j_1}-\tau^{dep}_{j_2}) \geq 0$ 
& Enforces $\operatorname{sign}(\theta_i^{h_1}-\theta_i^{h_2})=-\operatorname{sign}(\tau^{dep}_{j_1}-\tau^{dep}_{j_2})$, i.e., the earlier-departing unit must have the larger $\theta$ on trip $i$. \\

Constraint~(\ref{constraint:decouple same}) 
& decoupling 
& opposite direction 
& 0 
& $0 \geq 0$ 
& Inactive in this scenario. \\

\bottomrule
\end{tabular}
\end{adjustbox}

\caption{Detailed coefficient substitutions for all active and inactive coupling/decoupling scenarios generated by Constraints~(\ref{constraint:couple opposite})--(\ref{constraint:decouple same}). Throughout, unit $h_1$ is ahead of unit $h_2$ on trip $k$ if and only if $\theta_k^{h_1}<\theta_k^{h_2}$.}
\label{table:Situations in which constraints apply}
\end{sidewaystable}

\subsection{Characterisation of Ordering Feasibility}\label{sec:Ordering Feasibility}
\subsubsection{Local two-unit ordering feasibility under platform-feasible operations}\label{sec:equiv-station-feasibility}

This subsubsection proves the pairwise correctness of the local ordering layer under platform-feasible operations. Constraints~(\ref{constraint:couple opposite})--(\ref{constraint:decouple same}) encode the four coupling/decoupling station-event scenarios, while Constraints~(\ref{constraint:continuation lower})--(\ref{constraint:continuation upper}) encode pairwise formation continuation. We focus first on two units because each active local instance concerns a specific unit pair: coupling/decoupling induces a local precedence requirement on a single trip, whereas continuation imposes a cross-trip consistency requirement between two adjacent trips.

Consider either (i) a station event involving two units \(h_1\) and \(h_2\) that couple to form a departing trip \(j\) or decouple from an arriving trip \(i\), or (ii) two units traversing the same turn-round arc \(a=(i,j)\in \tilde A\). Let the platform order be the linear left-to-right order of the two units along the platform. We say that the local operation is feasible under platform-feasible operations if it can be realised using only allowed platform moves---arrival, departure, coupling or decoupling at the point of contact, and continuation without overtaking---and without requiring \(h_1\) and \(h_2\) to overtake or exchange positions. For case (ii), feasibility is understood together with the standing continuation assumption that no unit may be inserted between or removed from between the two units during the turn-round movement.

For same-direction coupling and decoupling events, all relevant arrival and departure times are assumed to be distinct; if ties occur, they are resolved by an exogenous operational priority order before the ordering constraints are generated. The five detailed pairwise results are stated and proved in Appendix~\ref{app:proof}: Lemmas~\ref{lem:equiv-coupling-opposite}--\ref{lem:equiv-decoupling-same} cover the four coupling/decoupling cases, and Lemma~\ref{lem:equiv-continuation} covers pairwise formation continuation. We then have the following result:

\begin{theorem}[Two-unit local feasibility versus ordering constraints]
\label{thm:station-feasibility-equiv}
Under the platform-feasible operation assumption of \S~\ref{sec:Model description}, the active local ordering constraints are necessary and sufficient for two-unit local feasibility. Specifically:
\begin{enumerate}[(i)]
    \item For each of the four coupling/decoupling scenarios, the corresponding active inequality among Constraints~(\ref{constraint:couple opposite})--(\ref{constraint:decouple same}) is necessary and sufficient for two-unit platform feasibility in that scenario;
    \item For pairwise formation continuation on a feasible turn-round arc, the active case of Constraints~(\ref{constraint:continuation lower})--(\ref{constraint:continuation upper}) is necessary and sufficient for continuation consistency under the no-resequencing and no-insertion/removal assumption.
\end{enumerate}
\end{theorem}

\begin{proof}
Part~(i) follows directly from Lemmas~\ref{lem:equiv-coupling-opposite}--\ref{lem:equiv-decoupling-same}. Part~(ii) follows from Lemma~\ref{lem:equiv-continuation}.
\end{proof}

Thus, a violated active local ordering instance is not merely a bookkeeping inconsistency in the order variables: it is a local station-operational infeasibility and therefore a blockage-causing local ordering conflict under the terminology adopted in this paper.

Theorem~\ref{thm:station-feasibility-equiv} establishes the correctness of each pairwise local requirement. For events involving more than two units, Constraints~(\ref{constraint:couple opposite})--(\ref{constraint:decouple same}) impose the corresponding pairwise feasibility conditions on every relevant unit pair. The remaining question is whether these induced pairwise requirements can be satisfied simultaneously by a single formation order on the affected trip. The next subsubsection answers this question through the acyclicity of the induced precedence digraph.

\subsubsection{Acyclicity and global feasibility of induced unit ordering}\label{sec:acyclicity-global}
Theorem~\ref{thm:station-feasibility-equiv} establishes necessary-and-sufficient local feasibility conditions for individual unit pairs. When a station event involves multiple units, the remaining issue is whether all active pairwise requirements can be satisfied simultaneously by a single formation order on the affected trip. We address this by encoding the induced pairwise requirements as a precedence digraph on the active units. The central concept is \emph{acyclicity}: it provides an exact consistency criterion for the induced ordering requirements, and it will later serve both to characterise multi-unit feasibility and to motivate the cycle-based strengthening and certification procedures used in the exact solution framework. We use \(h_1 \prec_k h_2\) to denote that unit \(h_1\) is required to be ahead of unit \(h_2\) on trip \(k\),
or equivalently \(\theta^{h_1}_k < \theta^{h_2}_k\). When the trip index is clear from context, we omit the subscript and write \(h_1 \prec h_2\). We will show that, e.g., a directed cycle such as \(h_1 \prec_k h_2 \prec_k h_3 \prec_k h_1\) means that no linear front--rear order can satisfy all induced requirements simultaneously.

For any trip \(i\), let $U_i:=\{h\in H:s_i^h=1\}$ denote the set of active units on trip \(i\) under the fixed integral assignment. For any distinct $h_1,h_2\in U_i$, define the induced pairwise precedence indicator $\rho_{h_1h_2}^i:=\mathbb{I}\{\theta_i^{h_1}<\theta_i^{h_2}\}$. Thus, $\rho_{h_1h_2}^i=1$ means that $h_1$ is ahead of $h_2$ in the coupled formation on trip $i$.

Whenever Constraints~(\ref{constraint:couple opposite})--(\ref{constraint:decouple same}) are active after coefficient substitution, each active instance fixes the sign of a difference of the form $\theta_k^{h_1}-\theta_k^{h_2}$ on a specific trip $k\in\{i,j\}$. Since active units have distinct integer positions, each such instance is equivalent to either $\theta_k^{h_1}<\theta_k^{h_2}$ or $\theta_k^{h_1}>\theta_k^{h_2}$. Accordingly, every active instance induces a strict precedence relation on the affected trip.

\begin{lemma}[Acyclicity and existence of a consistent formation order]\label{lem:acyclic_topological}
Fix a trip index $k$ and let $\prec_k$ be the strict precedence relation induced on a set $U\subseteq U_k$ by the active instances of Constraints~(\ref{constraint:couple opposite})--(\ref{constraint:decouple same}). Then the induced requirements on $U$ can be satisfied simultaneously by a single linear order if and only if $\prec_k$ is acyclic. In that case, any topological order of $\prec_k$ yields distinct positive integer positions satisfying all induced precedences on $U$.
\end{lemma}

\begin{proof}
This is the standard topological-ordering characterisation for finite directed graphs \citep{trotter2002combinatorics}. A finite directed graph admits a linear order consistent with all arcs if and only if it is acyclic. When $\prec_k$ is acyclic, any topological order of $U$ defines a ranking of the units, and assigning $\theta_k^h$ to be the rank of $h$ in that order yields distinct positive integer positions satisfying all induced precedences.
\end{proof}

Constraints~(\ref{constraint:continuation lower})--(\ref{constraint:continuation upper}) are not part of the induced precedence relation $\prec_k$. When active, they enforce $\theta_j^{h_1}-\theta_j^{h_2}=d_i d_j(\theta_i^{h_1}-\theta_i^{h_2})$, which is a cross-trip carry-over condition rather than a local precedence requirement on a single trip. We now apply this acyclicity criterion to the precedence digraph induced by a multi-unit coupling or decoupling event.

\subsubsection{A unified characterisation for multi-unit station feasibility}\label{sec:unified-station-feasibility}

Consider a coupling or decoupling station event on trip $k$, and let $U_k := \{h\in H:s^h_k=1\}$ be the active units on that trip. The active instances of Constraints~(\ref{constraint:couple opposite})--(\ref{constraint:decouple same}) induce pairwise precedences on the units involved in the event. Let E$_k:=\{(a,b)\in U_k\times U_k:\ a\prec_k b\}$, where $a\prec_k b$ denotes an induced precedence from \S~\ref{sec:acyclicity-global}, and define the directed graph $G_k:=(U_k,E_k)$. Any active unit that does not appear in an ordering instance is simply an isolated vertex of $G_k$. The cross-trip continuation constraints~(\ref{constraint:continuation lower})--(\ref{constraint:continuation upper}) are treated separately and are not part of this local precedence digraph. We assume that, for the platform-feasible events considered here, the multi-unit blockage logic decomposes into the pairwise precedence requirements generated by the corresponding unit pairs. Thus, once these pairwise requirements admit a common linear extension, the event can be realised without additional local restrictions.

\begin{theorem}[Unified characterisation for multi-unit feasibility]\label{thm:multiunit_station_feasibility}
Under the platform-feasible operating assumption, all active instances of Constraints~(\ref{constraint:couple opposite})--(\ref{constraint:decouple same}) induced by a coupling or decoupling event on trip $k$ are simultaneously satisfiable by a single formation order on $U_k$ if and only if $G_k$ is acyclic. In that case, any topological order of $G_k$ yields a feasible formation order on trip $k$.
\end{theorem}

\begin{proof}
Let $U\subseteq U_k$ denote the units that appear in at least one active ordering instance. By Lemma~\ref{lem:acyclic_topological}, the induced precedences on $U$ are simultaneously satisfiable if and only if they are acyclic. Since the units in $U_k\setminus U$ are isolated vertices of $G_k$, they can be inserted arbitrarily into any topological order of $U$, yielding a formation order on all of $U_k$. Conversely, if $G_k$ contains a directed cycle, then no linear order can satisfy all induced precedences simultaneously. The final equivalence with station feasibility follows from Theorem~\ref{thm:station-feasibility-equiv}, which identifies these induced precedences as the relevant pairwise feasibility conditions for the active coupling/decoupling event.
\end{proof}

\begin{figure}[htbp]
    \centering
    \includegraphics[width=\textwidth]
    {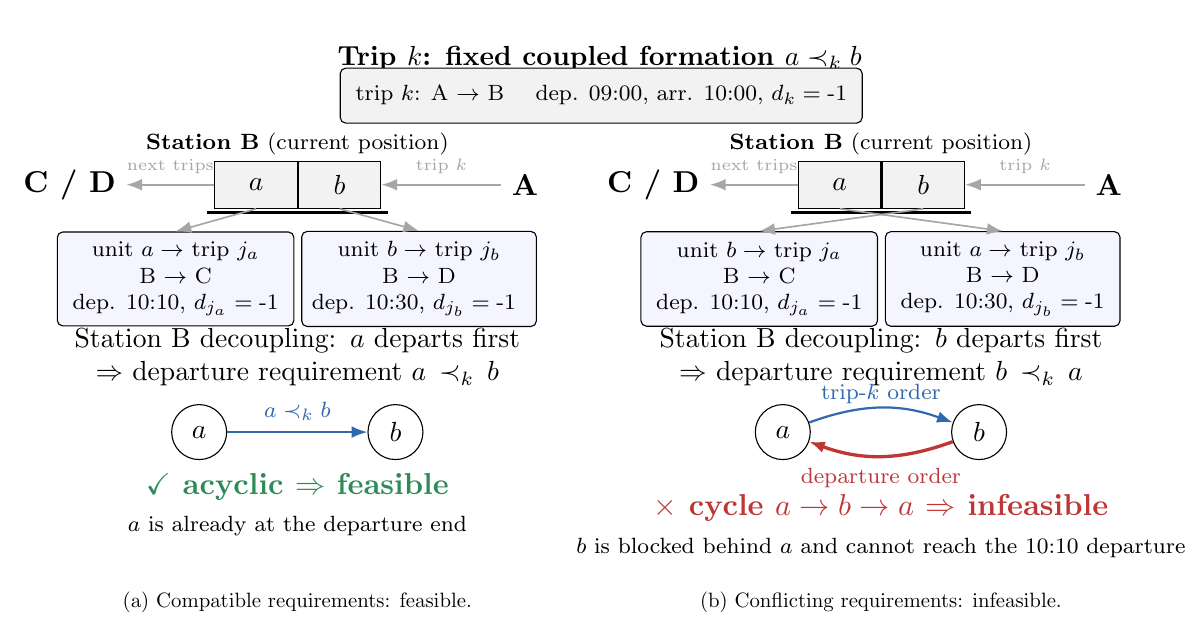}
    \caption{Two-unit illustration of Theorem~\ref{thm:multiunit_station_feasibility}. Units $a$ and $b$ operate trip $k$ as a coupled formation with $a$ ahead of $b$, i.e., $a\prec_k b$, and arrive at Station B at 10:00. In the station schematic, the units are currently located at B, with A to their right and the successor destinations C and D to their left. The two successor trips are fixed in both panels: trip $j_a$ runs from B to C and departs at 10:10, whereas trip $j_b$ runs from B to D and departs at 10:30. In (a), unit $a$ is assigned to $j_a$ and unit $b$ to $j_b$. The destination-side decoupling therefore also requires $a\prec_k b$, so the combined digraph is acyclic and the formation is feasible. In (b), the unit-to-trip assignments are exchanged: unit $b$ is assigned to $j_a$ and unit $a$ to $j_b$. The resulting requirement $b\prec_k a$ conflicts with the existing formation order $a\prec_k b$, creating the directed cycle $a\rightarrow b\rightarrow a$. Unit $b$ is consequently blocked behind unit $a$}
    \label{fig:three-unit-ordering-cases}
\end{figure}

Consider trip $k$, which is operated by two coupled train units $a$ and $b$. Suppose that the active requirements fix the formation on trip $k$ as $a\prec_k b$, so that unit $a$ is positioned ahead of unit $b$. Trip $k$ departs from Station A at 09:00 and arrives at Station B at 10:00. Two same-direction successor trips are fixed at Station B: trip $j_a$ operates from B to C and departs at 10:10, whereas trip $j_b$ operates from B to D and departs at 10:30. In Figure~\ref{fig:three-unit-ordering-cases}(a), unit $a$ is assigned to $j_a$ and unit $b$ is assigned to $j_b$. The destination-side decoupling therefore requires $a\prec_k b$, which coincides with the existing formation requirement on trip $k$. The induced precedence digraph contains only the arc $a\rightarrow b$ and is acyclic. Its topological order $(a,b)$ consequently gives a feasible formation in which unit $a$ can access the earlier departure without overtaking or blockage.

By contrast, Figure~\ref{fig:three-unit-ordering-cases}(b) keeps the timetable, directions, and destinations of $j_a$ and $j_b$ unchanged, but exchanges the unit-to-trip assignments: unit $b$ is assigned to the 10:10 departure $j_a$, while unit $a$ is assigned to the 10:30 departure $j_b$. The destination-side decoupling now requires $b\prec_k a$. Together with the existing requirement $a\prec_k b$ on trip $k$, this produces the directed cycle $a\rightarrow b\rightarrow a$. The formation $(a,b)$ leaves unit $b$ blocked behind unit $a$, whereas reversing the formation would violate the existing order on trip $k$. Hence, no linear formation order can satisfy both requirements, and no blockage-free decoupling is possible.
This two-unit cycle is the smallest obstruction characterised by Theorem~\ref{thm:multiunit_station_feasibility} and subsequently excluded by the activated cycle inequalities in \S~\ref{sec:activated-cycles}.

Although Figure~\ref{fig:three-unit-ordering-cases} uses only two train units for clarity, the same characterisation extends directly to an arbitrary number of active units. For any finite active set $U_k$, each active coupling or decoupling requirement induces a pairwise precedence arc in $G_k$. If the resulting digraph is acyclic, any topological order of $G_k$ provides a formation order satisfying all such requirements simultaneously. If $G_k$ contains a directed cycle of any length, no linear formation order can satisfy every precedence requirement in that cycle. Thus, the two-unit cycle shown in Figure~\ref{fig:three-unit-ordering-cases}(b) represents the smallest possible ordering obstruction, while Theorem~\ref{thm:multiunit_station_feasibility} applies without modification to coupling and decoupling events involving any number of train units.

\subsection{Formulation choice for the ordering layer}\label{sec:Equivalent local encodings and local LP strength}
The acyclicity result in Theorem~\ref{thm:multiunit_station_feasibility} characterises the trip-wise ordering object induced by active coupling and decoupling requirements. For a fixed active unit set $U_k$ on a trip $k$, this object can alternatively be represented by a position-assignment formulation. Binary variables $\lambda_k^{h,p}$ assign each active unit $h\in U_k$ to a position $p\in\{1,\ldots,|U_k|\}$, while standard prefix inequalities enforce the induced precedence requirements. At the integer level, the position-assignment formulation, the feasible position vectors, and the linear extensions of $G_k$ are equivalent. The formal correspondence is established in Proposition~\ref{prop:encoding-equivalence} in Appendix~\ref{app:local_encodings}.

For the LP comparison, let $\widehat Q_k^\lambda$ denote the natural LP relaxation of the fixed-trip position-assignment formulation defined in Appendix~\ref{app:local_encodings}. Let $P_{k,\theta}^{\mathrm{BM}}$ denote the projection onto the $\theta_k$-coordinates of the LP relaxation of the compact local Big-$M$ encoding in Constraints~(\ref{constraint:order2})-(\ref{constraint:order4}), restricted to the same fixed active set and the same induced precedence requirements.

\begin{theorem}[Fixed-trip LP dominance of the position-assignment encoding]
\label{thm:fixed_trip_lp_dominance}
Fix a trip $k$ and its active unit set $U_k$, and assume a valid local Big-$M$ calibration $M_k\geq |U_k|$. Then $\operatorname{proj}_{\theta_k} \bigl(\widehat Q_k^\lambda\bigr) \subseteq P_{k,\theta}^{\mathrm{BM}}$. The containment is strict in general. Hence, for a fixed active unit set, the position-assignment formulation provides a tighter local LP relaxation than the compact Big-$M$ comparison encoding.
\end{theorem}

\begin{proof}
See Appendix~\ref{app:local_encodings}.
\end{proof}

Theorem~\ref{thm:fixed_trip_lp_dominance} formalises the local relaxation-strength advantage of the position-assignment encoding. Nevertheless, this encoding is not adopted in the global formulation. In the full train unit scheduling problem, the active set $U_k$ is decision-dependent because it is determined by the unit-assignment and connection variables. Embedding the position-assignment formulation globally would therefore require variables for every trip, unit, and possible position, together with additional constraints linking them to the trip-assignment indicators. This would substantially increase the model size and couple the ordering layer more tightly to the unit-flow layer.

For this reason, the global model is written in the compact $(x,y,\theta,u)$ space. The position-assignment formulation is retained only as a fixed-trip benchmark for understanding the local ordering object and quantifying the strength sacrificed by the compact encoding. The proof of Theorem~\ref{thm:fixed_trip_lp_dominance} and the remaining fixed-trip equivalence results are provided in Appendix~\ref{app:local_encodings}. The next subsection uses the same acyclicity structure to derive compact activated cycle inequalities that partially recover this local strength without introducing the full position-assignment formulation globally.

\subsection{From acyclicity to activated cycle inequalities (ACIs)}\label{sec:activated-cycles}
The compact \((\theta,u)\)-based ordering layer avoids the size of the fixed-trip assignment formulation, but its LP relaxation can admit fractional comparison patterns. We therefore use the acyclicity characterisation of Theorem~\ref{thm:multiunit_station_feasibility} to derive activated cycle inequalities. These inequalities do not replace ordering certification; they serve as auxiliary LP strengthening for the coupling/decoupling precedence layer.

Recall that for distinct active units $a,b\in U_k$, the induced pairwise precedence indicator is $\rho_{ab}^k:=\mathbb{I}\{\theta_k^a<\theta_k^b\}$. If the local ordering requirements on trip $k$ are feasible, then every feasible formation order is a linear extension of the induced precedence digraph $G_k=(U_k,E_k)$, so no directed simple cycle can be realised. Consequently, for any directed simple cycle $C=(h_1\to h_2\to\cdots\to h_m\to h_1)$ on distinct units in $U_k$, where $m$ is the number of units in the cycle, the inequality
\begin{equation}
\sum_{r=1}^{m}\rho_{h_rh_{r+1}}^k \le m-1,
\qquad h_{m+1}:=h_1,
\label{eq:mcycle-rho}
\end{equation}
is valid.

In the full formulation, the active unit set is decision-dependent. Let $s_k^h$ denote whether unit $h$ is assigned to trip $k$. Then a globally valid activated version of (\ref{eq:mcycle-rho}) is
\begin{equation}
\sum_{r=1}^{m}\rho_{h_rh_{r+1}}^k
\le
m-1+\sum_{r=1}^{m}(1-s_k^{h_r}),
\qquad \forall k,\ \forall \text{ directed simple cycles } C.
\label{eq:activated-cycle-rho}
\end{equation}
When all vertices of $C$ are active, (\ref{eq:activated-cycle-rho}) reduces to the standard cycle inequality; if at least one vertex is inactive, the inequality becomes redundant.

Because the global model is written in the $(x,y,\theta,u)$ variables, it is convenient to express the same inequality directly in the implemented model variables. On active pairs, we have $\rho_{h_1h_2}^k=1-u_k^{h_1h_2}$, so (\ref{eq:activated-cycle-rho}) is equivalent to
\begin{equation}
\sum_{r=1}^{m} u^{k}_{h_rh_{r+1}}
\ge
\sum_{r=1}^{m} s^{k}_{h_r}-(m-1),\qquad h_{m+1}:=h_1, \forall k,\ \forall \text{ directed simple cycles } C.
\label{eq:activated-cycle-u}
\end{equation}

\begin{proposition}[Validity of activated cycle inequalities]\label{prop:activated-cycle-valid}
For any trip \(k\) and any directed simple cycle
\(C=(h_1\to h_2\to\cdots\to h_m\to h_1)\) on distinct units,
the activated cycle inequality (\ref{eq:activated-cycle-u}) is valid for the full ordered-unit formulation. When all units in \(C\)
are active, the inequality states that the corresponding cyclic ordering
pattern cannot be realised by any linear formation order.
\end{proposition}

\begin{proof}
If at least one vertex of \(C\) is inactive, then
\(\sum_{r=1}^{m} s^{k}_{h_r}\le m-1\), so the right-hand side is at most
zero and the inequality is redundant. If all vertices are active, then
\(u^{k}_{ab}=1\) means that \(a\) is behind \(b\), or equivalently that
\(a\) is not ahead of \(b\). Since a linear order cannot realise
\(h_1\prec h_2\prec\cdots\prec h_m\prec h_1\), at least one arc in the
cycle must be reversed, which is exactly
\(\sum_{r=1}^{m}u^{k}_{h_rh_{r+1}}\ge 1\). Hence the activated inequality
is valid. If the active original precedence implications additionally
force \(h_r\prec h_{r+1}\) for every arc of the cycle, then all
corresponding \(u\)-variables would have to be zero, contradicting the
cycle inequality. Therefore the original precedence implications together
with all such cycle inequalities rule out directed cycles.
\end{proof}

Activated cycle inequalities provide auxiliary LP strengthening for the coupling/decoupling precedence layer by encoding a direct consequence of the acyclicity characterisation. They are not the primary exactness mechanism of the full model. Exactness of TUSOU is enforced through ordering certification of the original local ordering logic, including the continuation layer, while dynamic separation of activated cycle inequalities is deferred to \S~\ref{sec:activated-cycles-enforcement}.

\subsection{Fixed-assignment orderability and projection exactness}\label{subsec:fixed_assignment_orderability}
Fix an integral assignment $x$. For each trip $k$, let $U_k:=\{h:s_k^h=1\}$ be the active units on trip $k$, and let $G_k=(U_k,E_k)$ be the precedence digraph induced by the active coupling/decoupling instances under $x$. A family $\Pi:=\{\pi_k:U_k\to\{1,\dots,|U_k|\}\}_{k}$ is called \emph{continuation-consistent} if:
\begin{enumerate}[(i)]
    \item Each $\pi_k$ is a topological order of $G_k$;
    \item For every active feasible arc $a=(i,j)\in\tilde A$ and every pair of units $h_1,h_2$ such that $x_a^{h_1}=x_a^{h_2}=1$; 
    
    we have $\pi_j(h_1)-\pi_j(h_2)=d_i d_j\bigl(\pi_i(h_1)-\pi_i(h_2)\bigr)$.
\end{enumerate}

The following result is stated under the platform-feasible abstraction adopted in this paper: local coupling/decoupling blockage logic decomposes into pairwise precedence requirements, same-direction time ties are resolved by exogenous priority before constraint generation, and active continuation arcs do not allow overtaking, resequencing, or insertion/removal between units that traverse the same arc.

\begin{corollary}[Fixed-assignment orderability and blockage-free ordering witnesses]\label{cor:existence-theta-under-certification}
Fix an integral assignment $x$. For each trip $k$, let $G_k=(U_k,E_k)$ be the precedence digraph induced by the active coupling/decoupling instances under $x$. Then the following are equivalent:
    \begin{enumerate}[(i)]
    \item there exist order and comparison variables \((\theta,u)\) satisfying the order-linking constraints and all active original ordering constraints among Constraints~(\ref{constraint:couple opposite})-(\ref{constraint:continuation upper});
    \item there exists a continuation-consistent family of trip-wise orders $\Pi=\{\pi_k\}_k$.
    \end{enumerate}
In particular, if every induced precedence digraph $G_k$ is acyclic and the resulting trip-wise orders can be chosen to be continuation-consistent on all active feasible arcs, then there exists a feasible order assignment $\theta$, and the corresponding schedule is blockage-free under the platform-feasible operating assumption.
\end{corollary}

\begin{proof}
(ii) $\Rightarrow$ (i): For each trip $k$ and each active unit $h\in U_k$, define $\theta_k^h:=\pi_k(h)$, and set $\theta_k^h:=0$ for $h\notin U_k$. Since each $\pi_k$ is a topological order of $G_k$, all active coupling/decoupling precedences on trip $k$ are satisfied. Since the family $\Pi$ is continuation-consistent, every active continuation instance also satisfies $\theta_j^{h_1}-\theta_j^{h_2}=d_i d_j\bigl(\theta_i^{h_1}-\theta_i^{h_2}\bigr)$. Thus all active original ordering constraints are satisfied. By Theorem~\ref{thm:station-feasibility-equiv}, Theorem~\ref{thm:multiunit_station_feasibility}, and Lemma~\ref{lem:equiv-continuation}, the resulting schedule is blockage-free under the platform-feasible operating assumption. The corresponding pairwise variables can be set as $u^{ab}_k=1$ if and only if $\theta^a_k>\theta^b_k$, so the order-linking constraints are also satisfied.

(i) $\Rightarrow$ (ii): Suppose there exists a feasible order assignment $\theta$. For each trip $k$, define $\pi_k$ to be the linear order induced by the values $\{\theta_k^h:h\in U_k\}$, i.e., $\pi_k(h_1)<\pi_k(h_2)$ if and only if $\theta_k^{h_1}<\theta_k^{h_2}$. Because every active unit receives a distinct value in $\{1,\dots,|U_k|\}$, each $\pi_k$ is a bijection from $U_k$ to $\{1,\dots,|U_k|\}$. Since $\theta$ satisfies all active coupling/decoupling constraints, each $\pi_k$ is a topological order of $G_k$. Since $\theta$ also satisfies all active continuation constraints, the family $\Pi=\{\pi_k\}_k$ is continuation-consistent. Hence (2) holds.
\end{proof}

Let \(\mathcal P\) denote the feasible region of the full ordered-unit formulation, and let \(\mathcal N\) denote the feasible region of its network-assignment layer: $\mathcal P:=\{(x,y,\theta,u):~(\ref{constraint:unit used})\text{--}(\ref{constraint:order1})\},  \mathcal N:=\{(x,y):(\ref{constraint:unit used})\text{--}(\ref{constraint:y 3}),(\ref{x}),(\ref{y})\}$. Define $X^{\mathrm{net}}:=\operatorname{proj}_x(\mathcal N)$. For any \(\bar{x}\in X^{\mathrm{net}}\), let \(L(\bar{x})\) be the set of active original ordering instances among Constraints~(\ref{constraint:couple opposite})\text{--}(\ref{constraint:continuation upper}), and define the set of ordering completions by $\Omega(\bar{x}):=\{(\theta,u):(\ref{constraint:order2})\text{--}(\ref{constraint:order4}),(\ref{constraint:order1}),\text{ and }L(\bar{x})\text{ hold}\}$. The projected orderable assignment set is $X^{\mathrm{ord}}:=\{\bar{x}\in X^{\mathrm{net}}:\Omega(\bar{x})\neq\emptyset\}$. 

\begin{proposition}[Projection exactness of the ordered-unit formulation]
\label{prop:projection_exactness}
Under the platform-feasible operating assumptions and the Big-\(M\) calibration used in the full formulation, $\operatorname{proj}_x(\mathcal P)=X^{\mathrm{ord}}$. Consequently, an integral assignment \(\bar{x}\) is the \(x\)-projection of a feasible single-stage ordered-unit schedule if and only if it admits at least one blockage-free ordering witness.
\end{proposition}

\begin{proof}
If \(\bar{x}\in\operatorname{proj}_x(\mathcal P)\), then there exist \(y,\theta,u\) such that \((\bar{x},y,\theta,u)\in\mathcal P\). Hence \((\bar{x},y)\in\mathcal N\), so \(\bar{x}\in X^{\mathrm{net}}\). Moreover, the same \((\theta,u)\) satisfies the order-linking constraints and all active original ordering instances induced by \(\bar{x}\), so \((\theta,u)\in\Omega(\bar{x})\). Thus \(\bar{x}\in X^{\mathrm{ord}}\).

Conversely, let \(\bar{x}\in X^{\mathrm{ord}}\). Then \(\bar{x}\in X^{\mathrm{net}}\), so there exists \(y\) with \((\bar{x},y)\in\mathcal N\). Since \(\Omega(\bar{x})\neq\emptyset\), there exists \((\theta,u)\in\Omega(\bar{x})\). The active ordering instances among Constraints~(\ref{constraint:couple opposite})-(\ref{constraint:continuation upper}) hold by definition of \(\Omega(\bar{x})\), and the inactive instances are relaxed by the calibrated Big-\(M\) terms. Therefore \((\bar{x},y,\theta,u)\in\mathcal P\), which implies \(\bar{x}\in\operatorname{proj}_x(\mathcal P)\).
\end{proof}

Proposition~\ref{prop:projection_exactness} shows that ordering certification is an exact membership test for the \(x\)-projection of the full ordered-unit feasible region. Thus rejecting \(\bar{x}\) means that the selected unit movements cannot be completed into any blockage-free ordered-unit schedule, rather than merely that a particular constructed ordering has failed.

\S~\ref{sec:Enhanced formulation} therefore establishes the structural basis of the paper. The ordered-unit ILP captures executable unit movements and local ordering requirements in a single stage; multi-unit coupling/decoupling feasibility is exactly characterised by acyclicity of the induced trip-wise precedence digraphs; and blockage-free schedules exist precisely when these trip-wise orders can be chosen to be continuation-consistent on the active feasible arcs. Activated cycle inequalities follow as valid strengthening consequences of the coupling/decoupling layer, but they do not by themselves enforce the full local ordering logic. This distinction motivates the exact solution framework developed next.

%% file: paper/5_alg.tex
Building on Corollary~\ref{cor:existence-theta-under-certification}, this section turns the structural characterisation of \S~\ref{sec:Enhanced formulation} into an exact certification-equipped branch-and-bound-and-cut framework. TUSOU searches over network-feasible unit assignments, uses ordering certification to decide whether an integer assignment \(\bar x\) admits at least one blockage-free ordering witness, and lazily recovers missing active ordering instances when a non-orderable assignment is encountered. Within this framework, ordering certification is the exactness mechanism, lazy recovery restores missing original ordering rows when a non-orderable assignment is encountered, and activated cycle inequalities strengthen the coupling/decoupling precedence relaxation.

\subsection{Fixed-assignment ordering certification and exact node processing}\label{sec:incumbent-certification}
We use ordering certification to mean the fixed-assignment test induced by Corollary~\ref{cor:existence-theta-under-certification}. Given an integral unit-assignment vector \(\bar x\), let \(s_i^h(\bar x)\) be the induced trip-assignment indicator of unit \(h\) on trip \(i\). Let \(\mathcal L^{\mathrm{cd}}(\bar x)\) denote the active coupling/decoupling instances among Constraints~(\ref{constraint:couple opposite})--(\ref{constraint:decouple same}), and let \(\mathcal L^{\mathrm{cont}}(\bar x)\) denote the active continuation instances among Constraints~(\ref{constraint:continuation lower})--(\ref{constraint:continuation upper}). These active instances are induced by the selected predecessor/successor arcs and by pairs of units traversing the same feasible arc under \(\bar x\). Define \(\mathcal L(\bar x):=\mathcal L^{\mathrm{cd}}(\bar x)\cup\mathcal L^{\mathrm{cont}}(\bar x)\). Let \(\mathcal O\) denote the order-domain and order-linking subsystem, i.e., Constraints~(\ref{constraint:order2})--(\ref{constraint:order4}) and~(\ref{constraint:order1}). We write \(z\models\mathcal C\) if \(z\) satisfies every constraint in the constraint set \(\mathcal C\). Then, for a fixed integral assignment \(\bar{x}\),
\begin{equation}
\Omega(\bar{x}):=\{(\theta,u):(\bar{x},\theta,u)\models \mathcal O\cup\mathcal L(\bar{x})\}.
\label{eq:fixed-assignment-ordering-set}
\end{equation}

For an integral assignment \(\bar{x}\), let \(U_k(\bar{x}):=\{h:s_k^h(\bar{x})=1\}\) be the active unit set on trip \(k\), and let \(q_k:=|U_k(\bar{x})|\). Let \(G_k(\bar{x})=(U_k(\bar{x}),E_k(\bar{x}))\) be the active coupling/decoupling precedence digraph induced by \(\mathcal L^{\mathrm{cd}}(\bar{x})\), where \((a,b)\in E_k(\bar{x})\) means that the active original ordering logic requires unit \(a\) to be ahead of unit \(b\) on trip \(k\). If every \(G_k(\bar{x})\) is acyclic, define the fixed-assignment ordering feasibility problem \(\mathcal C(\bar{x})\) as follows. For each \(k\), \(h\in U_k(\bar{x})\), and \(r=1,\ldots,q_k\), introduce \(\lambda_{khr}\in\{0,1\}\), and define \(p_k^h:=\sum_{r=1}^{q_k} r\lambda_{khr}\). The constraints of \(\mathcal C(\bar{x})\) are
\begin{align}
\sum_{r=1}^{q_k}\lambda_{khr} &= 1, && \forall k,\ h\in U_k(\bar{x}), \label{eq:cert_unit_position}\\
\sum_{h\in U_k(\bar{x})}\lambda_{khr} &= 1, && \forall k,\ r=1,\ldots,q_k, \label{eq:cert_position_unit}\\
p_k^a+1 &\le p_k^b, && \forall k,\ (a,b)\in E_k(\bar{x}), \label{eq:cert_precedence}\\
p_j^{h_1}-p_j^{h_2} &= d_i d_j\bigl(p_i^{h_1}-p_i^{h_2}\bigr), && \forall a=(i,j)\in\widetilde A,\ \forall h_1\neq h_2:\bar{x}_a^{h_1}=\bar{x}_a^{h_2}=1. \label{eq:cert_continuation}
\end{align}

\begin{algorithm}[H]
\caption{Ordering certification for a fixed assignment}
\label{alg:ordering_certification}
\KwIn{Integral unit-assignment vector \(\bar{x}\).}
\KwOut{\textsc{Accept} with an ordering witness \((\theta^\star,u^\star)\), or \textsc{Reject} with an infeasible active ordering subsystem.}
Construct \(U_k(\bar{x})\) and \(G_k(\bar{x})=(U_k(\bar{x}),E_k(\bar{x}))\) for every trip \(k\).\;
\If{some \(G_k(\bar{x})\) contains a directed cycle \(C\)}{Return \textsc{Reject} with the active original coupling/decoupling instances that generated \(C\).\;}
Build and solve the fixed-assignment feasibility problem \(\mathcal C(\bar{x})\) defined by Constraints~\eqref{eq:cert_unit_position}--\eqref{eq:cert_continuation}.\;
\If{\(\mathcal C(\bar{x})\) is feasible}{For each active unit \(h\in U_k(\bar{x})\), set \(\theta_k^{h,\star}:=p_k^h\), and set \(\theta_k^{h,\star}:=0\) for \(h\notin U_k(\bar{x})\).\; For every trip \(k\) and every distinct pair \(a,b\in H\), set \(u_k^{ab,\star}:=1\) if \(\theta_k^{a,\star}>\theta_k^{b,\star}\), and set \(u_k^{ab,\star}:=0\) otherwise.\; Return \textsc{Accept} with \((\theta^\star,u^\star)\).\;}
\Else{Return \textsc{Reject} with an infeasible subsystem of \(\mathcal C(\bar{x})\), mapped back to the corresponding active original ordering instances in \(\mathcal L(\bar{x})\). If no minimal subsystem is extracted, the full active set \(\mathcal L(\bar{x})\) can be used for lazy recovery.\;}
\end{algorithm}

By construction, Constraints~\eqref{eq:cert_unit_position}--\eqref{eq:cert_position_unit} assign every active unit to exactly one formation position and every formation position to exactly one active unit on each trip. Constraint~\eqref{eq:cert_precedence} enforces all active coupling/decoupling precedence requirements, while Constraint~\eqref{eq:cert_continuation} enforces all active continuation carry-over equalities. Hence \(\mathcal C(\bar{x})\) is feasible if and only if the active trip-wise precedence digraphs admit a continuation-consistent family of topological orders. By Corollary~\ref{cor:existence-theta-under-certification}, this is equivalent to \(\Omega(\bar{x})\neq\emptyset\). By Proposition~\ref{prop:projection_exactness}, it is also equivalent to \(\bar{x}\) belonging to the \(x\)-projection of the full ordered-unit feasible region. Thus Algorithm~\ref{alg:ordering_certification} accepts exactly the orderable assignments and rejects exactly the non-orderable assignments.

If \(\Omega(\bar x)\neq\emptyset\), TUSOU completes the candidate assignment with the certified witness \((\theta^\star,u^\star)\) and accepts it as a feasible ordered-unit schedule. Since the objective depends on deployment, passenger-service mileage, and empty-running decisions, and not on the ordering variables, replacing tentative relaxation values by the certified witness does not change the objective value. If \(\Omega(\bar x)=\emptyset\), the candidate assignment is rejected. Ordering certification identifies an infeasible active ordering subsystem, such as a directed precedence cycle or a continuation-inconsistent set of carry-over requirements. TUSOU then recovers the original formulation lazily by adding one or more globally valid original ordering instances associated with this subsystem. These recovered rows are original constraints of the ordered-unit formulation and therefore cannot remove any solution feasible for the full model.

\begin{proposition}[Exactness and finite termination of TUSOU] \label{prop:tusou_exact_termination}
Consider a fixed finite instance. Suppose TUSOU uses Algorithm~\ref{alg:ordering_certification} for fixed-assignment ordering certification, adds only original ordering rows from Constraints~(\ref{constraint:couple opposite})--(\ref{constraint:continuation upper}) or globally valid ACI cuts, and applies a complete branching rule over the finite integer variable set. If a rejected assignment \(\bar{x}\) yields no previously absent original ordering row for lazy recovery, TUSOU adds the valid assignment-exclusion cut
\[
\sum_{(h,a):\bar{x}_a^h=1}(1-x_a^h)
+
\sum_{(h,a):\bar{x}_a^h=0}x_a^h
\ge 1.
\label{eq:no_good_assignment}
\]
Then TUSOU finitely terminates. Every accepted incumbent is feasible for the full ordered-unit formulation. If TUSOU terminates with zero optimality gap, the returned schedule is optimal over the projected orderable assignment set \(X^{\mathrm{ord}}\), equivalently over the full ordered-unit formulation.
\end{proposition}

\begin{proof}
If Algorithm~\ref{alg:ordering_certification} accepts an integral assignment \(\bar{x}\), it returns an ordering witness \((\theta^\star,u^\star)\in\Omega(\bar{x})\). By Proposition~\ref{prop:projection_exactness}, the completed solution is feasible for the full ordered-unit formulation. If Algorithm~\ref{alg:ordering_certification} rejects \(\bar{x}\), then \(\Omega(\bar{x})=\emptyset\), so \(\bar{x}\notin X^{\mathrm{ord}}\). Rejecting such an assignment therefore cannot remove any feasible ordered-unit schedule.

Lazy recovery adds only original ordering rows of the full formulation, and ACI cuts are globally valid by Proposition~\ref{prop:activated-cycle-valid}. Hence these cuts cannot remove any full-model feasible solution. If no previously absent original ordering row is available, the no-good cut \eqref{eq:no_good_assignment} excludes only the already certified non-orderable assignment \(\bar{x}\), and is therefore valid for \(X^{\mathrm{ord}}\).

For a fixed finite instance, there are finitely many integer assignments, finitely many potential original ordering rows, and finitely many variables on which the complete branching rule can branch. Thus the search cannot generate an infinite sequence of distinct lazy recoveries, assignment exclusions, and branching decisions. TUSOU therefore finitely terminates. Since node bounds are computed from relaxations strengthened only by valid cuts, and incumbents are accepted only after certification, the reported optimality gap is valid for the full ordered-unit problem.
\end{proof}

To reduce per-node overhead, TUSOU maintains two LP templates with the same variable set. Both templates retain the unit-assignment variables, the auxiliary ordering variables, and the order-linking/comparison layer, because these variables support ACI separation and lazy recovery of original ordering rows. The Shallow template omits the expensive local-ordering Constraints~(\ref{constraint:couple opposite})--(\ref{constraint:decouple same}) and Constraints~(\ref{constraint:continuation lower})--(\ref{constraint:continuation upper}) during LP bounding, except for globally recovered ordering rows that have already been added lazily. The Deep template includes the local-ordering layer for ordering-enriched bounding. Thus, Shallow does not ignore ordering feasibility; it delays the most expensive local-ordering rows until they are needed, while every integer assignment candidate is still certified by \eqref{eq:fixed-assignment-ordering-set} before it can be accepted.

The Deep template is activated once node depth reaches a prescribed threshold \(D_{\min}\), scaled with instance size as
\begin{equation}
D_{\min}=\lceil \alpha n_x\rceil,
\end{equation}
where \(n_x\) is the number of arc decision variables and \(\alpha\in(0,1)\) is fixed across instances. This trigger affects only LP bounding. Exactness is preserved because the timing of LP-level local-ordering enforcement does not affect the acceptance rule: every integer assignment candidate must pass ordering certification, regardless of whether the current LP bound was computed with the Shallow or Deep template.

At each node, TUSOU solves an LP relaxation to obtain a valid bound. During LP processing, globally valid ACI cuts may be separated for the coupling/decoupling precedence layer, as described in \S~\ref{sec:activated-cycles-enforcement}. Whenever an integer assignment candidate is encountered, TUSOU invokes ordering certification. If a witness \((\theta^\star,u^\star)\) is returned, the candidate is accepted and completed with that witness. If ordering certification proves \(\Omega(\bar x)=\emptyset\), the candidate is rejected and missing original ordering instances from an infeasible active subsystem are added lazily. 

\subsection{Separation and branch-and-bound-and-cut enforcement of ACI cuts}\label{sec:activated-cycles-enforcement}
The ACI cuts introduced in \S~\ref{sec:Enhanced formulation} are globally valid for the coupling/decoupling precedence layer, but their full family is too large to include explicitly. TUSOU therefore separates them dynamically during LP processing. These cuts are used only as auxiliary relaxation strengthening: they do not apply to the cross-trip continuation Constraints~(\ref{constraint:continuation lower})--(\ref{constraint:continuation upper}), and they do not replace ordering certification at integer assignment candidates.

At an LP node, separation is performed trip-wise on a restricted candidate set. For a fixed trip \(k\), let \(s_k^h\) and \(u_k^{pq}\) denote the current LP values of the trip-assignment and pairwise comparison variables. We first select \(S_k:=\{h\in H:s_k^h\ge \tau_s\}\) and keep at most \(N_{\max}\) units with the largest \(s_k^h\)-values. If fewer than two units remain, no ACI cut is separated for this trip. On the selected units, we construct a complete directed graph with arc length $w_{pq}^k:=u_{pq}^k+(1-s_k^p)$. For a directed simple cycle \(C=(h_1\to h_2\to\cdots\to h_m\to h_1)\), the activated cycle inequality~(\ref{eq:activated-cycle-u}) is violated exactly when $L_k(C):=\sum_{r=1}^m w_{h_r h_{r+1}}^k<1,\qquad h_{m+1}:=h_1$. The separation routine therefore searches for low-weight directed simple cycles on \(S_k\) and adds the most violated ACI cuts subject to a cut budget. Here is Algorithm~\ref{alg:sep-activated-cycles}.

\begin{algorithm}[H]
\small
\caption{Separation of ACI cuts for a fixed trip $k$}
\label{alg:sep-activated-cycles}
\KwIn{Fractional LP solution $(s,u)$ at a node; threshold $\tau_s$; candidate limit $N_{\max}$; cut budget $M_{\max}$; violation tolerance $\varepsilon_{\mathrm{cut}}$.}
\KwOut{A set of violated ACI cuts of the form Constraint~(\ref{eq:activated-cycle-u}).}

Select $S_k$ as the $N_{\max}$ units with largest $s_k^h$ among those satisfying $s_k^h\ge \tau_s$.\;

\If{$|S_k|<2$}{
Stop and return the empty cut set.\;
}

Construct the directed graph $D_k=(S_k,A_k)$ with $A_k:=\{(p,q):p,q\in S_k,\ p\ne q\}$\;

Assign each arc $(p,q)\in A_k$ the length $w_{pq}^k:=u_k^{pq}+(1-s_k^p)$\;

Initialise an empty cycle pool $\mathcal{P}_k$.\;

Search for directed simple cycles $C$ in $D_k$ with length $L_k(C):=\sum_{(p,q)\in C} w_{pq}^k < 1-\varepsilon_{\mathrm{cut}}$ \;

Add the detected violated cycles to $\mathcal{P}_k$ and sort them by decreasing violation $1-L_k(C)$\;

Add up to $M_{\max}$ cuts corresponding to the most violated cycles in $\mathcal{P}_k$.\;

\ForEach{selected cycle $C=(h_1\to\cdots\to h_m\to h_1)$}{
Add the activated cycle inequality
\[
\sum_{r=1}^{m} u_k^{h_rh_{r+1}}
\ge
\sum_{r=1}^{m}s_k^{h_r}-(m-1),
\qquad h_{m+1}:=h_1.
\]
}
\end{algorithm}

The restricted separation routine is heuristic. If no violated ACI cut is found on the selected candidate set, TUSOU proceeds without adding an ACI cut at that LP node. This has no effect on correctness, because every added ACI cut is globally valid for the ordered-unit formulation and the acceptance of integer assignment candidates is governed solely by ordering certification in \S~\ref{sec:incumbent-certification}. If a candidate admits an ordering witness, it can be accepted; if no ordering witness exists, the candidate is rejected and missing original ordering instances from an infeasible active subsystem are recovered lazily as described in \S~\ref{sec:incumbent-certification}.

In implementation, we do not instantiate any \(\lambda\)-variables or dense position-occupancy matrices. LP bounding, ACI separation, and branching are carried out in the compact \((x,y,\theta,u)\)-based representation. ACI cuts use only the assignment and pairwise comparison layer, and are separated only for coupling/decoupling precedence patterns; continuation consistency is handled by ordering certification rather than by cycle cuts.

%% file: paper/6_Experiment.tex
\subsection{Computational study overview}\label{sec:comp_overview}
This section evaluates the ordered-unit formulation and the TUSOU exact solution framework on real-world-derived instances. The computational study has three purposes. First, it assesses the ability of the proposed single-stage model to produce certified blockage-free unit schedules on real-world  passenger-rail instances. Second, it quantifies the role of ordering certification by recording how often assignment-feasible integer candidates fail the original local ordering logic. Third, it examines the computational effect of enforcing the local-ordering layer at different stages of the branch-and-bound tree and benchmarks TUSOU against a direct full-model Gurobi baseline using a set of TUSOU variants by varying the degree of local ordering enforcement (“template modes”).

All TUSOU variants use the same exact certification mechanism, lazy recovery of violated original ordering instances, and activated cycle inequalities as auxiliary LP-strengthening cuts. The ablation therefore varies only the timing of local-ordering enforcement during LP bounding: Always, Depth, and Shallow. The Gurobi baseline solves the full ordered-unit formulation directly.

At fractional nodes, TUSOU applies maximum-fractionality branching to the binary unit-arc variables. Specifically, each fractional variable \(x_a^h\) is assigned the score
\[
    \operatorname{score}(x_a^h)
    =
    0.5-\left|x_a^h-0.5\right|,
\]
and a variable with the largest score is selected, subject to an integrality tolerance of \(10^{-6}\). Thus, variables whose LP values are closest to \(0.5\) are prioritised. Ties are first resolved using the fractionality of the associated aggregate trip-unit assignment, after which the deterministic variable order is retained. Branching creates the two child disjunctions \(x_a^h=0\) and \(x_a^h=1\).

We employ a two-phase node-search strategy. Phase~1 uses depth-priority node selection and terminates once the first assignment passing ordering certification is obtained. Phase~2 resumes the remaining search queue from this certified incumbent and applies best-bound node selection to improve the global bound. Node bounds are obtained from LP relaxations solved using the persistent HiGHS interface. The relative optimality-gap tolerance is set to \(10^{-3}\), while the integrality and node-comparison tolerances are set to \(10^{-6}\).

Primal improvement heuristics are disabled during Phase~1 and enabled during Phase~2. Local branching with a Hamming-radius parameter of \(24\) and a RINS-like fixing heuristic with tolerance \(10^{-6}\) are invoked every \(50\) processed fractional nodes up to depth \(40\), with a \(20\)-second time limit for each heuristic subproblem. Any solution generated by these heuristics is subjected to the same fixed-assignment ordering certification as a regular integer candidate before being accepted as an incumbent.

TUSOU was implemented in Python~3.13. Its branch-and-bound-and-cut layer was built on top of pybnb, and LP relaxations were solved with HiGHS. All experiments were run on the University of Leeds high-performance computing system on the same class of CPU nodes with  AMD EPYC~9634 processors and each run was restricted to 58 solver threads. Direct full-model baselines were solved with Gurobi~13.0.1 using the same node class and thread limit. 

\subsection{Instances and data}\label{sec:instances_data}
We use five real-world-derived instances based on TransPennine Express operations. The first instance, ASC-Base, is extracted from the Anglo-Scottish route and represents the base corridor case. The second instance, TPE-Combined, combines services from the Anglo-Scottish route and the South TransPennine route, resulting in a larger timetable and connection network. The third instance, ASC-Scaled, uses the same station set, trip set, and feasible connection structure as ASC-Base, but increases the number of available units and creates a more unit-intensive stress instance. The fourth and fifth instance, NTP and NTP2, are extracted from the North TransPennine route. The fifth instance is the largest in terms of unit-indexed arcs and is used to test scalability on a denser connection network.

Table~\ref{tab:instance-summary} summarises the instance sizes. For consistency with the implemented unit-indexed connection network, the arc counts are reported after expansion over available individual train units, rather than as unexpanded timetable-level connection possibilities.

\begin{table}[t]
\centering
\caption{Summary of the real-world-derived instances.}
\label{tab:instance-summary}
\begin{tabular}{lrrrrr}
\toprule
Instance & \#Stations & \#Trips & Direction & \#Unit-indexed arcs & \#Units \\
\midrule
ASC-Base & 8 & 32 & double & 4,360 & 10/10/10/10 \\
TPE-Combined & 12 & 73 & double & 13,040 & 10/10/10/10 \\
ASC-Scaled & 8 & 32 & double & 8,720 & 20/20/20/20 \\
NTP & 6 & 92 & double & 39,645 & 15/10/10/10 \\
NTP2 & 8 & 130 & double & 101,680 & 20/20/20/20 \\
\bottomrule
\end{tabular}
\end{table}

\subsection{Controlled ablation experiments}\label{sec:exp}
The computational comparison is conducted on the five instances described in \S~\ref{sec:instances_data}. We compare four configurations for each instance. The first three are TUSOU configurations corresponding to the LP-template enforcement modes Always, Depth, and Shallow. The fourth is a direct Gurobi baseline that solves the full ordered-unit formulation.

The three TUSOU template modes differ in how the local-ordering layer is enforced during LP bounding. Always enforces the full local-ordering template at every node. Depth uses a shallow template in the upper part of the search tree and activates the deeper local-ordering template after a prescribed depth. Shallow uses the shallow template throughout the tree and relies on ordering certification and lazy recovery to enforce the original ordering logic exactly. For each configuration, Table~\ref{tab:expB} reports wall-clock time, final optimality gap, explored branch-and-bound nodes, the number of activated cycle inequalities added, and the number of integer assignment candidates seen and rejected by ordering certification.

\begin{table}[htp]
\centering
\scriptsize
\setlength{\tabcolsep}{2pt}
\renewcommand{\arraystretch}{0.95}
\caption{Computational comparison on the real-world-derived instances}
\label{tab:expB}
\begin{tabular}{l l l r r r r c c}
\hline
Instance 
& Solver 
& \shortstack{Template\\mode} 
& Time (s) 
& Gap (\%) 
& Nodes
& ACI cuts 
& \shortstack{Candidates\\(seen/rej)} 
& Rej.\% \\
\hline
ASC-Base & TUSOU & Always  & 329.73  & 0.00 & 76 & 0  & 2/0    & 0.0 \\
ASC-Base & TUSOU & Depth    & 101.19  & 0.00 & 246 & 11 & 17/11  & 64.7 \\
\best{ASC-Base} & \best{TUSOU} & \best{Shallow}  & \best{1.91}   & \best{0.00} & \best{151} & \best{8}  & \best{10/8}   & \best{80.0} \\
ASC-Base & Gurobi & - & 92.12 & 0.00 & - & -- & -- & -- \\
\hline
TPE-Combined & TUSOU & Always & $>86400^{\ddagger}$ & --   & -- & -- & --/-- & -- \\
TPE-Combined & TUSOU & Depth & $>86400^{\ddagger}$ & --   & --  & -- & --/-- & -- \\
\best{TPE-Combined} & \best{TUSOU} & \best{Shallow} & \best{89.83} & \best{0.00} & \best{2939} & \best{7} & \best{9/7} & \best{77.8} \\
TPE-Combined & Gurobi & - & 1066.18 & 0.00 & - & -- & -- & -- \\
\hline
ASC-Scaled & TUSOU & Always & $4702.90^{\dagger}$ & 13.65 & 5 & 0  & 1/0      & 0.0 \\
ASC-Scaled & TUSOU & Depth & $3762.71^{\dagger}$ & 12.36 & 829 & 0   & 1/0      & 0.0 \\
\best{ASC-Scaled} & \best{TUSOU} & \best{Shallow} & \best{12.92} & \best{0.00} & \best{646} & \best{6} & \best{8/6}   & \best{75.0} \\
ASC-Scaled & Gurobi & - & 548.73 & 0.00 & - & -- & -- & -- \\
\hline
NTP & TUSOU & Always & $>86400^{\ddagger}$ & --   & -- & -- & --/-- & -- \\
NTP & TUSOU & Depth & $>86400^{\ddagger}$ & --   & -- & -- & --/-- & -- \\
\best{NTP} & \best{TUSOU} & \best{Shallow} & \best{457.90} & \best{0.00} & \best{1368} & \best{6} & \best{9/6}   & \best{66.7} \\
NTP & Gurobi & - & 9285.99 & 0.00 & - & -- & -- & -- \\
\hline
NTP2 & TUSOU & Always & $>86400^{\ddagger}$ & --   & -- & -- & --/-- & -- \\
NTP2 & TUSOU & Depth & $>86400^{\ddagger}$ & --   & -- & -- & --/-- & -- \\
\best{NTP2} & \best{TUSOU} & \best{Shallow} & \best{2597.59} & \best{0.00} & \best{101} & \best{1} & \best{2/1}   & \best{50.0} \\
NTP2 & Gurobi & - & $>86400^{\ddagger}$ & -- & - & -- & -- & -- \\
\hline
\end{tabular}
\vspace{1mm}
\begin{minipage}{0.98\textwidth}
\footnotesize
\emph{Notes.} Bold rows indicate the fastest configuration attaining zero reported gap for each instance. For TUSOU rows, Cand. reports integer assignment candidates submitted to fixed-assignment ordering certification; rejected means that ordering certification proves \(\Omega(\bar x)=\emptyset\). All TUSOU rows use ACI separation and ordering certification; Gurobi solves the full ordered-unit model directly without TUSOU callbacks. $^\dagger$ stopped with a positive final gap; $^\ddagger$ exceeded the 86,400-second time limit.

\end{minipage}
\end{table}
\subsection{Computational performance and component ablation}
\subsubsection{Overall performance and Gurobi comparison}
Table~\ref{tab:expB} shows that selective local-ordering enforcement remains computationally important within the certification-equipped framework. On ASC-Base, all three TUSOU configurations prove optimality, and Shallow is the fastest at 1.91 seconds, compared with 329.73 seconds for Always, 101.19 seconds for Depth, and 92.12 seconds for the direct Gurobi baseline. On TPE-Combined, Always and Depth exceed the 86,400-second time limit, whereas Shallow solves the instance to zero reported gap in 89.83 seconds; the direct Gurobi baseline solves the same instance in 1066.18 seconds. On ASC-Scaled, Always and Depth stop with positive final gaps of 13.65\% and 12.36\%, respectively, while Shallow closes the gap in 12.92 seconds; the direct Gurobi baseline takes 548.73 seconds. On NTP, Always and Depth again exceed the time limit, while Shallow solves the instance in 457.90 seconds, compared with 9285.99 seconds for Gurobi. On NTP2, Always, Depth, and the direct Gurobi baseline all exceed the time limit, whereas Shallow solves the instance to zero reported gap in 2597.59 seconds.

The calibrated full ordered-unit formulation is strong enough for the direct Gurobi baseline to solve the first four instances, often at the root node. However, on the largest NTP2 instance, the direct full-model baseline exceeds the time limit. TUSOU Shallow is faster than the direct full-model baseline on all five instances and is the only tested approach that solves NTP2 within the time limit. These results indicate the effectiveness of the TUSOU architecture: it provides ordering certification diagnostics, supports selective local-ordering enforcement, and yields a structure-aware exact framework for the ordered-unit problem.

\subsubsection{Effect of selective local-ordering enforcement}
The comparison across Always, Depth, and Shallow highlights the computational value of separating exactness from uniform local-ordering enforcement. Always enforces the full local-ordering layer at every node and becomes computationally expensive on the larger TPE-Combined, NTP, and NTP2 instances. Depth reduces this burden by delaying full local-ordering enforcement, but it still exceeds the time limit on TPE-Combined, NTP, and NTP2 and stops with a positive gap on ASC-Scaled. Shallow, by contrast, relies most strongly on ordering certification and lazy recovery, and is the most reliable TUSOU template across the five reported instances. It is the only TUSOU template that solves all five instances to zero reported gap. These results support the central design of TUSOU: exactness is enforced at integer assignment candidates through ordering certification, while the expensive local-ordering layer need not be imposed uniformly during LP bounding.

\subsubsection{Role of activated cycle inequalities}
In the Shallow runs, ACI separation is active on all five instances and generates 8 cuts on ASC-Base, 7 on TPE-Combined, 6 on ASC-Scaled, 6 on NTP, and 1 on NTP2. Although these counts are modest, ACI cuts make a distinct structural contribution to TUSOU. First, they transfer a necessary condition for trip-wise orderability—the absence of directed precedence cycles—into the LP relaxation, allowing cyclic comparison patterns to be eliminated before a complete integer assignment is submitted to ordering certification. Second, each ACI is globally valid and remains available throughout the search, so one separated inequality may strengthen multiple subsequent nodes. In this respect, ACI separation complements lazy recovery, which is triggered only after ordering certification encounters a non-orderable integer assignment.

The modest number of generated cuts also indicates that ACI separation is selective rather than burdensome: TUSOU introduces this additional ordering structure only where the current relaxation exhibits a relevant cyclic pattern. ACI cuts therefore provide a low-volume, structure-specific strengthening layer, while ordering certification retains responsibility for exactness and for detecting more general infeasibility involving continuation consistency.

\subsection{Ordering certification and blockage elimination}
We next examine the role of ordering certification. A rejected candidate is an integer assignment candidate \(\bar x\) for which ordering certification proves \(\Omega(\bar x)=\emptyset\). Thus, the selected unit movements admit no continuation-consistent family of trip-wise orders and no blockage-free ordering witness under the platform-feasible operating assumptions. Such candidates are feasible with respect to the relaxed network-level search model but station-infeasible at the ordering level.

Table~\ref{tab:certification-diagnostics} aggregates the ordering certification diagnostics across the terminated TUSOU configurations reported in Table~\ref{tab:expB}. Across these nine configurations, ordering certification rejects 39 of 59 integer assignment-candidate encounters, corresponding to an overall rejection rate of 66.1\%. In the Shallow configurations alone, which rely most directly on ordering certification and lazy recovery, 28 of 38 encounters are rejected, corresponding to a rejection rate of 73.7\%. These rejection rates show that non-orderable assignments are common in the relaxed search space and that ordering feasibility is not suitable to be treated as post-processing should global optimality be guaranteed.

\begin{table}[t]
\centering
\small
\caption{Aggregated ordering certification diagnostics across terminated TUSOU configurations with candidate diagnostics.}
\label{tab:certification-diagnostics}
\begin{tabular}{lrrrr}
\toprule
Instance & Terminated configurations & Candidates seen & Rejected & Rejection rate \\
\midrule
ASC-Base & 3 & 29 & 19 & 65.5\% \\
TPE-Combined & 1 & 9 & 7 & 77.8\% \\
ASC-Scaled & 3 & 10 & 6 & 60.0\% \\
NTP & 1 & 9 & 6 & 66.7\% \\
NTP2 & 1 & 2 & 1 & 50.0\% \\
\midrule
Total & 9 & 59 & 39 & 66.1\% \\
\bottomrule
\end{tabular}
\end{table}

These diagnostics also clarify why Shallow is effective. Because Shallow does not enforce all local-ordering constraints during LP bounding, it may encounter integer assignment candidates that are network-feasible but not station orderable. Ordering certification filters these candidates much more efficiently while maintaining the exactness of B\&B: candidates that admit a blockage-free ordering witness can be accepted, while candidates with \(\Omega(\bar x)=\emptyset\) are rejected and used to recover missing original ordering instances. The results therefore support the main algorithmic design of TUSOU: avoid the cost of enforcing the complete local-ordering layer at every node, but certify every integer assignment candidate before it can become an accepted schedule.

\subsection{Representative certified schedule}
To complement the aggregate computational results, we inspect one representative certified solution from the TPE-Combined real-world-derived instance. The purpose of this subsection is not to introduce an additional synthetic test case, but to show how the schedules produced by TUSOU instantiate the local ordering mechanisms developed in \S~\ref{sec:Enhanced formulation}.

\begin{sidewaysfigure}[p]
    \centering
    \includegraphics[width=0.95\textheight]{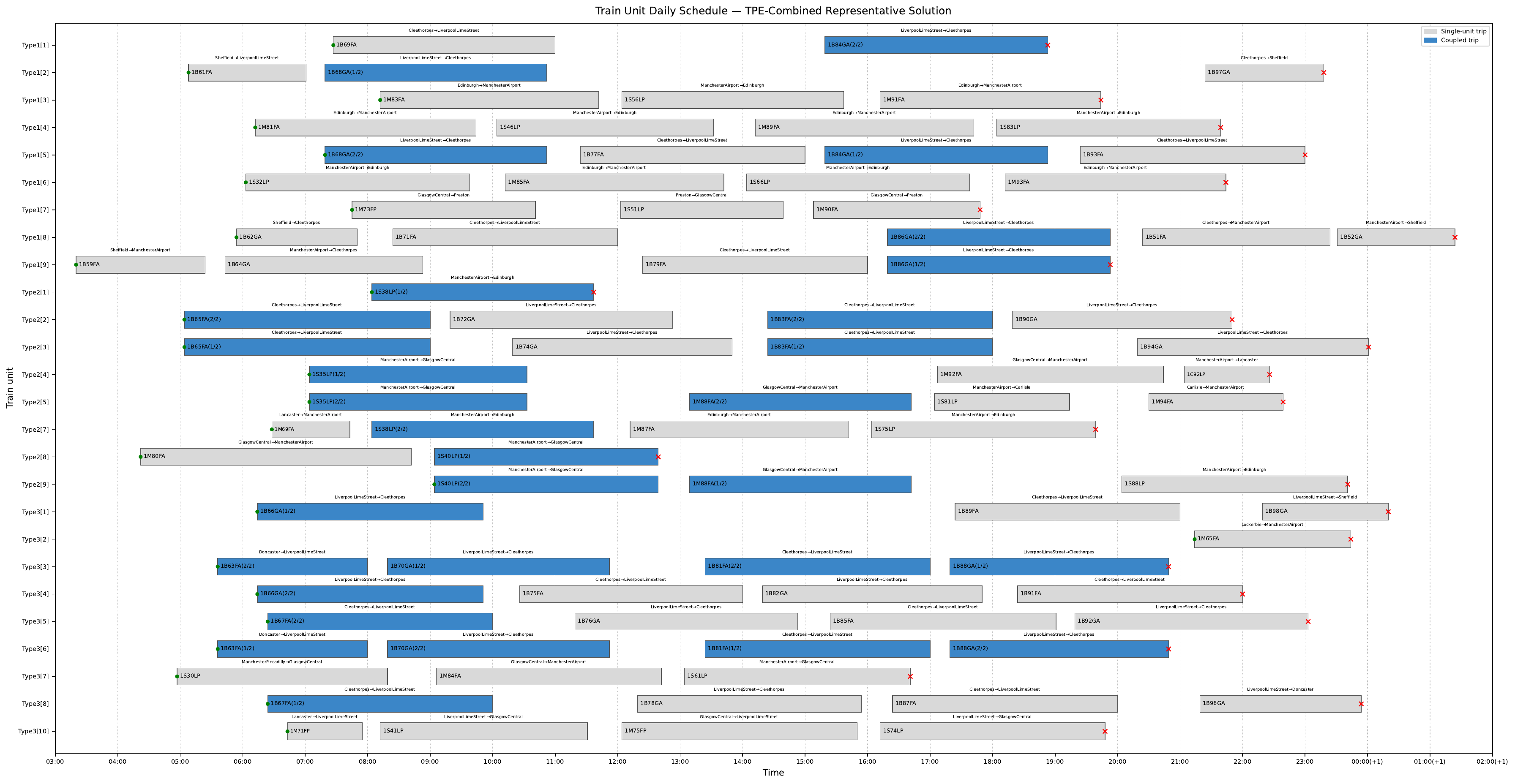}
    \caption{Complete certified unit schedule for the representative TPE-Combined solution. Grey bars denote single-unit trips and blue bars denote coupled trips. For a coupled trip, the label \(i(r/q)\) indicates that the unit occupies position \(r\) in a \(q\)-unit formation.}
    \label{fig:tpe_rep_schedule}
\end{sidewaysfigure}

Figure~\ref{fig:tpe_rep_schedule} shows the complete certified unit schedule for this representative solution. Grey bars denote trips operated by a single unit, whereas blue bars denote trips operated by coupled formations. For coupled trips, the label \(i(r/q)\) indicates that the corresponding unit occupies position \(r\) in a formation of size \(q\), where position 1 is counted from the leading side in the direction of travel. The schedule is certified by TUSOU in the sense that all active continuation instances and all active coupling/decoupling ordering instances among Constraints~(\ref{constraint:couple opposite})--(\ref{constraint:continuation upper}) are satisfied. Thus, the figure represents an executable unit-level schedule rather than only a circulation-feasible assignment.

We next inspect three representative multi-unit station events from the same certified solution. 
These events are active ordering instances selected from the TPE-Combined schedule, rather than additional synthetic test cases. 
They illustrate how the certification procedure checks the original local ordering logic in realistic corridor operations. The complete list of inspected local certificates is reported in Appendix~\ref{app:certificates}.

First, trip \texttt{1S35LP} illustrates same-direction decoupling. 
In the certified solution, \texttt{Type2[5]} continues from \texttt{1S35LP} to \texttt{1M88FA}, which departs at 13:09, whereas \texttt{Type2[4]} continues to \texttt{1M92FA}, which departs later at 17:07. 
For this realised directional case, the earlier-departing unit must be placed on the platform-accessible side. 
The certified order assigns \(\theta=2\) to \texttt{Type2[5]} and \(\theta=1\) to \texttt{Type2[4]}, so the active same-direction decoupling instance of Constraint~(11) is satisfied and the split can be executed without overtaking.

Second, trip \texttt{1M88FA} illustrates a station event where a coupled formation is both created and later split. 
The trip is formed by coupling \texttt{Type2[5]}, arriving from \texttt{1S35LP} at 10:33, with \texttt{Type2[9]}, arriving from \texttt{1S40LP} at 12:39. 
In the certified formation, \texttt{Type2[5]} has \(\theta=2\) and \texttt{Type2[9]} has \(\theta=1\), which satisfies the realised same-direction coupling order. 
The same order also supports the later decoupling: \texttt{Type2[5]} continues to \texttt{1S81LP}, departing at 17:04, while \texttt{Type2[9]} continues to \texttt{1S88LP}, departing at 20:04. 
Thus the unit that must leave first is already on the accessible side, and the active decoupling requirement is also satisfied.

\begin{figure}[htbp]
    \centering
    \begin{subfigure}[t]{0.5\textwidth}
        \centering
        \includegraphics[width=\linewidth]{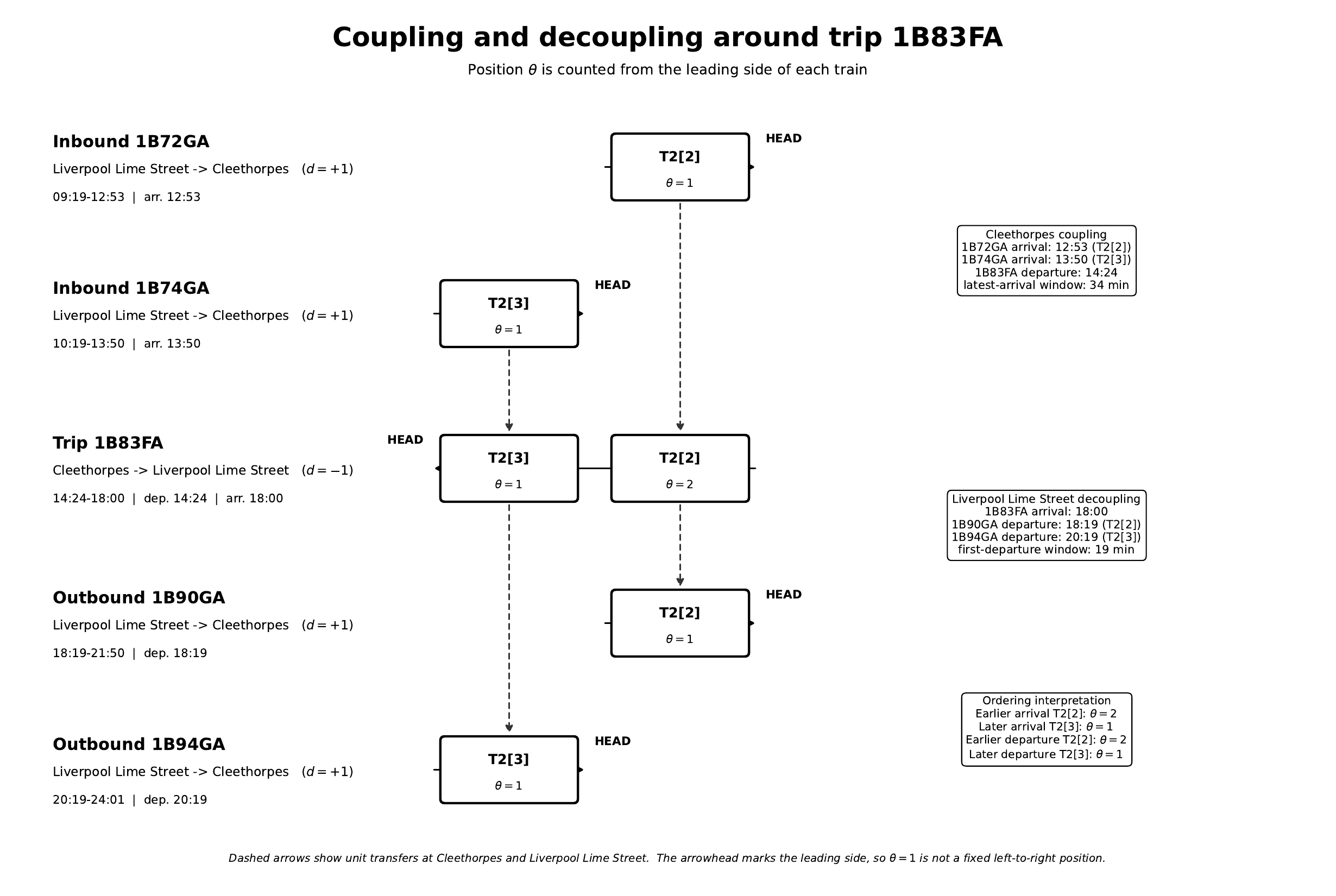}
        \caption{An example of Coupling and decoupling at Cleethorpes}
        \label{fig:tpe1}
    \end{subfigure}\hfill
    \begin{subfigure}[t]{0.5\textwidth}
        \centering
        \includegraphics[width=\linewidth]{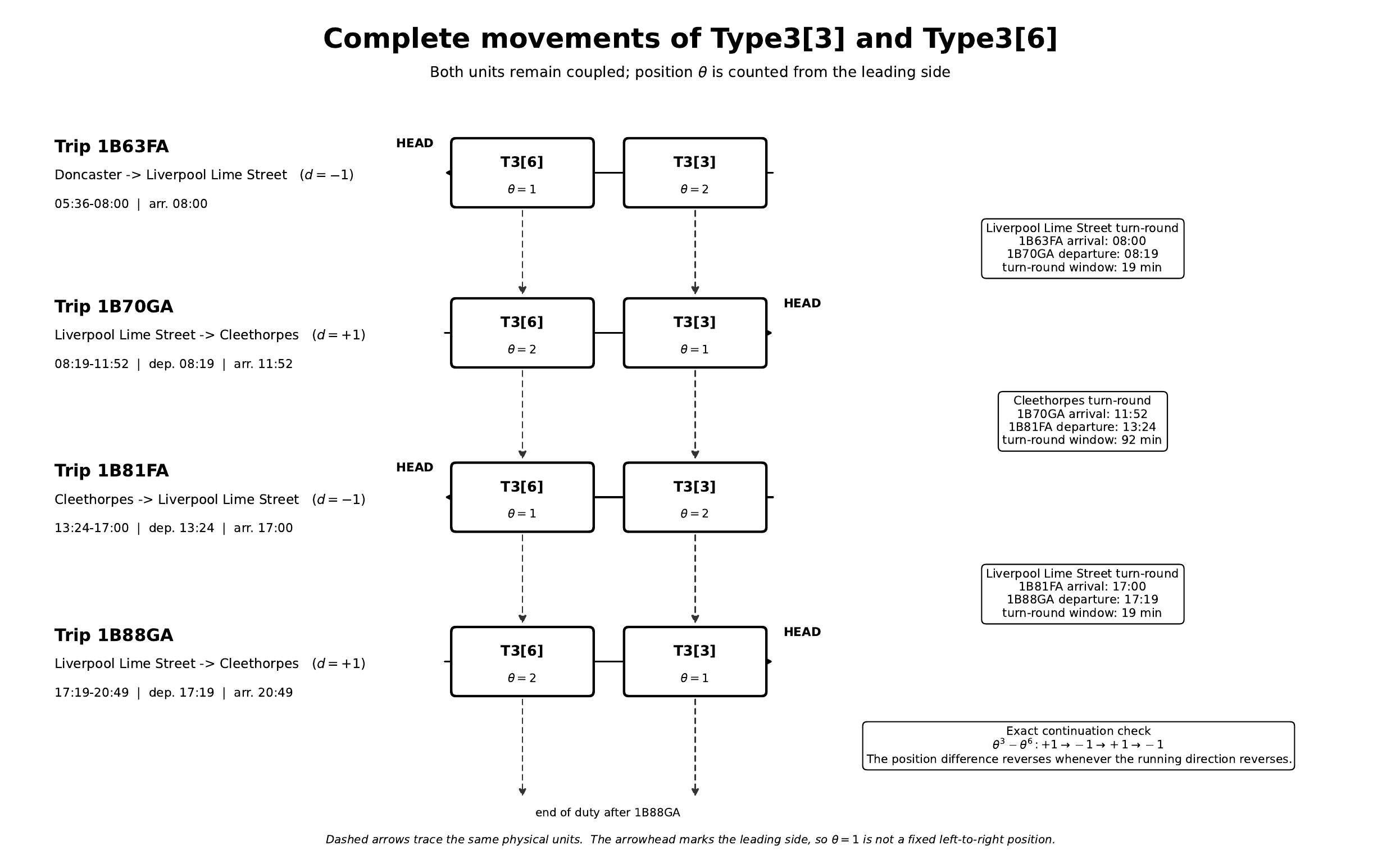}
        \caption{An example of Formation continuation}
        \label{fig:tpe2}
    \end{subfigure}\hfill
\end{figure}

Third, as shown in Figure~\ref{fig:tpe1}, trip \texttt{1B83FA} illustrates a repeated coupling--decoupling pattern. 
The formation is created by coupling \texttt{Type2[2]}, arriving from \texttt{1B72GA} at 12:53, and \texttt{Type2[3]}, arriving from \texttt{1B74GA} at 13:50. 
The certified order places \texttt{Type2[2]} at \(\theta=2\) and \texttt{Type2[3]} at \(\theta=1\), satisfying the realised same-direction coupling condition. 
The same formation then decouples into \texttt{1B90GA}, operated by \texttt{Type2[2]} and departing at 18:19, and \texttt{1B94GA}, operated by \texttt{Type2[3]} and departing at 20:19. 
Because the earlier-departing unit is again positioned on the accessible side, the certified order satisfies the corresponding active decoupling condition.

The same solution also illustrates continuation consistency. For example, Figure~\ref{fig:tpe2} shows \texttt{Type3[3]} and \texttt{Type3[6]} jointly operate the sequence $\texttt{1B63FA} \rightarrow \texttt{1B70GA} \rightarrow \texttt{1B81FA} \rightarrow \texttt{1B88GA}.$ The directions alternate along this chain. On \texttt{1B63FA}, the certified order is $\texttt{Type3[6]} \prec \texttt{Type3[3]}$, whereas on \texttt{1B70GA} the order is reversed: $\texttt{Type3[3]} \prec \texttt{Type3[6]}$. The order reverses again on \texttt{1B81FA} and then again on \texttt{1B88GA}. This is exactly the behaviour imposed by Constraints~(\ref{constraint:continuation lower})-(\ref{constraint:continuation upper}): when two units traverse the same feasible arc and the direction changes, the pairwise position difference changes sign; when the direction is preserved, the pairwise position difference is preserved.

Taken together, the highlighted events show how the certified solution satisfies both parts of the ordering logic. The three coupling and decoupling events above satisfy their active direction-aware precedence requirements, while the multi-trip continuation chain satisfies the pairwise carry-over equalities. A candidate violating any one of these requirements would remain feasible at the network-assignment level but would be station-infeasible under platform-feasible operations.

\subsection{Summary of computational findings}
The reported TUSOU experiments lead to three main conclusions. First, the degree  of local- ordering enforcement is critical. Enforcing the full local-ordering layer at every search node is often computationally expensive, whereas the Shallow strategy, which relies on ordering certification and lazy recovery for exactness, is the most reliable tested TUSOU configuration. It solves all five reported instances to zero reported gap, including the largest NTP2 instance.

Second, the direct full-model Gurobi baseline confirms that the calibrated full ordered-unit formulation is strong on the first four instances, but it does not solve the largest NTP2 instance within the time limit. TUSOU Shallow remains faster than the Gurobi baseline on all five instances and provides explicit diagnostics on non-orderable assignment candidates.

Third, ordering certification is essential. Across the nine terminated TUSOU configurations with candidate diagnostics, 39 of 59 integer assignment-candidate encounters are rejected after full certification of the original local ordering logic. In the Shallow configurations alone, 28 of 38 encounters are rejected. These results show that many assignment-feasible integer candidates are not executable because their selected unit movements admit no continuation-consistent family of trip-wise orders and hence no blockage-free ordering witness.

%% file: paper/7_conclusion.tex
This paper studied the Train Unit Scheduling Problem under platform-feasible operations, where station feasibility depends on the within-formation order of identified train units and where coupling and decoupling must be executable without unpermitted on-platform resequencing. We formulated a single-stage ordered-unit integer linear programming model that jointly determines identified unit movements, turn-round reconfiguration, and within-formation positions. Under the modelled platform restrictions, each feasible integer solution directly yields an executable unit-level schedule, rather than a network-feasible assignment that still requires post hoc ordering repair. The main structural result is an exact fixed-assignment characterisation of blockage-free orderability: active coupling and decoupling decisions induce trip-wise precedence digraphs, while active continuation arcs impose pairwise carry-over consistency. A fixed assignment is orderable if and only if these local precedence digraphs admit a continuation-consistent family of topological orders. This result provides the theoretical basis for ordering certification, which decides whether a candidate assignment \(\bar{x}\) admits at least one blockage-free ordering witness, i.e., whether \(\Omega(\bar{x})\neq\emptyset\).

Building on this characterisation, we developed TUSOU, an exact certification-equipped branch-and-bound-and-cut framework for the ordered-unit formulation. TUSOU searches over network-feasible unit-level assignments, but accepts an integer assignment only when ordering certification returns a feasible ordering witness. Non-orderable assignments are rejected and used to recover missing active ordering instances lazily. This separates exactness from uniform enforcement of all local-ordering rows during LP bounding. Within this framework, ordering certification and lazy recovery provide the exactness mechanism, while activated cycle inequalities serve as auxiliary strengthening for the coupling/decoupling precedence layer.

The computational study on five real-world-derived TransPennine Express instances supports this design while also showing that the calibrated full formulation is strong. The direct Gurobi full-model baseline proves optimality on ASC-Base, TPE-Combined, ASC-Scaled, and NTP, but exceeds the 86,400-second time limit on NTP2. The Shallow TUSOU template is the most reliable tested configuration: it solves all five instances to zero reported gap and is faster than the Gurobi baseline on all five instances. Across the nine terminated TUSOU configurations with candidate diagnostics, ordering certification rejects 39 of 59 integer assignment-candidate encounters; in the Shallow configurations alone, it rejects 28 of 38 encounters. These results indicate that orderability cannot be safely treated as post-processing and should be embedded in the optimisation process.

Several extensions remain open. First, the current platform-feasible operating rules are represented through exogenous permissions and direction-dependent ordering constraints; future work could incorporate richer station-specific infrastructure data, limited shunting resources, dwell-time-dependent resequencing permissions, and more detailed platform-side topology. Second, scalability on larger networks and more heterogeneous fleets may benefit from stronger relaxations, stronger separation routines, alternative extended formulations for selected station events, and hybrid formulations that trade memory for tighter bounds. Third, the ordered-unit framework could be integrated with timetable design, disruption recovery, and robustness analysis, where explicit unit order may be important not only for nominal-day executability but also for recovery options after perturbations.

%% file: paper/12_notation.tex
\begin{table}[H]
\centering
\scriptsize
\setlength{\tabcolsep}{2pt}
\renewcommand{\arraystretch}{0.95}
\caption{Main notation of the formulation with unit ordering.}
\label{table:2}
\begin{tabularx}{\textwidth}{p{0.2\textwidth} X}
\toprule
Symbol & Meaning \\
\midrule
\multicolumn{2}{l}{\textbf{Sets and indices}} \\
$T$ & set of train unit types \\
$F$ & set of coupling-compatible families \\
$F_j$ & set of families feasible for trip $j$ \\
$H_t$ & set of available individual units of type $t$ \\
$H$ & set of all individual train units, $H=\bigcup_{t\in T}H_t$ \\
$\tilde N$ & set of trip nodes \\
$0$ & Artificial source node representing the sign-on of train units \\
$\infty$& Artificial sink node representing the sign-off of train units \\
$N$ & set of all nodes, $N=\tilde N\cup\{0,\infty\}$ \\
$\tilde A$ & set of turn-round arcs between trip nodes \\
$A^{\text{on}},A^{\text{off}},A$ & sign-on arc set, sign-off arc set, and complete arc set \\
$\tilde N^t$ & set of trips that can be operated by type $t$ \\
$\tilde A^t$ & set of feasible turn-round arcs for type $t$ \\
$E^t$ & set of type-$t$ empty-running arcs \\
$A^h$ & set of arcs for individual unit $h$ \\
$A_j^h$ & the set of arcs of unit $h$ that originate from trip node $j$ \\
\midrule
\multicolumn{2}{l}{\textbf{Parameters}} \\
$W_r,\ r=1,2,3$ & weights of the three terms in the objective \\
$d_i$ & direction parameter of trip $i$ \\
$\tau_i^{\text{dep}},\tau_i^{\text{arr}}$ & departure time and arrival time of trip $i$ \\
$\mu_j^t$ & mileage of trip $j$ when operated by a unit of type $t$ (excluding empty-running) \\
$b^t$ & available fleet size of unit type $t$ \\
$\kappa^t$ & capacity of unit type $t$ \\
$r_j$ & passenger demand on trip $j$ \\
$v^f$ & maximum number of units allowed for family $f$ \\
$n^t$ & number of cars in a unit of type $t$ \\
$\bar q_i$ & valid upper bound on the formation size of trip $i$, defined by $\bar q_i=\max_{f\in F_i}v^f$ \\
$M_i$ & trip-specific Big-$M$ constant used in the order-linking constraints \\
$M^{\text{oppC}}_{j,i_1,i_2}$ & Big-$M$ constant for opposite-direction coupling into trip $j$ \\
$M^{\text{sameC}}_{j,i_1,i_2}$ & Big-$M$ constant for same-direction coupling into trip $j$ \\
$M^{\text{oppD}}_{i,j_1,j_2}$ & Big-$M$ constant for opposite-direction decoupling from trip $i$ \\
$M^{\text{sameD}}_{i,j_1,j_2}$ & Big-$M$ constant for same-direction decoupling from trip $i$ \\
$M_{ij}^{\mathrm{cont}}$ & Big-$M$ constant for continuation consistency constraints on arcs between trips $i$ and $j$ \\
\midrule
\multicolumn{2}{l}{\textbf{Decision variables}} \\
$x_a^h\in\{0,1\}$ & equals 1 if unit $h$ traverses arc $a$ \\
$y_j^f\in\{0,1\}$ & equals 1 if trip $j$ is operated by family $f$ \\
$\theta_i^h\in\mathbb{Z}_{\ge 0}$ & position of unit $h$ in the formation of trip $i$, counted from the leading side in the travel direction of trip $i$; $\theta_i^h=0$ if $h$ is not assigned to $i$ \\
$u_i^{h_1h_2}\in\{0,1\}$ & pairwise order indicator; when both units are active on trip $i$, $u_i^{h_1h_2}=1$ iff $\theta_i^{h_1}>\theta_i^{h_2}$ \\
\midrule
\multicolumn{2}{l}{\textbf{Auxiliary notation}} \\
$\delta_+^h(i),\delta_-^h(i)$ & sets of feasible arcs of unit $h$ departing from and arriving at node $i$ \\
$s_i^h$ & assignment indicator; equals 1 if unit $h$ operates trip $i$ \\
$q_i$ & number of units assigned to trip $i$ (formation size on trip $i$) \\
$\Theta_i^h$ & explicit admissible domain of $\theta_i^h$, namely $\{0,1,\ldots,\bar q_i\}$ \\
\bottomrule
\end{tabularx}
\end{table}

%% file: paper/13_proof.tex
\subsection{Detailed proofs of the two-unit feasibility lemmas}
\begin{lemma}[Opposite-direction coupling]\label{lem:equiv-coupling-opposite}
Consider a coupling event that forms trip $j$ from two incoming units associated with trips $i_1$ and $i_2$, and suppose the two units arrive from opposite directions, i.e., $d_{i_1}-d_{i_2}\neq 0$. Then the coupling is platform-feasible under platform-feasible operations if and only if the induced departure formation on trip $j$ satisfies the sign condition $d_j(d_{i_1}-d_{i_2})(\theta_j^{h_1}-\theta_j^{h_2})\ge 0$, which is exactly Constraint~(\ref{constraint:couple opposite}) in its active case.
\end{lemma}

\begin{proof}
Because the units arrive from opposite directions, each unit can enter the platform only from its respective end, so the platform order at the moment of coupling is determined by the arrival directions and cannot be swapped without overtaking. 
The coupling is therefore feasible under platform-feasible operations if and only if the required front--rear order implied by $d_j$ matches this fixed platform order, which is precisely enforced by the sign of $(\theta_j^{h_1}-\theta_j^{h_2})$ in Constraint~(\ref{constraint:couple opposite}). 
If the sign is violated, realising the formation would require the two units to exchange positions on the platform, i.e., overtaking.
\end{proof}

\begin{lemma}[Same-direction coupling]\label{lem:equiv-coupling-same}
Consider a coupling event that forms trip $j$ from two incoming units associated with trips $i_1$ and $i_2$, and suppose the two units arrive from the same direction, i.e., $d_{i_1}+d_{i_2}\neq 0$. Then the coupling is platform-feasible under platform-feasible operations if and only if the departure formation on trip $j$ satisfies the consistency condition $d_j(d_{i_1}+d_{i_2})(\theta_j^{h_1}-\theta_j^{h_2})(\tau_{i_1}^{arr}-\tau_{i_2}^{arr}) \ge 0$, which is exactly Constraint~(\ref{constraint:couple same}) in its active case.
\end{lemma}

\begin{proof}
When both units arrive from the same end, the earlier-arriving unit occupies the inner platform position relative to the common arrival side, and without overtaking this relative platform order cannot be reversed before coupling. 
Therefore, the coupling is feasible under platform-feasible operations if and only if the required order on trip $j$ is consistent with the arrival-time order, with a possible reversal depending on whether $d_j$ matches the common arrival direction, which is captured by the sign of $(\theta_j^{h_1}-\theta_j^{h_2})(\tau_{i_1}^{arr}-\tau_{i_2}^{arr})$ in Constraint~(\ref{constraint:couple same}). 
If violated, the formation would require swapping the two units on the platform, which is impossible without overtaking.
\end{proof}

\begin{lemma}[Opposite-direction decoupling]\label{lem:equiv-decoupling-opposite}
Consider a decoupling event where an arriving formation on trip $i$ decouples into two outgoing units that subsequently operate trips $j_1$ and $j_2$, and suppose the two resulting directions are opposite, i.e., $d_{j_1}-d_{j_2}\neq 0$. Then the decoupling is platform-feasible under platform-feasible operations if and only if the arrival formation on trip $i$ satisfies the sign condition $d_i(d_{j_1}-d_{j_2})(\theta_i^{h_1}-\theta_i^{h_2})\le 0$, which is exactly Constraint~(\ref{constraint:decouple opposite}) in its active case.
\end{lemma}

\begin{proof}
Because the two outgoing trips depart in opposite directions, each unit must leave the platform from a different end, so feasibility under platform-feasible operations requires that each unit is already positioned on the correct side of the decoupling relative to its departure end. 
This is possible if and only if the arrival formation order on trip $i$ is consistent with the two departure directions, which is exactly the sign restriction on $(\theta_i^{h_1}-\theta_i^{h_2})$ enforced by Constraint~(\ref{constraint:decouple opposite}). 
If the sign is violated, at least one unit would have to pass the other to reach its departure end.
\end{proof}

\begin{lemma}[Same-direction decoupling]\label{lem:equiv-decoupling-same}
Consider a decoupling event where an arriving formation on trip $i$ decouples into two outgoing units that subsequently operate trips $j_1$ and $j_2$, and suppose the two resulting directions are the same, i.e., $d_{j_1}+d_{j_2} \neq 0$. Then the decoupling is platform-feasible under platform-feasible operations if and only if the arrival formation on trip $i$ satisfies the temporal-order consistency condition $d_i(d_{j_1}+d_{j_2})(\theta^{h_1}_i-\theta^{h_2}_i)(\tau^{\mathrm{dep}}_{j_1}-\tau^{\mathrm{dep}}_{j_2})\ge 0$, which is exactly Constraint~(\ref{constraint:decouple same}) in its active case.
\end{lemma}

\begin{proof}
When both units depart in the same direction, the unit closer to the departure end can leave first without passing the other, while the unit behind cannot depart earlier without overtaking. 
Hence, the decoupling is feasible under platform-feasible operations if and only if the spatial order in the arriving formation on trip $i$ is consistent with the departure-time order of $j_1$ and $j_2$ up to the direction encoding, which is enforced by the sign of $(\theta_i^{h_1}-\theta_i^{h_2})(\tau^{dep}_{j_1}-\tau^{dep}_{j_2})$ in Constraint~(\ref{constraint:decouple same}). 
If violated, the scheduled first-departing unit would be blocked and would need to pass the other unit on the platform.
\end{proof}

\begin{lemma}[Formation continuation with pairwise order consistency]\label{lem:equiv-continuation}
Consider a feasible arc $a=(i,j)\in \tilde A$ and two units $h_1,h_2$ such that $x_a^{h_1}=x_a^{h_2}=1$. Assume that no off-platform resequencing is allowed on $a$, and that no unit may be inserted between or removed from between $h_1$ and $h_2$ along the continuation. Then the continuation is feasible under platform-feasible operations if and only if $\theta_j^{h_1}-\theta_j^{h_2}=d_i d_j(\theta_i^{h_1}-\theta_i^{h_2})$, which is exactly the active case of Constraints~(\ref{constraint:continuation lower})--(\ref{constraint:continuation upper}).
\end{lemma}

\begin{proof}
Under the stated continuation assumption, the physical spacing and left-to-right order of $h_1$ and $h_2$ remain unchanged between the completion of trip $i$ and the start of trip $j$. If $d_i=d_j$, the leading end is the same on both trips, so the pairwise position difference measured by $\theta$ is preserved. If $d_i=-d_j$, the leading end reverses, so the same physical configuration is read from the opposite end and the pairwise position difference changes sign. Conversely, if the displayed equality holds, the continuation can be realised without overtaking, resequencing, or any insertion/removal between the two units. Hence the condition is necessary and sufficient.
\end{proof}

\subsection{Equivalent fixed-trip encodings and local LP strength}
\label{app:local_encodings}

This appendix formalises the fixed-trip ordering encodings referred to in \S~\ref{sec:Equivalent local encodings and local LP strength}. Fix a trip \(k\) and let \(G_k=(U_k,E_k)\) be the precedence digraph induced by the active coupling/decoupling instances on that trip. Continuation constraints are excluded here because they are cross-trip carry-over conditions rather than local precedence requirements on a single trip. Let \(q_k:=|U_k|\).

Define the set of feasible linear extensions of \(G_k\) by
\[
\Sigma_k :=
\left\{
\sigma:U_k\rightarrow \{1,\ldots,q_k\}\; \middle| \;
\sigma \text{ is a bijection and }
\sigma(h_1)<\sigma(h_2) \text{ whenever } (h_1,h_2)\in E_k
\right\}.
\]
By Lemma~\ref{lem:acyclic_topological}, \(\Sigma_k\neq\emptyset\) if and only if \(G_k\) is acyclic.

For any \(\sigma\in\Sigma_k\), define the corresponding position vector by
\[
\theta^h_k(\sigma):=\sigma(h),\quad h\in U_k,
\qquad
\theta^h_k(\sigma):=0,\quad h\notin U_k.
\]
Let
\[
Q^\theta_k := \operatorname{conv}\{\theta_k(\sigma):\sigma\in\Sigma_k\}.
\]

For any \(\sigma\in\Sigma_k\), define the pairwise precedence vector \(\rho_k(\sigma)\) by
\[
\rho^{k}_{h_1h_2}(\sigma)
:=
\mathbf{1}\{\sigma(h_1)<\sigma(h_2)\},
\qquad
h_1,h_2\in U_k,\; h_1\neq h_2.
\]
Let
\[
Q^\rho_k := \operatorname{conv}\{\rho_k(\sigma):\sigma\in\Sigma_k\}.
\]
Thus, \(Q^\rho_k\) is the convex hull of the pairwise precedence vectors induced by the linear extensions of \(G_k\).

Let $P_{k,\theta}^{\mathrm{BM}}$ denote the projection onto the $\theta_k$-coordinates of the corresponding local Big-$M$ LP
relaxation, including the same active precedence requirements. 
An alternative fixed-trip encoding uses position-assignment variables
\[
\lambda^{h,p}_k\in\{0,1\},
\qquad h\in U_k,\quad p\in\{1,\ldots,q_k\},
\]
where \(\lambda^{h,p}_k=1\) indicates that unit \(h\) occupies position \(p\) on trip \(k\). The assignment constraints are
\begin{align}
\sum_{p=1}^{q_k}\lambda^{h,p}_k &= 1,
&& \forall h\in U_k, \label{eq:app_lambda_unit}\\
\sum_{h\in U_k}\lambda^{h,p}_k &= 1,
&& \forall p\in\{1,\ldots,q_k\}. \label{eq:app_lambda_position}
\end{align}
The position variables are recovered by the projection
\begin{equation}
\theta^h_k
=
\sum_{p=1}^{q_k}p\,\lambda^{h,p}_k,
\quad \forall h\in U_k,
\qquad
\theta^h_k=0,\quad \forall h\notin U_k.
\label{eq:app_lambda_projection}
\end{equation}
To enforce an induced precedence \(h_1\prec_k h_2\), equivalently \((h_1,h_2)\in E_k\), without Big-\(M\), impose the standard prefix inequalities
\begin{equation}
\sum_{p=1}^{r}\lambda^{h_1,p}_k
\ge
\sum_{p=1}^{r}\lambda^{h_2,p}_k,
\qquad
\forall (h_1,h_2)\in E_k,\quad
\forall r\in\{1,\ldots,q_k-1\}.
\label{eq:app_lambda_prefix}
\end{equation}
Let \(Q^\lambda_k\) denote the convex hull of all integer \((\lambda_k,\theta_k)\) satisfying \eqref{eq:app_lambda_unit}--\eqref{eq:app_lambda_prefix} together with \eqref{eq:app_lambda_projection}. For the comparison of the natural LP relaxations, let $\widehat Q_k^\lambda$ denote the polyhedron obtained from equationss~\eqref{eq:app_lambda_unit}--\eqref{eq:app_lambda_prefix} and the projection equations~\eqref{eq:app_lambda_projection} by replacing$\lambda_k^{h,p}\in\{0,1\}$ with $0\leq \lambda_k^{h,p}\leq 1$.

The global optimisation model in the main text is written and solved in the \((x,y,\theta,u)\) variables, with \(s\) used as derived shorthand for trip assignment. The auxiliary variables \(\rho\) and \(\lambda\) are introduced only as local, trip-wise devices for analysing the ordering object and comparing local formulation strength.

\begin{proposition}[Equivalent fixed-trip encodings]
\label{prop:encoding-equivalence}
Fix a trip \(k\) and its induced precedence digraph \(G_k=(U_k,E_k)\). The linear extensions in \(\Sigma_k\), the feasible integer position vectors generating \(Q^\theta_k\), and the feasible integer assignment matrices satisfying \eqref{eq:app_lambda_unit}--\eqref{eq:app_lambda_prefix} are in one-to-one correspondence. Moreover, \(Q^\theta_k\) is exactly the projection of \(Q^\lambda_k\) onto the \(\theta\)-coordinates under \eqref{eq:app_lambda_projection}, and \(Q^\rho_k\) is the convex hull of the induced pairwise precedence vectors.
\end{proposition}

\begin{proof}
Each \(\sigma\in\Sigma_k\) defines a unique assignment matrix \(\lambda_k(\sigma)\) by setting \(\lambda^{h,\sigma(h)}_k=1\) for every \(h\in U_k\) and all other entries to zero. This matrix satisfies the assignment constraints and the prefix inequalities because \(\sigma\) respects every precedence in \(E_k\). The projection \eqref{eq:app_lambda_projection} then gives the corresponding position vector \(\theta_k(\sigma)\).

Conversely, every feasible integer assignment matrix satisfying \eqref{eq:app_lambda_unit}--\eqref{eq:app_lambda_position} defines a bijection from \(U_k\) to \(\{1,\ldots,q_k\}\). The prefix inequalities \eqref{eq:app_lambda_prefix} ensure that this bijection respects every induced precedence in \(E_k\), and hence defines a linear extension in \(\Sigma_k\). Taking convex hulls yields the stated projection and convex-hull relationships.
\end{proof}

\begin{proof}[Proof of Theorem~\ref{thm:fixed_trip_lp_dominance}] Take any $(\lambda_k,\theta_k)\in\widehat Q_k^\lambda$. By the assignment constraints and equations~\eqref{eq:app_lambda_projection}, each $\theta_k^h$ is a convex combination of the positions $\{1,\ldots,q_k\}$, where $q_k=|U_k|$. Hence $1\leq \theta_k^h\leq q_k \qquad \forall h\in U_k$. For every induced precedence $(h_1,h_2)\in E_k$, the prefix inequalities imply $\sum_{p=1}^{r}\lambda_k^{h_1,p}
  \geq
  \sum_{p=1}^{r}\lambda_k^{h_2,p},
  \qquad r=1,\ldots,q_k-1$. Using
\[
  \theta_k^h
  =
  q_k-
  \sum_{r=1}^{q_k-1}\sum_{p=1}^{r}\lambda_k^{h,p},
\]
we obtain $\theta_k^{h_1}\leq\theta_k^{h_2}$. Thus the projected position vector satisfies all corresponding active precedence requirements of the compact formulation.

It remains to show that the relaxed comparison variables can be chosen consistently. For any ordered pair $h_1\neq h_2$, let$\delta=\theta_k^{h_1}-\theta_k^{h_2}$. Because $\delta\in[1-q_k,q_k-1]$ and $M_k\geq q_k$, the relaxed Big-$M$ comparison constraint admits any value
\[
u_k^{h_1h_2}\in
\left[
\max\left\{0,\frac{\delta+1}{M_k}\right\},
\min\left\{1,\frac{\delta+M_k-1}{M_k}\right\}
\right].
\]
This interval is nonempty, so every point in $\operatorname{proj}_{\theta_k}(\widehat Q_k^\lambda)$ belongs to $P_{k,\theta}^{\mathrm{BM}}$.

To establish strictness, consider $U_k=\{1,2\}$, $E_k=\emptyset$, $q_k=M_k=2$. The position-assignment LP satisfies $\theta_k^1+\theta_k^2=3$. However, $\theta_k^1=\theta_k^2=1, u_k^{12}=u_k^{21}=\frac12$ is feasible for the local Big-$M$ LP relaxation of Constraint~(16).
Since $(1,1)$ does not satisfy $\theta_k^1+\theta_k^2=3$, the containment is strict.
\end{proof}

%% file: paper/16_certificates.tex
\begin{table}[hp]
\centering
\scriptsize
\setlength{\tabcolsep}{2pt}
\renewcommand{\arraystretch}{0.95}
\caption{Representative ordering certificates for multi-unit coupling and decoupling events in the TPE-Combined solution. Source and sink sign-on/sign-off arcs are not listed, since the table focuses on trip-to-trip station events governed by the local ordering logic.}
\label{tab:tpe_order_certificates}
\begin{tabular}{l l p{0.37\textwidth} p{0.37\textwidth}}
\hline
Trip & Event & Active movements & Certified ordering feasibility \\
\hline
\texttt{1B65FA} 
& Decoupling 
& \texttt{Type2[2]} continues to \texttt{1B72GA} departing 09:19; \texttt{Type2[3]} continues to \texttt{1B74GA} departing 10:19. 
& Same-direction decoupling. The earlier-departing unit \texttt{Type2[2]} has the larger position, \(\theta=2\), while \texttt{Type2[3]} has \(\theta=1\). Hence the active instance of Constraint~(11) is satisfied. \\

\texttt{1B66GA} 
& Decoupling 
& \texttt{Type3[4]} continues to \texttt{1B75FA} departing 10:26; \texttt{Type3[1]} continues to \texttt{1B89FA} departing 17:24. 
& Same-direction decoupling. The earlier-departing unit \texttt{Type3[4]} has \(\theta=2\), while \texttt{Type3[1]} has \(\theta=1\), so the unit that must leave first is placed on the accessible side. \\

\texttt{1B67FA} 
& Decoupling 
& \texttt{Type3[5]} continues to \texttt{1B76GA} departing 11:19; \texttt{Type3[8]} continues to \texttt{1B78GA} departing 12:19. 
& Same-direction decoupling. The earlier-departing unit \texttt{Type3[5]} has \(\theta=2\), and \texttt{Type3[8]} has \(\theta=1\). The certified order therefore satisfies the temporal-order condition. \\

\texttt{1S35LP} 
& Decoupling 
& \texttt{Type2[5]} continues to \texttt{1M88FA} departing 13:09; \texttt{Type2[4]} continues to \texttt{1M92FA} departing 17:07. 
& Same-direction decoupling. The earlier-departing unit \texttt{Type2[5]} has \(\theta=2\), while \texttt{Type2[4]} has \(\theta=1\). Thus the planned split can be executed without overtaking. \\

\texttt{1B68GA} 
& Decoupling 
& \texttt{Type1[5]} continues to \texttt{1B77FA} departing 11:24; \texttt{Type1[2]} continues to \texttt{1B97GA} departing 21:24. 
& Same-direction decoupling. The earlier-departing unit \texttt{Type1[5]} has \(\theta=2\), while \texttt{Type1[2]} has \(\theta=1\), satisfying the active ordering requirement. \\

\texttt{1M88FA} 
& Coupling 
& \texttt{Type2[5]} arrives from \texttt{1S35LP} at 10:33; \texttt{Type2[9]} arrives from \texttt{1S40LP} at 12:39. 
& Same-direction coupling. The earlier-arriving unit \texttt{Type2[5]} is assigned the larger position, \(\theta=2\), and \texttt{Type2[9]} has \(\theta=1\). The certified departure formation satisfies the active instance of Constraint~(9). \\

\texttt{1M88FA} 
& Decoupling 
& \texttt{Type2[5]} continues to \texttt{1S81LP} departing 17:04; \texttt{Type2[9]} continues to \texttt{1S88LP} departing 20:04. 
& Same-direction decoupling. The earlier-departing unit \texttt{Type2[5]} has \(\theta=2\), while \texttt{Type2[9]} has \(\theta=1\). Therefore the split is blockage-free. \\

\texttt{1B83FA} 
& Coupling 
& \texttt{Type2[2]} arrives from \texttt{1B72GA} at 12:53; \texttt{Type2[3]} arrives from \texttt{1B74GA} at 13:50. 
& Same-direction coupling. The earlier-arriving unit \texttt{Type2[2]} has \(\theta=2\), while \texttt{Type2[3]} has \(\theta=1\). The certified order satisfies the temporal-arrival condition. \\

\texttt{1B83FA} 
& Decoupling 
& \texttt{Type2[2]} continues to \texttt{1B90GA} departing 18:19; \texttt{Type2[3]} continues to \texttt{1B94GA} departing 20:19. 
& Same-direction decoupling. The earlier-departing unit \texttt{Type2[2]} has \(\theta=2\), and the later-departing unit \texttt{Type2[3]} has \(\theta=1\). Hence the active departure-order condition is satisfied. \\

\texttt{1B84GA} 
& Coupling 
& \texttt{Type1[1]} arrives from \texttt{1B69FA} at 11:00; \texttt{Type1[5]} arrives from \texttt{1B77FA} at 15:00. 
& Same-direction coupling. The earlier-arriving unit \texttt{Type1[1]} has \(\theta=2\), while \texttt{Type1[5]} has \(\theta=1\). The resulting coupled formation is therefore platform-feasible. \\

\texttt{1B86GA} 
& Coupling 
& \texttt{Type1[8]} arrives from \texttt{1B71FA} at 12:00; \texttt{Type1[9]} arrives from \texttt{1B79FA} at 16:00. 
& Same-direction coupling. The earlier-arriving unit \texttt{Type1[8]} has \(\theta=2\), while \texttt{Type1[9]} has \(\theta=1\). The certified formation order satisfies the active arrival-order requirement. \\
\hline
\end{tabular}
\end{table}

%% file: main.bib
@article{cacchiani2012nominal,
  title={Nominal and robust train timetabling problems},
  author={Cacchiani, V and Toth, P},
  journal={European Journal of Operational Research},
  volume={219},
  number={3},
  pages={727--737},
  year={2012},
  publisher={Elsevier}
}

@article{lin2014two,
  title={A two-phase approach for real-world train unit scheduling},
  author={Lin, Z and Kwan, RSK},
  journal={Public Transport},
  volume={6},
  pages={35--65},
  year={2014},
  publisher={Springer}
}

@article{cacchiani2010solving,
  title={Solving a real-world train-unit assignment problem},
  author={Cacchiani and Caprara, A and Toth, P},
  journal={Mathematical Programming},
  volume={124},
  pages={207--231},
  year={2010},
  publisher={Springer}
}

@article{schrijver1993minimum,
  title={Minimum circulation of railway stock},
  author={Schrijver, A},
  journal={Cwi Quarterly},
  volume={6},
  number={3},
  pages={205--217},
  year={1993},
  publisher={Citeseer}
}

@article{alfieri2006efficient,
  title={Efficient circulation of railway rolling stock},
  author={Alfieri, A and Groot, R and Kroon, L and Schrijver, A},
  journal={Transportation Science},
  volume={40},
  number={3},
  pages={378--391},
  year={2006},
  publisher={INFORMS}
}

@article{fioole2006rolling,
  title={A rolling stock circulation model for combining and splitting of passenger trains},
  author={Fioole, P-J and Kroon, L and Mar{\'o}ti, G and Schrijver, A},
  journal={European Journal of Operational Research},
  volume={174},
  number={2},
  pages={1281--1297},
  year={2006},
  publisher={Elsevier}
}

@article{peeters2008circulation,
  title={Circulation of railway rolling stock: a branch-and-price approach},
  author={Peeters, M and Kroon, L},
  journal={Computers \& Operations Research},
  volume={35},
  number={2},
  pages={538--556},
  year={2008},
  publisher={Elsevier}
}

@article{freling2005shunting,
  title={Shunting of passenger train units in a railway station},
  author={Freling, R and Lentink, RM and Kroon, LG and Huisman, D},
  journal={Transportation Science},
  volume={39},
  number={2},
  pages={261--272},
  year={2005},
  publisher={INFORMS}
}

@article{kroon2008shunting,
  title={Shunting of passenger train units: an integrated approach},
  author={Kroon, LG and Lentink, RM and Schrijver, A},
  journal={Transportation Science},
  volume={42},
  number={4},
  pages={436--449},
  year={2008},
  publisher={INFORMS}
}

@article{lin2016branch,
  title={A branch-and-price approach for solving the train unit scheduling problem},
  author={Lin, Zhiyuan and Kwan, Raymond SK},
  journal={Transportation Research Part B: Methodological},
  volume={94},
  pages={97--120},
  year={2016},
  publisher={Elsevier}
}

@article{lei2022resolution,
  title={Resolution of coupling order and station level constraints in train unit scheduling},
  author={Lei, Li and Kwan, Raymond SK and Lin, Zhiyuan and Copado-Mendez, Pedro J},
  journal={Public Transport},
  volume={14},
  number={1},
  pages={27--61},
  year={2022},
  publisher={Springer}
}

@article{borndorfer2016integrated,
  title={Integrated optimization of rolling stock rotations for intercity railways},
  author={Bornd{\"o}rfer, Ralf and Reuther, Markus and Schlechte, Thomas and Waas, Kerstin and Weider, Steffen},
  journal={Transportation Science},
  volume={50},
  number={3},
  pages={863--877},
  year={2016},
  publisher={INFORMS}
}

@techreport{grimm2017propagation,
  author = {Grimm, Boris and Borndörfer, Ralf and Reuther, Markus and Schade, Stanley and Schlechte, Thomas},
  title = {A Propagation Approach to Acyclic Rolling Stock Rotation Optimization},
  institution = {Zuse Institute Berlin},
  year = {2017},
  number = {ZIB Report 17-24},
  address = {Berlin, Germany},
  month = {May},
  url = {https://opus4.kobv.de/opus4-zib/frontdoor/index/index/docId/6319}
}

@inproceedings{grimm2019cut,
  title={A cut separation approach for the rolling stock rotation problem with vehicle maintenance},
  author={Grimm, Boris and Bornd{\"o}rfer, Ralf and Reuther, Markus and Schlechte, Thomas},
  booktitle={19th Symposium on Algorithmic Approaches for Transportation Modelling, Optimization, and Systems (ATMOS 2019)},
  pages={1--1},
  year={2019},
  organization={Schloss Dagstuhl--Leibniz-Zentrum f{\"u}r Informatik}
}

@inproceedings{grimm2023assignment,
  title={Assignment based resource constrained path generation for railway rolling stock optimization},
  author={Grimm, Boris and Bornd{\"o}rfer, Ralf and Bushe, Julian},
  booktitle={23rd Symposium on Algorithmic Approaches for Transportation Modelling, Optimization, and Systems (ATMOS 2023)},
  pages={13--1},
  year={2023},
  organization={Schloss Dagstuhl--Leibniz-Zentrum f{\"u}r Informatik}
}

@article{gao2020branch,
  title={A branch-and-price approach for trip sequence planning of high-speed train units},
  author={Gao, Yuan and Schmidt, Marie and Yang, Lixing and Gao, Ziyou},
  journal={Omega},
  volume={92},
  pages={102150},
  year={2020},
  publisher={Elsevier}
}

@article{gao2022weekly,
  title={Weekly rolling stock planning in Chinese high-speed rail networks},
  author={Gao, Yuan and Xia, Jun and D’Ariano, Andrea and Yang, Lixing},
  journal={Transportation Research Part B: Methodological},
  volume={158},
  pages={295--322},
  year={2022},
  publisher={Elsevier}
}

@article{wagenaar2017rolling,
  title={Rolling stock rescheduling in passenger railway transportation using dead-heading trips and adjusted passenger demand},
  author={Wagenaar, Joris and Kroon, Leo and Fragkos, Ioannis},
  journal={Transportation Research Part B: Methodological},
  volume={101},
  pages={140--161},
  year={2017},
  publisher={Elsevier}
}

@article{grimm2025comparison,
  title={A comparison of two models for rolling stock scheduling},
  author={Grimm, Boris and Hoogervorst, Rowan and Bornd{\"o}rfer, Ralf},
  journal={Transportation Science},
  volume={59},
  number={5},
  pages={1101--1129},
  year={2025},
  publisher={INFORMS}
}

@article{cacchiani2013lagrangian,
  title={A Lagrangian heuristic for a train-unit assignment problem},
  author={Cacchiani, Valentina and Caprara, Alberto and Toth, Paolo},
  journal={Discrete Applied Mathematics},
  volume={161},
  number={12},
  pages={1707--1718},
  year={2013},
  publisher={Elsevier}
}

@article{cacchiani2019effective,
  title={An effective peak period heuristic for railway rolling stock planning},
  author={Cacchiani, Valentina and Caprara, Alberto and Toth, Paolo},
  journal={Transportation Science},
  volume={53},
  number={3},
  pages={746--762},
  year={2019},
  publisher={INFORMS}
}

@article{cacchiani2012railway,
  title={Railway rolling stock planning: Robustness against large disruptions},
  author={Cacchiani, Valentina and Caprara, Alberto and Galli, Laura and Kroon, Leo and Mar{\'o}ti, G{\'a}bor and Toth, Paolo},
  journal={Transportation Science},
  volume={46},
  number={2},
  pages={217--232},
  year={2012},
  publisher={INFORMS}
}

@article{borndorfer2021deutsche,
  title={Deutsche bahn schedules train rotations using hypergraph optimization},
  author={Bornd{\"o}rfer, Ralf and E{\ss}er, Thomas and Frankenberger, Patrick and Huck, Andreas and Jobmann, Christoph and Krostitz, Boris and Kuchenbecker, Karsten and Mohrhagen, Kai and Nagl, Philipp and Peterson, Michael and others},
  journal={INFORMS Journal on Applied Analytics},
  volume={51},
  number={1},
  pages={42--62},
  year={2021},
  publisher={INFORMS}
}

@article{cacchiani2013integer,
  title={On integer polytopes with few nonzero vertices},
  author={Cacchiani, Valentina and Caprara, Alberto and Mar{\'o}ti, G{\'a}bor and Toth, Paolo},
  journal={Operations Research Letters},
  volume={41},
  number={1},
  pages={74--77},
  year={2013},
  publisher={Elsevier}
}

@article{cadarso2011robust,
  title={Robust rolling stock in rapid transit networks},
  author={Cadarso, Luis and Mar{\'\i}n, {\'A}ngel},
  journal={Computers \& Operations Research},
  volume={38},
  number={8},
  pages={1131--1142},
  year={2011},
  publisher={Elsevier}
}

@article{cadarso2014improving,
  title={Improving robustness of rolling stock circulations in rapid transit networks},
  author={Cadarso, Luis and Mar{\'\i}n, {\'A}ngel},
  journal={Computers \& Operations Research},
  volume={51},
  pages={146--159},
  year={2014},
  publisher={Elsevier}
}

@article{lin2017train,
  title={Train unit scheduling guided by historic capacity provisions and passenger count surveys},
  author={Lin, Zhiyuan and Barrena, Eva and Kwan, Raymond SK},
  journal={Public Transport},
  volume={9},
  number={1},
  pages={137--154},
  year={2017},
  publisher={Springer}
}

@article{lin2016local,
  title={Local convex hulls for a special class of integer multicommodity flow problems},
  author={Lin, Zhiyuan and Kwan, Raymond SK},
  journal={Computational Optimization and Applications},
  volume={64},
  number={3},
  pages={881--919},
  year={2016},
  publisher={Springer}
}

@article{lin2020avoiding,
  title={Avoiding unnecessary demerging and remerging of multi-commodity integer flows},
  author={Lin, Zhiyuan and Kwan, Raymond SK},
  journal={Networks},
  volume={76},
  number={2},
  pages={206--231},
  year={2020},
  publisher={Wiley Online Library}
}

@article{lusby2017branch,
  title={A branch-and-price algorithm for railway rolling stock rescheduling},
  author={Lusby, Richard M and Haahr, J{\o}rgen Thorlund and Larsen, Jesper and Pisinger, David},
  journal={Transportation Research Part B: Methodological},
  volume={99},
  pages={228--250},
  year={2017},
  publisher={Elsevier}
}

@article{nielsen2012rolling,
  title={A rolling horizon approach for disruption management of railway rolling stock},
  author={Nielsen, Lars Kj{\ae}r and Kroon, Leo and Mar{\'o}ti, G{\'a}bor},
  journal={European Journal of Operational Research},
  volume={220},
  number={2},
  pages={496--509},
  year={2012},
  publisher={Elsevier}
}

@article{wagenaar2017maintenance,
  title={Maintenance appointments in railway rolling stock rescheduling},
  author={Wagenaar, Joris C and Kroon, Leo G and Schmidt, Marie},
  journal={Transportation Science},
  volume={51},
  number={4},
  pages={1138--1160},
  year={2017},
  publisher={INFORMS}
}

@book{trotter2002combinatorics,
  title={Combinatorics and partially ordered sets},
  author={Trotter, William T},
  year={2002},
  publisher={Johns Hopkins University Press}
}
